\documentclass[mod]{jour}

\usepackage[english]{babel}
\usepackage{epsfig,psfrag}
\usepackage{amsfonts}
\usepackage{amsmath}
\usepackage{amssymb}
\usepackage{multicol}
\numberwithin{equation}{section} 
\numberwithin{theorem}{section} 
\numberwithin{proposition}{section} 
\numberwithin{lemma}{section} 
\numberwithin{corollary}{section} 
\numberwithin{definition}{section} 
\numberwithin{conjecture}{section} 
\numberwithin{remark}{section}

\newcommand{\R}{{\mathbb{R}}}
\newcommand{\Z}{{\mathbb{Z}}}
\newcommand{\N}{{\mathbb{N}}}

\newcommand{\Tplato}{{\mathfrak{T}}}
\newcommand{\Cplato}{{\mathfrak{C}}}
\newcommand{\Oplato}{{\mathfrak{O}}}
\newcommand{\Dplato}{{\mathfrak{D}}}
\newcommand{\Iplato}{{\mathfrak{I}}}

\newcommand{\K}{{\mathcal{K}}}
\newcommand{\Kfour}{{\mathcal{K}_4}}
\newcommand{\KPi}{{\mathcal{K}^P_i}}
\newcommand{\KPuno}{{\mathcal{K}^P_1}}
\newcommand{\KPdue}{{\mathcal{K}^P_2}}
\newcommand{\KPtre}{{\mathcal{K}^P_3}}
\newcommand{\KPj}{{\mathcal{K}^P_j}}
\newcommand{\KPh}{{\mathcal{K}^P_h}}
\newcommand{\Kclo}{\overline{\mathcal{K}}}

\newcommand{\er}{\mathsf{e}_r}
\newcommand{\es}{\mathsf{e}_s}
\newcommand{\eM}{\mathsf{e}_M}
\newcommand{\eV}{\mathsf{e}_V}
\newcommand{\euno}{\mathsf{e}_1}
\newcommand{\edue}{\mathsf{e}_2}
\newcommand{\etre}{\mathsf{e}_3}
\newcommand{\ej}{\mathsf{e}_j}
\newcommand{\ealpha}{\mathsf{e}_\alpha}
\newcommand{\ebeta}{\mathsf{e}_\beta}
\newcommand{\eort}{\mathsf{e}_\perp}
\newcommand{\oldnu}{\mathsf{e}_q}
\newcommand{\doldnu}{\dot{\mathsf{e}}_q}
\newcommand{\olda}{b}

\newcommand{\A}{{\mathcal{A}}}

\newcommand{\fundD}{{\rm D}}
\newcommand{\fundS}{{\rm S}}

\newcommand{\RR}{{\mathcal{R}}}

\def\accauno{H^1_T(\R,\mathcal{X})}
\def\accaunoR3{H^1_T(\R,\mathbb{R}^3)}
\def\lama{\Lambda^{{\sf (a)}}}
\def\lamTG{\Lambda_G}
\def\lamTP{\Lambda^P}
\def\lamoa{\Lambda_0^{\sf (a)}}
\def\lamoasim{\Lambda_0^{\sf (a)}\hskip -0.15cm / \hskip -0.15cm\sim}

\def\faceref{F}

\def\asseuno{\xi_1}
\def\assedue{\xi_2}
\def\assetre{\xi_3}
\def\assej{\xi_j}

\def\condA{{\sf (a)}}
\def\condB{{\sf (b)}}
\def\condC{{\sf (c)}}

\def\unitsphere{\mathsf{S}^2}

\def\Cquad{C_\diamond}

\newcommand{\partgen}{\mathsf{P}_1}
\newcommand{\particle}{\mathsf{P}}

\def\genp{u_1}
\def\hatgenp{\hat{u}_1}

\def\checkgenp{\check{u}_1}
\def\dgenp{\dot{u}_1}
\def\genpv{v_1}

\def\usn{{\rm u}^{(\sigma,n)}}
\def\usuno{{\rm u}^{(\sigma,1)}}
\def\usngen{{\rm u}^{(\sigma,n)}_1}
\def\vnun{{\rm v}^{(\nu,n)}}
\def\vnungen{{\rm v}_1^{(\nu,n)}}
\def\dvnungen{\dot{\rm v}_1^{(\nu,n)}}
\def\vnuuno{{\rm v}^{(\nu,1)}}
\def\vnuunogen{{\rm v}_1^{(\nu,1)}}

\def\usunugen{{\rm u}^{(\sigma_u,n_u)}_1}
\def\usunusgen{{\rm u}^{(\sigma_u,n_u,s)}_1}
\def\usunugenzero{{\rm u}^{(\sigma_u,n_u,0)}_1}
\def\usunugenuno{{\rm u}^{(\sigma_u,n_u,1)}_1}

\def\usunuuno{{\rm u}^{(\sigma_u,n_u,1)}}

\def\typloop{u}
\def\tloophat{\hat{u}}
\def\sloop{v}
\def\sloopgen{v_1}

\def\npm{\mathsf{n}^\pm}
\def\npiu{\mathsf{n}^+}
\def\nmeno{\mathsf{n}^-}

\def\minloop{u_*}
\def\hatminloop{\hat{u}_*}
\def\minloopgen{u_{*,1}}
\def\hatminloopgen{\hat{u}_{*,1}}
\def\dhatminloopgen{\dot{\hat{u}}_{*,1}}
\def\dminloopgen{\dot{u}_{*,1}}

\def\minloopPi{u_*^{P,i}}
\def\minloopPj{u_*^{P,j}}
\def\motionloopPi{v_*^{P,i}}
\def\motionloopPigen{v_{*,1}^{P,i}}

\begin{document}
\title{\bf Platonic Polyhedra, Topological Constraints 
and Periodic Solutions of the Classical $N$-Body Problem}
\titlerunning{Platonic polyhedra and periodic orbits of the $N$--body problem}
\author{G. Fusco\inst{1} \and G.~F. Gronchi\inst{2}\and P. Negrini\inst{3}}
\institute{Dip. di Matematica Pura ed Applicata, Universit\`a di L'Aquila
\and
Dip. di Matematica, Universit\`a di Pisa
\and
Dip. di Matematica, Universit\`a di Roma {\em`La Sapienza'}
}
\date{Received: date / Revised version: date}
\maketitle

\begin{multicols}{2}{
\quad\\
\quad\\
\quad\\
\quad\\
\quad\\
\quad\\
\quad\\
\quad\\
\quad\\
\quad\\
\quad\\

\noindent{\em Propositum est mihi, Lector, hoc libello demonstrare, quod
Creator Optimus Maximus, in creatione mundi huius mobilis, et dispositione
coelorum, ad illa quinque regularia corpora, inde a Pythagora et Platone, ad
nos utque, celebratissima respexerit, atque ad illorum naturam coelorum
numerum, proportiones, et motuum rationem accommodaverit.}
\ (J. Kepler, {\em Myst. Cosm.} \cite{kepler})}
\end{multicols}
\begin{abstract}
We prove the existence of a number of smooth periodic motions
$\minloop$ of the classical Newtonian $N$--body problem which, up to a
relabeling of the $N$ particles, are invariant under the rotation
group ${\cal R}$ of one of the five Platonic polyhedra.  The number
$N$ coincides with the order $|{\cal R}|$ of ${\cal R}$ and the
particles have all the same mass. Our approach is variational and
$\minloop$ is a minimizer of the Lagrangean action $\A$ on a suitable
subset $\K$ of the $H^1$ $T$--periodic maps $u:\R\to\R^{3N}$. The set
$\K$ is a cone and is determined by imposing to $u$ both topological
and symmetry constraints which are defined in terms of the rotation
group ${\cal R}$. There exist infinitely many such cones $\K$, all
with the property that $\A|_{\K}$ is coercive. For a certain number of
them, using level estimates and local deformations, we show that
minimizers are free of collisions and therefore classical solutions of
the $N$--body problem with a rich geometric--kinematic structure.
\end{abstract}


\newpage     

\centerline{\bf\Large List of symbols}
\begin{table}[th]
\label{tab:symbols}
\begin{center}
\begin{tabular}{l|r}
{\bf symbol}    &{\bf meaning} \\
\hline &\\ 
$\Tplato$, $\Cplato$, $\Oplato$, $\Dplato$, $\Iplato$
                &the five Platonic polyhedra \\
%
${\cal T}, {\cal O}, {\cal I}$ &symmetry groups of rotations
                                of the Platonic polyhedra\\ 
${\cal Q}_{\cal T}, {\cal Q}_{\cal O}, {\cal Q}_{\cal I}$ 
&Archimedean polyhedra, Section~\ref{subsec:chartop}\\
$\mathcal{X}$        &configuration space\\
$\mathfrak{S}$  &loops with collisions\\
$\Gamma$ 
&set of the axes of rotation of the group ${\cal R}\in \{ {\cal T}, {\cal O},
{\cal I}\}$\\
$\tilde{\cal R}$  &reflection group associated to ${\cal R}$\\
$\typloop$      &generic element of $H^1_T(\R,R^{3N})$ and\\
$\partgen$      &generating particle\\
$\tau_1$        &trajectory of the generating particle\\
$\genp$         &generating loop for $\typloop$\\
$\checkgenp$, $\hatgenp$    &variations of $\genp$\\
$\minloop$,  $\minloopgen$     &minimizer and related generating loop\\
%
%
$\lama$, $\lamoa$        &loop spaces\\
$\lamoasim$     &space of equivalence classes of loops in $\lamoa$\\
$\K$, $\KPi$, $\tilde{\K}^\nu$           &cones of loops\\
%
%
%
$\unitsphere$   &unit sphere in $\R^3$\\
$\ej$             &direction vector of the axis $\assej, j=1,2,3$\\
$\N$            &natural numbers excluding 0\\
$\usngen$        &map associated to a periodic sequence $(\sigma,n)$ of
                 domains $D_k$, see (\ref{genpsigman}) \\
$\vnungen$       &map associated to a periodic sequence $(\nu,n)$ of 
                 vertexes of ${\cal Q}_{\cal R}$, see (\ref{vnun}) \\
$\fundD$        &fundamental domain   \\
$\fundS_1$, $\fundS_2$, $\fundS_3$ &faces of $\fundD$ \\
$D_k$          &domains obtained by applying the elements of 
               $\tilde{\cal R}$ to $\fundD$ \\
$S_k$, $S^k$   &faces of $D_k$\\
$R_{S}$        &reflection with respect to the plane of the face $S$\\
$R \mapsto R^S$  &bijection defined in ${\cal R}\setminus\{I\}$ for a given $S$, cfr. (\ref{eq:RSRx})
\end{tabular} 
\end{center}
\end{table}
\newpage

%

\section{Introduction}
\label{sec:intro}

In the last few years many interesting periodic motions of the
classical Newtonian $N$--body problem have been discovered as
minimizers of the {\em Action functional}
\[
\A : \lamTG \to \R\cup \{+\infty\}\,,
\]
\begin{equation}
\A(u) = \int_0^T\Bigl(\frac{1}{2}\sum_{h=1}^Nm_h|\dot{u}_h|^2 +
\sum_{1\leq h<k\leq N}\frac{m_h\,m_k}{|u_h-u_k|} \Bigr)\,dt
\label{action_functional}
\end{equation}
on loop spaces $\lamTG \subset \accauno$, of $T$--periodic motions,
equivariant with respect to the action of a suitably chosen group $G$
\cite{FT04},\cite{chenciner03}.  We denote by $u=(u_1,\dots,u_N):\R \to
\mathcal{X}$ a typical element of $\lamTG$, by $u_h:\R \to \R^3$ the motion
of the mass $m_h$
and by $\mathcal{X}\subset \R^{3N}$ the configuration space
\[
\mathcal{X} = \Bigl\{x=(x_1,\ldots, x_N)\in\R^{3N} : \sum_{h=1}^N m_h x_h = 0
\Bigr\}\ .
\]
$\accauno$ denotes the Sobolev space of $L^2$ $T$--periodic maps
$u:\R\to\mathcal{X}$ with $L^2$ first derivative.

The interest in this classical problem was revived by the discovery of the now
famous Eight \cite{moore93}, \cite{CM00}: the rather surprising fact that
three equal masses can move periodically one after the other with a time shift
of $T/3$ on the same fixed planar trajectory which has the shape of a
symmetric eight.  Another remarkable motion is the Hip--Hop \cite{CV00},
\cite{DTW83} where four equal masses, at intervals of $T/2$, coincide
alternatively with the vertexes of two tetrahedra, each one symmetric of the
other with respect to the center of mass.

The lack of coercivity of the action functional on the whole set
$\accauno$ of $T$--periodic motions and the fact that, on the basis of
Sundman's estimates \cite{wintner}, \cite{saari}, collisions give a bounded
contribution to the action are the main mathematical obstructions in
the search for $T$--periodic motions of the $N$--body problem as
minimizers of the action. One of the main results of the research
effort developed in the last ten years is a quite
systematic way to deal with these obstructions. A basic observation
concerning the problem of coercivity is that the action functional
(\ref{action_functional}) is coercive when restricted to a loop space
$\lamTG$ of motions that possess suitable symmetries. The idea to
impose symmetries to obtain coercivity was introduced in \cite{DGM}
and used in \cite{cotizelati}, and it is
considered in a general abstract context in \cite{FT04}, where also a
necessary and sufficient condition for coercivity is given.

The original motivation for our work was aesthetical: we wondered about the
existence of new periodic motions which we could compare in perfection and
beauty with the Eight, the Hip--Hop and the other interesting motions that
have been recently discovered, see \cite{simo}, \cite{chen}, \cite{terrvent}.

\begin{figure}
\psfragscanon
\psfrag{O}{$O$}
\psfrag{M}{$M$}
\psfrag{V}{$V$}
\psfrag{rM}{$\eM$}
\psfrag{rV}{$\eV$}
\psfrag{ra}{$\ealpha$}
\psfrag{rb}{$\ebeta$}
\psfrag{xi1}{$\asseuno$}
\psfrag{xi2}{$\assedue$}
\psfrag{xi3}{$\assetre$}
\centerline{\epsfig{figure=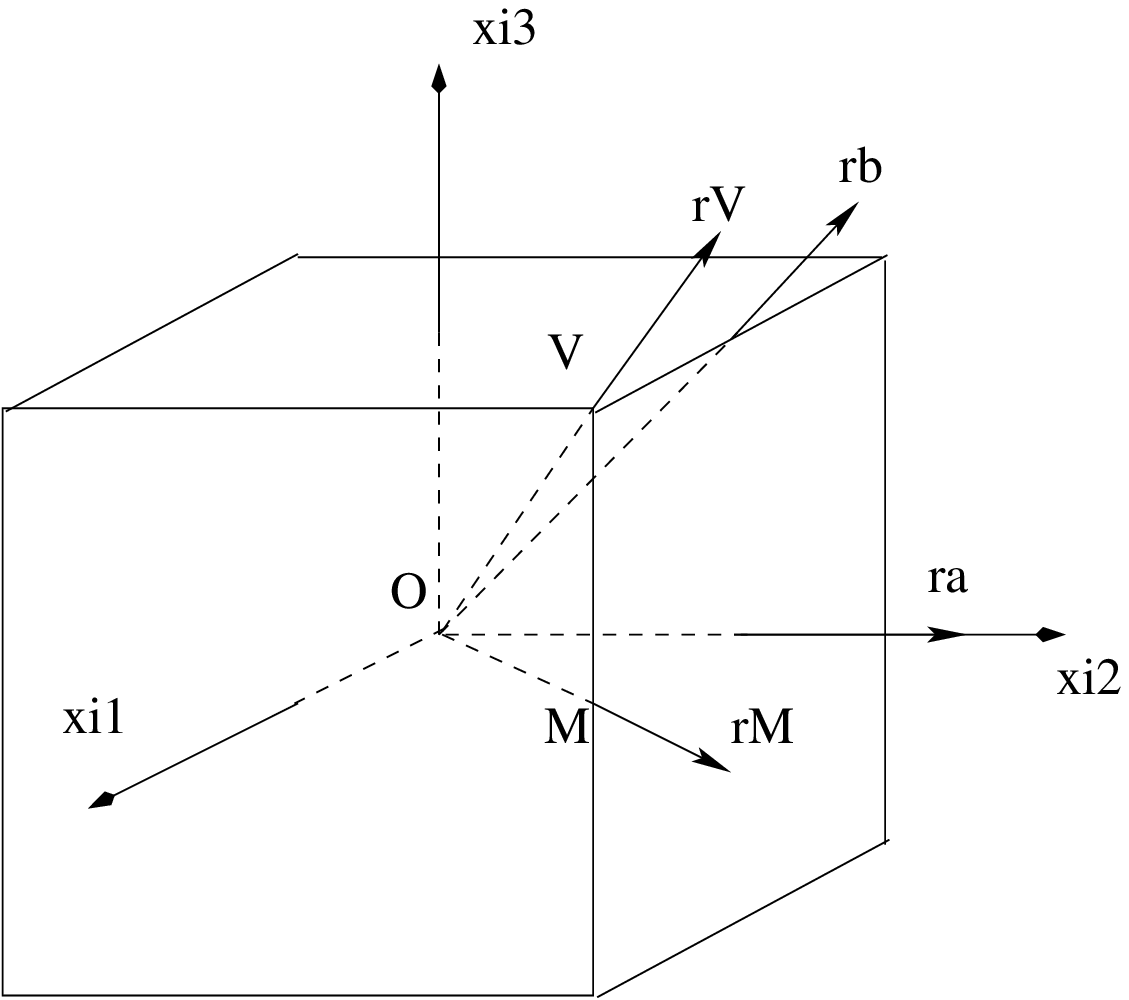,width=6cm}}
\psfragscanoff
\caption{The unit vectors $\eM, \eV, \ealpha, \ebeta$ for $P=\Cplato$.}
\label{fig:reference}
\end{figure}

Let $P$ be one of the five Platonic polyhedra, that is $P\in\{\Tplato$,
$\Cplato$, $\Oplato$, $\Dplato$, $\Iplato\}$ where $\Tplato, \Cplato, \Oplato,
\Dplato, \Iplato$ stand for {\em Tetrahedron}, {\em Cube}, {\em Octahedron},
{\em Dodecahedron} and {\em Icosahedron}.  Let ${\cal T}, {\cal O}, {\cal I}$
be the groups of rotations of $\Tplato$, of $\Cplato$ and $\Oplato$, of
$\Dplato$ and $\Iplato$ respectively, and denote by ${\cal R}\in\{{\cal T},
{\cal O}, {\cal I}\}$ the group of rotations of $P$.  Let $K = 4,6,8,12,20$ be
the number of the faces of $P$ and $H=3,4,3,5,3$ the number of the vertexes of
each face of $P$.  Let $\faceref$ be one of the faces of $P$, $L$ one of the
sides of $\faceref$ and $M$ the middle point of $L$. Consider a right--handed
orthogonal frame $O \asseuno \assedue \assetre$, with the origin in the center
of $P$, the axis $\asseuno$ oriented from $O$ to the center of $\faceref$
and the axes $\assedue, \assetre$ such that $M$ lies in the half--plane
$\{\assetre = 0, \assedue>0\}$. We let $V$ be the vertex of $L$ in the
half--space $\assetre>0$ (see Figure~\ref{fig:reference}). We denote by $\ej$
the direction vector of $\assej, j=1,2,3$.  The question that was at the
origin of the present work is the following:

\smallbreak\noindent ({\cal Q})\quad {\em Does it exist a $T$--periodic motion
of the classical Newtonian $N$--body problem with $N=HK=12,24,60$ equal masses
that satisfies conditions \condA, \condB, \condC\ below?}

\begin{itemize}
\item[\condA] If $\genp :\R\to\R^3$,
  $\genp(t+T)=\genp(t)\ \forall t\in\R$, is the motion of one of the $N$
  particles, called the {\em generating particle}, the motion of the
  $N$ particles is determined by a bijection:
\begin{equation}
\{2,\dots,N\}\ni j\to R_j\in{\cal
  R}\setminus \{I\}: u_j = R_j \genp\ .
\label{partgen}
\end{equation} 

\item[\condB] Associated to each face of $P$ there are $H$
  particles that move one after the other on the same trajectory with a time
  shift of $T/H$.

\item[\condC] The motion of the generating particle satisfies
\[
\genp(t) = S_3 \genp(-t)
\quad\quad \forall t\in\R\,,
\]
where $S_3$ is the reflection with respect to the plane $\assetre = 0$.
\end{itemize}
We remark that \condA, \condB, \condC\ imply in particular that the trajectory
\[
\tau_1=\{x\in\R^3:\exists\, t \mbox{ with } x=\genp(t) \}
\]
has all the symmetries of $\faceref$.

\begin{figure}[h]
\centerline{$a)$\psfig{figure=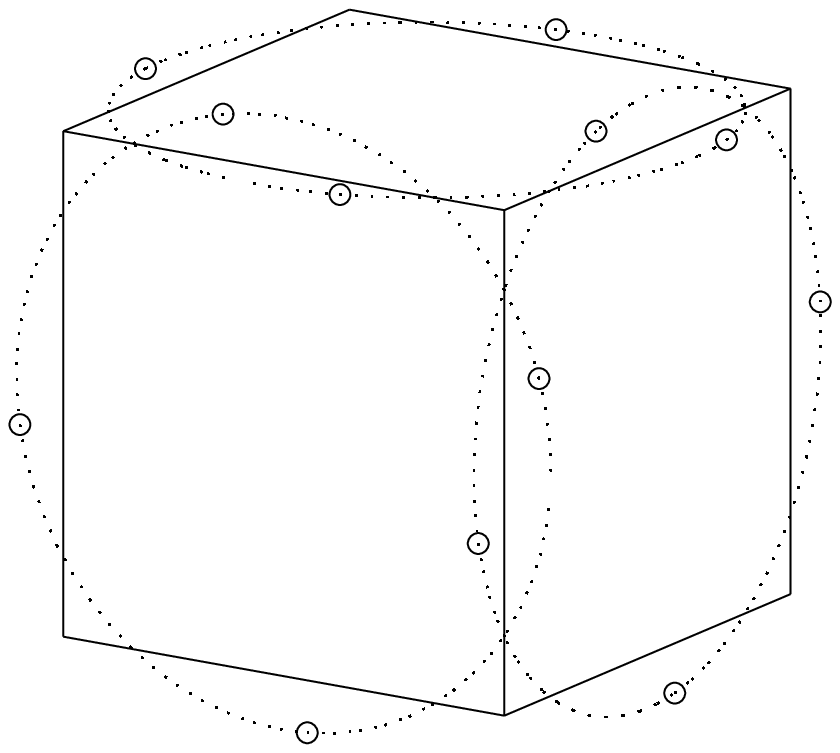,width=4.5cm}
\hskip 1.5cm
$b)$\epsfig{figure=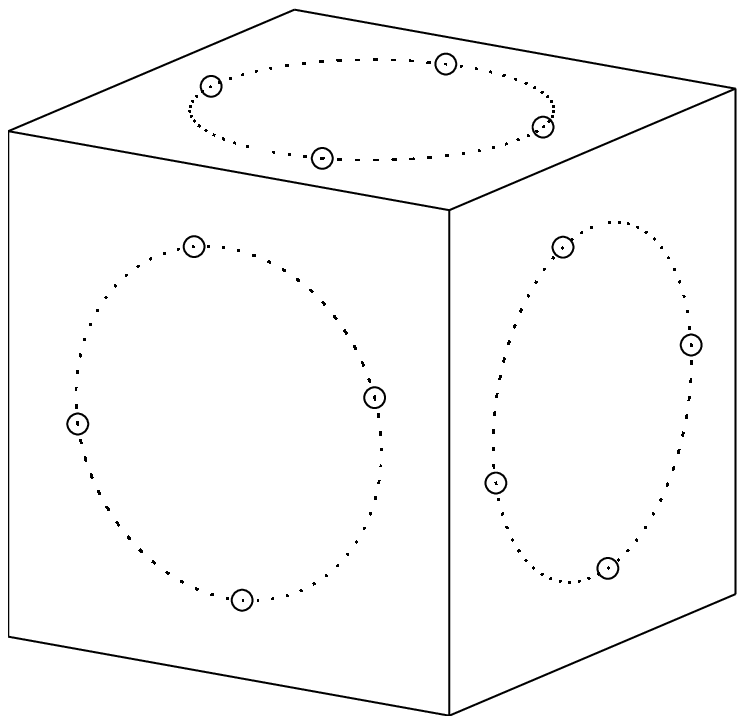,width=4cm}}
\caption{Loops in the space $\lamTP$ for $P=\Cplato$.}
\label{fig:motions_cube}
\end{figure}

We denote by $\lamTP\subset \accauno$ the subset of $T$--periodic maps that
satisfy conditions \condA, \condB, \condC.  Figure~\ref{fig:motions_cube}
visualizes possible structures of motions that satisfy \condA, \condB, \condC\
for the case $P=\Cplato$. One checks immediately that the action functional
(with $m_j=1,\, j=1\ldots N$)
\begin{equation}
\A(u) = \frac{N}{2}\int_0^T \biggl(|\dgenp|^2 + \sum_{R_j\in{\cal
R}\setminus\{I\}}\frac{1}{|(R_j-I)\genp|}\biggr)\,dt
\label{action_formula}
\end{equation}
is not coercive on the loop space $\lamTP$.  Indeed, if we set $\genp^\lambda
= \genp + \lambda \euno$, $\lambda>0$, $\euno=(1,0,0)$, we have
\small
\[
\A (u) - \A (u^\lambda) = \frac{N}{2}\int_0^T \sum_{R_j\in{\cal
R}\setminus\{I\}} \left(\frac{1}{|(R_j-I)\genp|} - \frac{1}{|(R_j-I)\genp +
\lambda(R_j-I)\euno|} \right)\,dt
\]
\normalsize
and therefore $\A(u) - \A(u^\lambda) > 0$ for $\lambda$ large enough.  This is
also a consequence of the abstract coercivity condition formulated in
\cite{FT04}.  In \cite{FT04} a motion $u:\R\to\mathcal{X}$ is said to
be {\em equivariant} with respect to the action of a finite group $G$ if
\begin{equation}
\rho(g)\;u_{\sigma(g^{-1})(i)} (t) = u_i(\tau(g)(t))\qquad \forall\,
g\in G,\ \forall\, t\in\R,\ \forall\, i\in\{1\dots N\}\ .
\label{equivariance}
\end{equation}
Here $\rho:G\to O(3)$, $\tau:G\to O(2)$ are orthogonal representations of $G$
and $\sigma:G\to S_N$ is a homomorphism of $G$ to $S_N$, the group of
permutations of $\{1\dots N\}$.\footnote{If the masses of the $N$ particles
are not all equal it is required that $\sigma(g)(i) = j \Rightarrow
m_i=m_j$.}  If we let $\lamTG\subset \accauno$ be the set of loops that
fulfill the equivariance condition (\ref{equivariance}), then the coercivity
condition formulated and proven in \cite{FT04} is the following:

\begin{theorem}
The action functional $\A$ is coercive on $\lamTG$ if and only
if
\begin{equation}
\mathcal{X}^G = \{0\}
\label{coercond}
\end{equation}
where $\mathcal{X}^G\subset \R^{3N}$ is the subset of the
configuration space invariant under $G$, that is
\[
\mathcal{X}^G = \bigl\{ x\in\mathcal{X}:\rho(g)x_{\sigma(g^{-1})(i)} = x_i,
\forall\, g\in G,\ \forall\, i\in\{1\dots N\}\bigr\}\ .
\]
\label{coerteor}
\end{theorem}

For the loop space $\lamTG=\lamTP$, defined by \condA, \condB, \condC,
condition (\ref{coercond}) is not satisfied. Indeed the nonzero vector
\[
x=(x_1,\ldots,x_N), \qquad x_1=\euno\,, x_j = R_j \euno,\ j=2,\ldots,N
\]
belongs to $\mathcal{X}^G$.

The fact that for $\lamTG = \lamTP$ condition (\ref{coercond}) is violated
does not exclude {\em a priori} a positive answer to question ({\cal
Q}). Actually, in spite of the non--coercivity of $\A$ on $\lamTP$, motions
defined by \condA, \condB, \condC\ may correspond to local minimizers, that
may exist even though (\ref{coercond}) is not satisfied.  We also remark that
for a loop $u(t)$ that satisfies (\ref{coercond}) we necessarily have
\begin{equation}
\bar u = \frac{1}{T}\int_0^T u(t)\,dt = 0\ .
\label{zeromean}
\end{equation}
This follows from (\ref{equivariance}), that implies
\[
\int_0^T u_i(t)\,dt = \int_0^T u_i(\tau(g)(t))\,dt = \rho(g)\int_0^T
u_{\sigma(g^{-1})(i)}(t)\,dt\,,
\]
which is equivalent to $\bar u_i =
\rho(g)\bar{u}_{\sigma(g^{-1})(i)}$. Clearly (\ref{zeromean}) poses
strong geometric restrictions on the motion and we may expect that
many complex and interesting motions with a rich geometric--kinematic
structure correspond to local minimizers that, as the loops defined by
\condA, \condB, \condC, do not need to satisfy condition
(\ref{zeromean}).  Is it possible to detect some of these local
minimizers?  Let $\mathfrak{S}\subset\lamTG$ be the subset of the
loops that present collisions:
\[
\mathfrak{S} =\{u\in\lamTG: \exists\, t_c\in \R, h\neq k\in\{1\dots
N\}: u_h(t_c) = u_k(t_c) \}\ .
\]
This is a well defined subset of $\lamTG$, closed in the $C^0$
topology.  We focus on {\em open cones} $\K\subset\lamTG$ with the property
\begin{equation}
\partial \K \subset \mathfrak{S} 
\label{bordercone}
\end{equation}
where $\partial \K$ is the $C^0$ boundary of $\K$.  The idea is that $\A$ can
be coercive on $\K$ even though (\ref{coercond}) is not satisfied; if we
are able to prove that a minimizer $\minloop$ of $\A|_{\Kclo}$ exists and is
collision free, then automatically we have $\minloop\in \K$ and therefore a
genuine solution of the $N$--body problem. In the following we discuss
non--trivial situations where the above ideas can be successfully applied.
Indeed we show the existence of new $T$--periodic solutions of the classical
$N$--body problem with a rather rich and complex structure.  In particular we
give a positive answer to question ({\cal Q}).
These results are precisely stated in Theorems \ref{teo4body}, \ref{teoplato},
\ref{platorb_exist} below.

\noindent
 We remark that restricting the action to a cone $\K\subset \lamTG$ that
satisfies (\ref{bordercone}) corresponds to the introduction of {\em
topological} constraints beside the {\em symmetry} constraints imposed by the
equivariance condition (\ref{equivariance}).  The idea to obtain coercivity by
introducing topological constraints, that restrict the action to subsets $\K$
satisfying condition (\ref{bordercone}), goes back to Poincar\'e
\cite{poincare} and has been exploited in \cite{gordon77}, \cite{venturelli},
\cite{terrvent}.  In the proof that minimizers $\minloop$ of $\A|_{\Kclo}$,
for the considered cones $\K$, are free of collisions we take advantage of
ideas of various authors \cite{marchal01}, \cite{chenciner03}, \cite{ventPhD}.
Moreover we need to overcome the extra difficulty of dealing with a
topological constraint, that does not allow general perturbations of
$\minloop\in\partial\K$, but only those that move $\minloop$ inside $\K$.  For
this reason Marchal's idea of averaging the action on a sphere or the average
on a suitable circle, leading to the definition of the rotating circle
property in \cite{FT04}, can not be applied in our context. The paper is
organized as follows: in Sections~\ref{sec:cones},~\ref{sec:platocones_one}
and \ref{sec:platocones_two} we define different kinds of cones $\K$ that
satisfy (\ref{bordercone}), and we prove the coercivity of $\A|_{\K}$. In
Section~\ref{sec:coll} we prove Theorems~\ref{teo4body}, \ref{teoplato},
\ref{platorb_exist} by showing that, for certain cones $\K$, minimizers
$\minloop\in\K$ are collision free.  In Section~\ref{sec:tot_coll} we exclude
total collisions by means of level estimates. In Section~\ref{sec:part_coll}
we exclude partial collisions via local perturbations. In
Section~\ref{sec:conj} we present conjectures and numerical experiments, and
prove the existence of $T$--periodic motions that violate (\ref{zeromean}).

\section{An example of cone $\K$}
\label{sec:cones}
We consider the set of $T$-periodic loops of $N=4$ unit masses
that satisfy

\begin{equation}
\left\{
\begin{array}{l}
u_1(t) = S_3 u_1(-t)\cr
u_1(\frac{T}{4} +t) = S_2 u_1(\frac{T}{4} -t)\cr
\end{array}\right.
\label{4body_cond1}
\end{equation}
\begin{equation}
\left\{
\begin{array}{l}
u_2(t) = R_3u_1(t)\cr
u_3(t) = R_1u_1(t)\cr
u_4(t) = R_2u_1(t)\cr
\end{array}\right.
\hskip 1.5cm
\label{4body_cond2}
\end{equation}
where $S_j$ is the reflection with respect to the plane $\assej = 0$, $j =
1,2,3$, and $R_j$ is the rotation of $\pi$ around the axis $\assej$, $j
=1,2,3$ (cfr. Figure~\ref{fig:fourbody}).
%
%
The loop space defined by (\ref{4body_cond1}) and
(\ref{4body_cond2}) does not satisfy condition (\ref{coercond}). Indeed the
vector $x \in \mathcal{X}\subset \R^{12}$ defined by
\[
x_1 = x_3  = \euno\,,
\quad
x_2 = x_4  = -\euno
\]
is in $\mathcal{X}^G$ and therefore, by Theorem~\ref{coerteor}, $\mathcal{A}$
is not coercive on $\Lambda_G$ as one can also directly verify.  Following the
approach outlined above, we regain coercivity by restricting the set of
allowed loops to the cone
\begin{equation}
\Kfour = \{u \in \Lambda_G: u_{11}(0)u_{11}(T/4) < 0\}\ .
\label{conecond4b}
\end{equation}
We denote with $(u_{11},u_{12},u_{13})$ the components of $u_1$.
\begin{proposition}
The cone $\Kfour$ satisfies condition (\ref{bordercone}) and
$\mathcal{A}|_{\Kfour}$ is coercive.
\label{conecoerc}
\end{proposition}
\begin{proof}
From (\ref{conecond4b}) it follows $[u_{11}(0) - u_{11}(T/4)]^2$ $>
(u_{11}(0))^2 + (u_{11}(T/4))^2$. This and (\ref{4body_cond1}), that implies
$u_{13}(0) = u_{12}(T/4) = 0$, yield
\begin{equation}
\left(|u_{1}(0)|^2 + |u_{1}(T/4)|^2\right)^{\frac{1}{2}} <
\left|u_{1}(0) - u_{1}(T/4)\right|.
\label{u10T4_bound}
\end{equation}
Assume $\genp^k$ and $t^k$, $k\in\N$, are sequences such that
$\Vert\genp^k\Vert_{C^0} = |\genp^k(t^k)| \to +\infty$ as $k\to +\infty$.
Then, if the sequence $\genp^k(0)$ is bounded, we have
$\lim_{k\to+\infty}|\genp^k(0) - \genp^k(t^k)| = +\infty$.  On the other hand,
if the sequence $\genp^k(0)$ is unbounded, (\ref{u10T4_bound}) implies that
$\lim_{k\to+\infty}|\genp^k(0) - \genp^k(T/4)| = +\infty$. This proves
coercivity.

\begin{figure}[ht]
\centerline{\psfig{figure=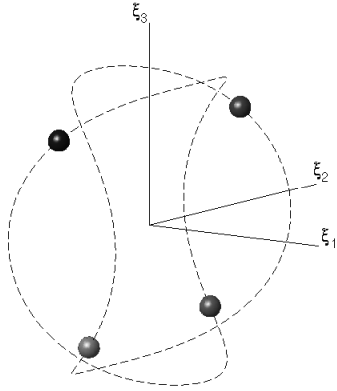,width=4.5cm}}
\caption{Structure of a motion (the map $\minloop$ in Theorem~\ref{teo4body})
that satisfies (\ref{4body_cond1}), (\ref{4body_cond2}).}
\label{fig:fourbody}
\end{figure}

\noindent To complete the proof we observe that $u$ belongs to
$\partial \Kfour$ if and only if one of the following is true: $u_{11}(0) =
0;\, u_{11}(T/4) = 0$.  As already remarked $u_{13}(0) = 0$. Therefore
$u_{11}(0) = 0\Rightarrow u_{1}(0) = u_{12}(0) \edue = u_4(0)$ by
(\ref{4body_cond2})$_3$. The other case is analogous.

\rightline{\small$\square$}
\end{proof}

Proposition~\ref{conecoerc} is the first step in the proof of the
following Theorem, that we establish in Section~\ref{sec:coll}:
\begin{theorem}
There exists a $T$-periodic solution $\minloop\in\Kfour$ of the classical
Newtonian $4$-body problem.
\label{teo4body}
\end{theorem}
In Figure~\ref{fig:fourbody} we show the geometry of the solution $\minloop$
computed numerically.

\section{Cones $\K$ and Platonic Polyhedra I.}
\label{sec:platocones_one}

The space $\lamTP$ characterized by $\condA, \condB, \condC$\ includes maps
$u:\R\to\R^{3N}$, such that the map $\genp:\R\to\R^3$ of the
generating particle satisfies
\begin{equation}
u_{11}(t_m) = \min_{t\in\R} u_{11}(t) > 0\,,
\label{minu01}
\end{equation}
in contrast with (\ref{zeromean}).  From (\ref{minu01}) it follows that, if
$u$ is a real motion of the $N$--body problem, then, at time $t_m$, the $N$
masses of the system are distributed on both sides of the plane $\asseuno
=u_{11}(t_m)$. Indeed (\ref{minu01}) and \condA\ imply that some of the masses
lie in the half--space $\asseuno <u_{11}(t_m)$. From this and (\ref{minu01}),
that implies $\ddot u_{11}(t_m)\geq 0$, it follows $u_{j1}(t_m) > u_{11}(t_m)$
for some $j\neq 1$.  This suggests that in the search of local minimizers of
$\left.\A\right|_{\lamTP}$ one should concentrate on situations of the type
shown in Figure~\ref{fig:motions_cube}$a$ and disregard the ones sketched in
Figure~\ref{fig:motions_cube}$b$.  Moreover it is clear from the picture that
a continuous deformation from a loop $u\in\lamTP$ of the type in
Figure~\ref{fig:motions_cube}$a$ into a loop of the type of
Figure~\ref{fig:motions_cube}$b$ and vice-versa can not be done without
collisions.  It is also clear that if the $H^1$--norm 
of a map $u \in \lamTP$, corresponding to Figure~\ref{fig:motions_cube}$a$,
tends to infinity, also the diameter of the orbit of each particle tends to
infinity and we have coercivity in spite of the fact that (\ref{coercond}) is
violated.  Based on these observations we now define cones $\K \subset \lamTP$
that satisfy (\ref{bordercone}).  Given unit vectors $\er, \es$, $\er\neq \pm
\es$, we let $\widehat{\er \es}$ be the angle
\[
\widehat{\er \es} = \{x\in\R^3: x = a \er + b \es, a,b>0 \}\ .
\]
Let $r_M$ ($r_V$) be the line through $OM$ ($OV$) and let $\eM$ ($\eV$) be the
unit vectors directed as $OM$ ($OV$). Observe that, beside $\asseuno$ and
$r_M$ ($\asseuno$ and $r_V$), on the plane $\asseuno r_M$ ($\asseuno r_V$)
there is at least another axis $r_\alpha$ ($r_\beta$) of some rotations in
${\cal R}\setminus \{I\}$.  We choose $r_\alpha \not\in \{\asseuno, r_M\}$
($r_\beta \not\in \{\asseuno, r_V\}$) and its direction vector $\ealpha$
($\ebeta$) such that $M \in \widehat{\euno \ealpha}$ ($V \in \widehat{\euno
\ebeta}$) and the measure $\psi_\alpha$ ($\psi_\beta$) of the angle
$\widehat{\euno \ealpha}$ ($\widehat{\euno \ebeta}$) is minimum. We also let
$\phi_M$ ($\phi_V$) be the measure of $\widehat{\euno \eM}$ ($\widehat{\euno
\eV}$).

We define
\begin{equation}
\begin{array}{ll}
\KPuno &= \bigl\{u\in\lamTP: \genp(0) \in\widehat{\eM \ealpha},\ \ 
\genp(\frac{T}{2H}) \in \widehat{\euno \eV}\bigr\}\,,\cr
\end{array}
\label{cone1}
\end{equation}
\begin{equation}
\begin{array}{ll}
\KPdue &= \bigl\{u\in\lamTP: \genp(0) \in\widehat{\eM \ealpha},\ \  
\genp(\frac{T}{2H}) \in \widehat{\eV \ebeta}\bigr\}\,,\cr
\end{array}
\label{cone2}
\end{equation}
\begin{equation}
\begin{array}{ll}
\KPtre &= \bigl\{u\in\lamTP: \genp(0) \in\widehat{\euno \eM},\ \ 
\genp(\frac{T}{2H}) \in \widehat{\eV \ebeta}\bigr\}\ .\cr
\end{array}
\label{cone3}
\end{equation}

\begin{figure}[h]
\centerline{{\small $\KPuno$}
\epsfig{figure=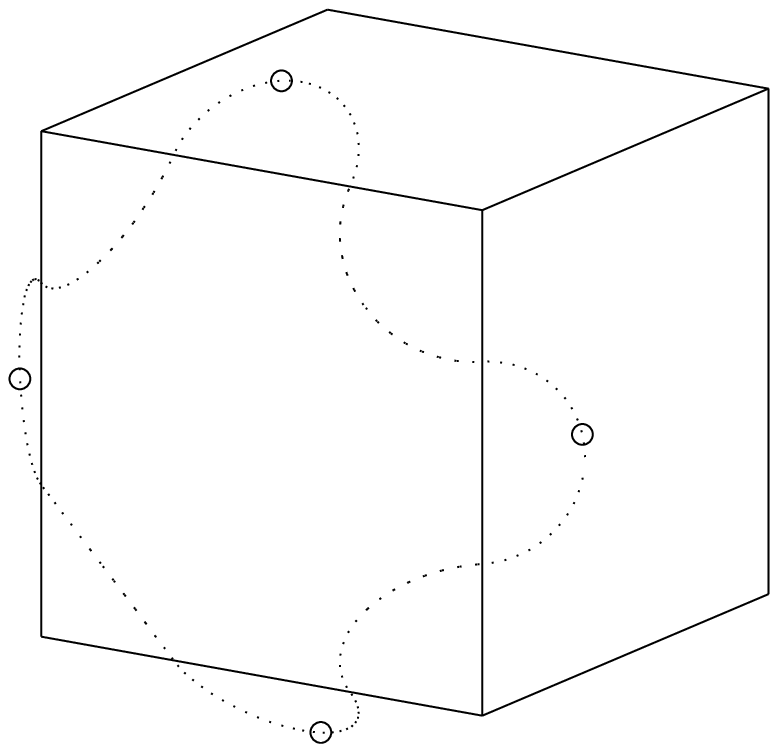,width=3.1cm}\hskip 0.2cm
{\small $\KPdue$}
\epsfig{figure=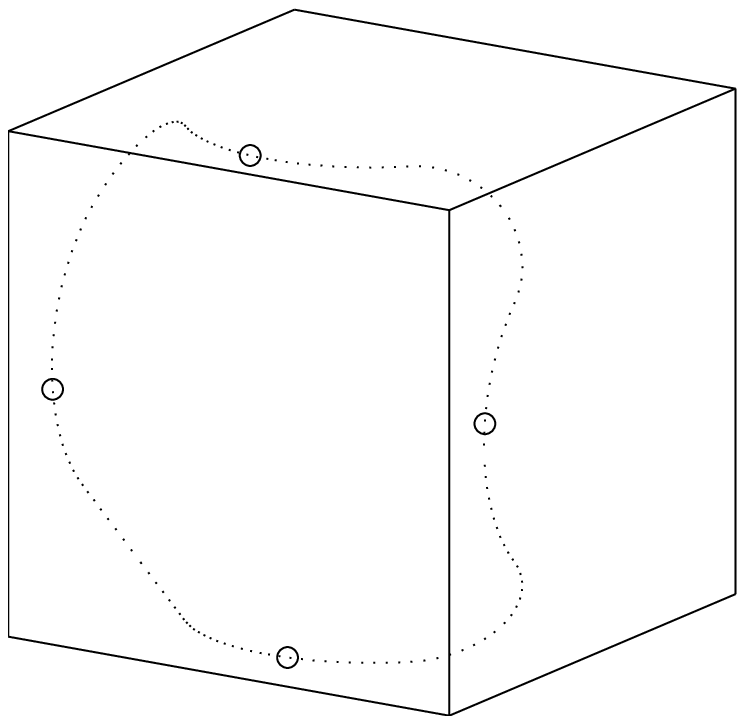,width=3.1cm}\hskip 0.2cm
{\small $\KPtre$}
\epsfig{figure=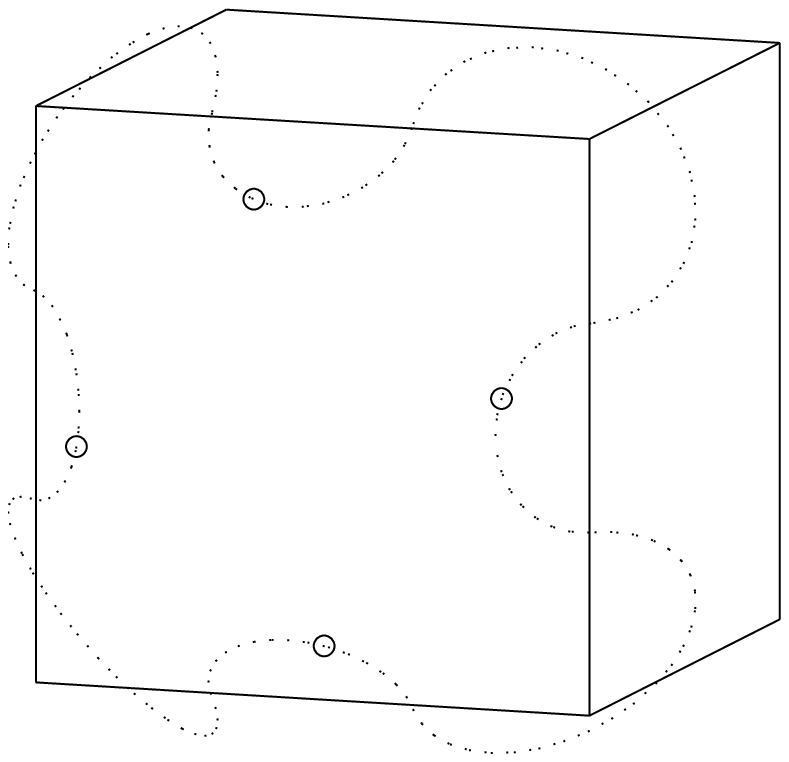,width=3.1cm}
}
\caption{The trajectory $\tau_1$ for typical loops in the
cones $\KPuno, \KPdue, \KPtre$ respectively, in the case $P=\Cplato$.}
\label{figcones}
\end{figure}

\noindent The definition of $\KPuno$ is illustrated in Figure~\ref{figcones}
for the case $P=\Cplato$.  We remark that
\begin{equation}
\KPj \cap \KPh = \emptyset \quad \mbox{ for } j\neq h\ .
\label{disjointcones}
\end{equation}
\begin{proposition}
The cone $\KPi\subset \lamTP$, $i=1,2,3$, satisfies (\ref{bordercone})
and $\A|_{\KPi}$ is coercive.
\label{border&coerc}
\end{proposition}
\begin{proof}

Let $R\in {\cal R}$ be the rotation of angle $2\pi/H$ around $\asseuno$. Then
\condB\ implies $R \genp(t-\frac{T}{H}) = \genp(t)$ for all $t$, and in
particular $R\genp (-\frac{T}{2H}) = \genp(\frac{T}{2H})$. Hence from \condC\
it follows $R S_3\genp(\frac{T}{2H}) = \genp(\frac{T}{2H})$. Since $R S_3$
coincides with the reflection with respect to the plane $\asseuno r_V$ we
conclude that, for $\typloop \in \lamTP$, $\genp(\frac{T}{2H})$ lies on the
plane $\asseuno r_V$. Moreover by \condC, $\genp(0)$ lies in the plane
$\asseuno r_M$. Therefore $\typloop \in\partial\KPi$ implies that either
$\genp(0)$ or $\genp(\frac{T}{2H})$ belongs to the boundary of one of the
angles in the definition of $\KPi$, that is to the axis $r$ of some rotation
in ${\cal R}\setminus\{ I\}$.  Thus by \condA\ we have a collision of the
generating particle with all the other particles associated to the maximal
cyclic group of the rotations with axis $r$. This proves (\ref{bordercone}).

\noindent To show coercivity of $\A|_{\KPi}$ we observe that from the
definition of $\KPi$ and \condB\ it follows the existence of constants
$c^P_i>0$, depending only on $P$ and $i=1,2,3$, such that at least one of the
two following inequalities holds true:
\begin{equation}
\Bigl|\genp(\frac{T}{H}) - \genp(0)\Bigr| = |R \genp(0) - \genp(0)|
 \geq c^P_1|\genp(0)|\,,
\label{lowerest0}
\end{equation}
\begin{equation}
\Bigl|\genp(\frac{T}{2H}) - \genp(-\frac{T}{2H})\Bigr| = \Bigl|R
 \genp(-\frac{T}{2H}) - \genp(-\frac{T}{2H})\Bigr| \geq
 c^P_i\Bigl|\genp(-\frac{T}{2H})\Bigr|\,,
\label{lowerestT2H}
\end{equation}
for all $u\in\KPi$.  If $\bar t\in(0,T)$ is such that $|\genp(\bar t)| =
\Vert\genp\Vert_{C^0}$, then from (\ref{lowerest0}), (\ref{lowerestT2H}) we
have either
\[
\begin{array}{ll}
2T^{1/2}{\left(\int_0^T |\dgenp|^2\right)}^{1/2} &\geq
\left|\genp(\frac{T}{H})-\genp(0)\right| + \left|\genp(\bar
t)-\genp(0)\right|\cr
&\geq c^P_i|\genp(0)| + \bigl|\Vert\genp\Vert_{C^0} - |\genp(0)|\bigr|\cr
\end{array}
\]
or
\[
\begin{array}{ll}
2T^{1/2}{\left(\int_0^T |\dgenp|^2\right)}^{1/2} &\geq
|\genp(\frac{T}{2H})-\genp(-\frac{T}{2H})| + \left|\genp(\bar
t)-\genp(-\frac{T}{2H})\right| \cr
&\geq c^P_i\left|\genp(-\frac{T}{2H})\right| +
\left|\Vert\genp\Vert_{C^0} - |\genp(-\frac{T}{2H})|\right|
\cr
\end{array}
\]
and coercivity follows.

\rightline{\small$\square$}
\end{proof}

On the basis of Proposition~\ref{border&coerc} we shall prove the following
theorem (Section~\ref{sec:coll}):
\begin{theorem}
Given $P\in\{\Tplato,\Cplato,\Oplato,\Dplato,\Iplato\}$ let $\KPi \subset
\lamTP$, $i=1,2,3$, be the cones defined in (\ref{cone1}), (\ref{cone2}),
(\ref{cone3}):
\begin{itemize}
\item[(i)] there exists a minimizer $\minloopPi\in\KPi$ of
$\A|_{\KPi}$ , $i=1,2,3$, and $\minloopPi$ is a smooth $T$--periodic
solution of the classical Newtonian $N$--body problem ($N=12$ for
$P=\Tplato$; $N=24$ for $P=\Cplato, \Oplato$; $N=60$ for
$P=\Dplato,\Iplato$).
\item[(ii)] $\minloopPi \neq \minloopPj\,, i\neq j$.
\end{itemize}
\label{teoplato}
\end{theorem}
$T$--periodic solutions of the classical Newtonian $N$--body problem
for $N=|{\cal R}|$ particles, ${\cal R}\in\{{\cal T}, {\cal O}, {\cal
I}\}$, satisfying condition \condA\ and the coercivity condition
(\ref{zeromean}) have already appeared in \cite{ferrario} and
\cite{moore06}. We stress that we do not require condition
(\ref{zeromean}), see also Theorem~\ref{teo:not_meanzero} below.
\begin{figure}[h]
\leftline{\epsfig{figure=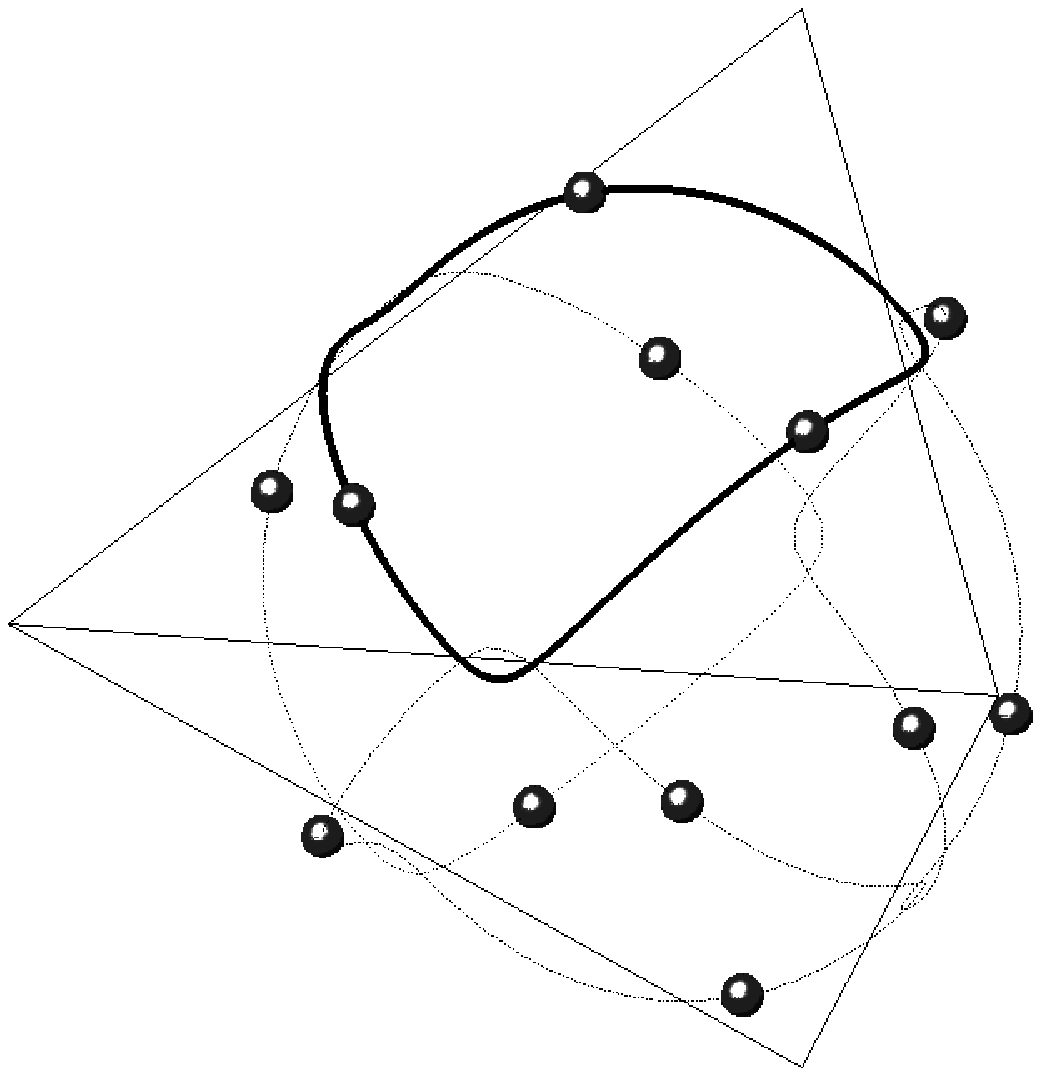,width=5cm}
\hskip 0.1cm $P=\Tplato$
\epsfig{figure=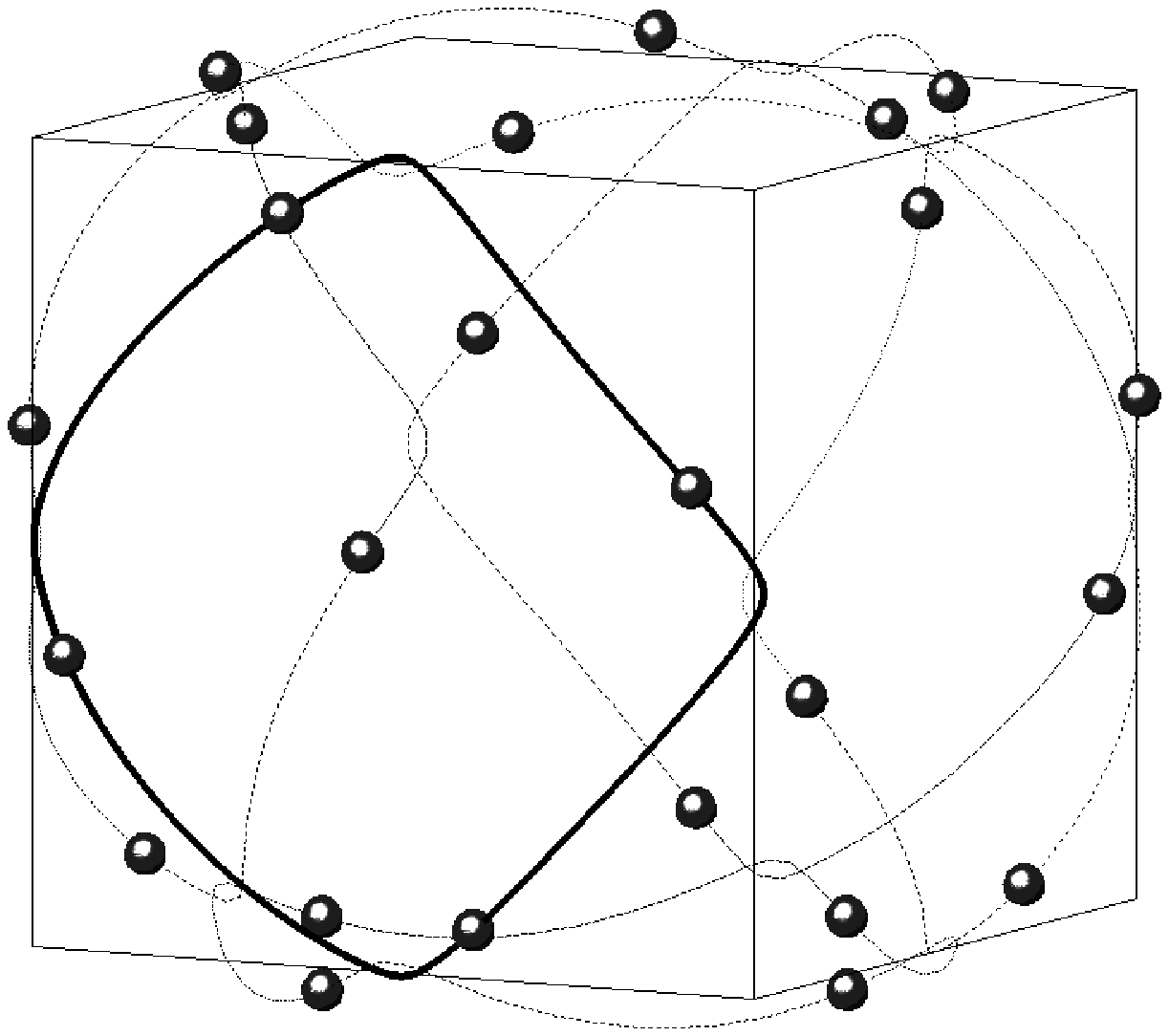,width=5.5cm}
\hskip 0.1cm$P=\Cplato$}
\vskip 0.2cm
\centerline{\epsfig{figure=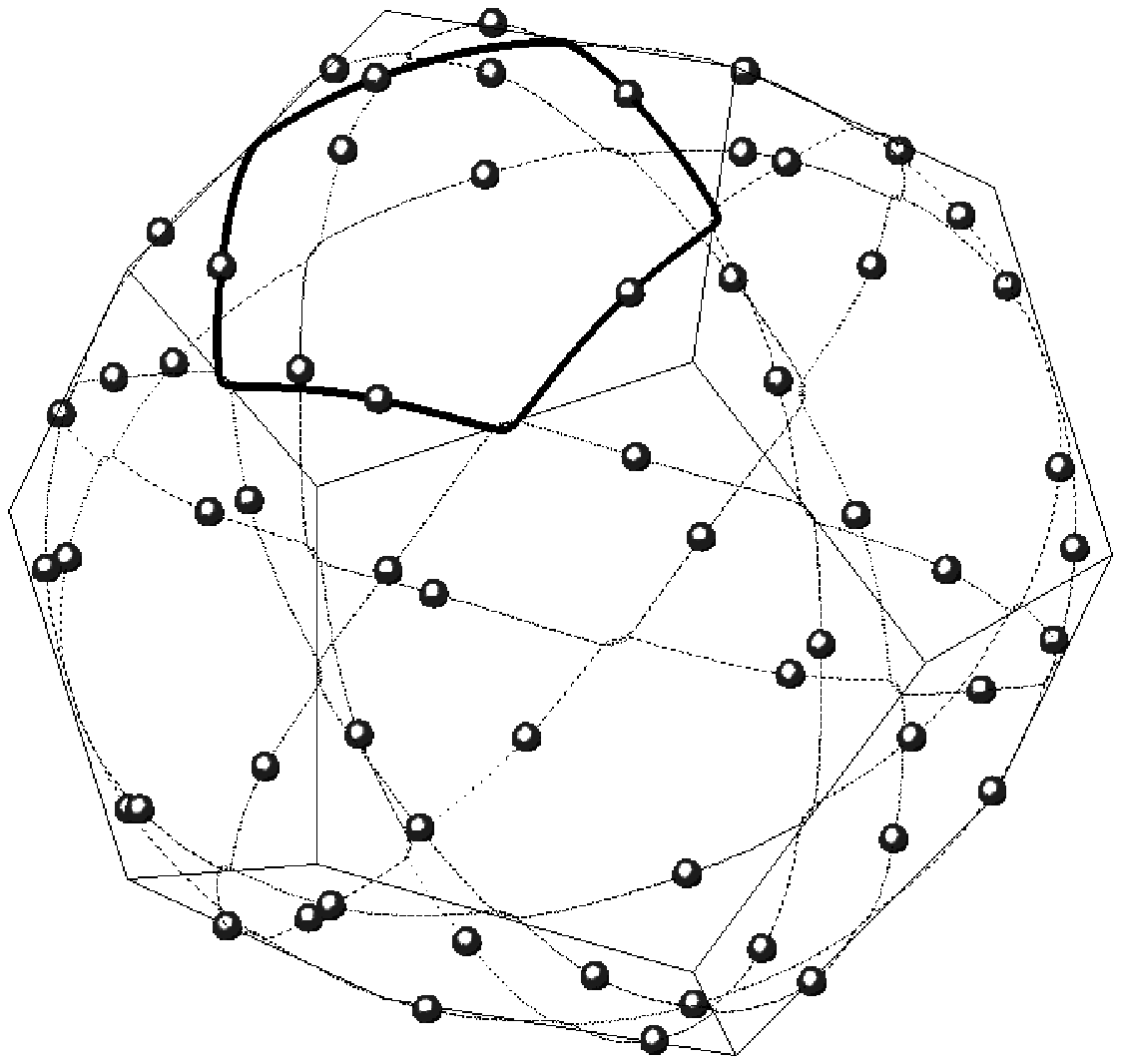,width=6.5cm}
\hskip 0.1cm$P=\Dplato$
}
\caption{Geometry of minimizers $\minloop^{P,1}$ for the cones $\KPuno$, with
$P = \Tplato, \Cplato, \Dplato$.}
\label{orbits_K123}
\end{figure}
Figure~\ref{orbits_K123} shows numerical simulations of the trajectory
of the particles for some of the minimizers $\minloopPi$ in
Theorem~\ref{teoplato}. 

\section{Cones $\K$ and Platonic Polyhedra II.}
\label{sec:platocones_two}

Definitions~(\ref{cone1}), (\ref{cone2}), (\ref{cone3}) give only three
examples of cones $\K$ satisfying (\ref{bordercone}) and ensuring coercivity
of $\A|_\K$. Actually for each choice of ${\cal R}\in\{{\cal T}, {\cal O},
{\cal I}\}$ there are infinitely many cones $\K$ with these properties.  For
each $R\in{\cal R}\setminus\{I\}$, let $r(R)$ be the axis of rotation of $R$
and define
\[
\Gamma = \bigcup_{R\in{\cal R}\setminus \{I\}} r(R) \subset \R^3\ .
\]
We denote by $\lama \subset \accauno$ the space of $T$--periodic maps $u:\R
\to \mathcal{X}$ defined by \condA\ when the map $\genp \in \accaunoR3$,
representing the motion of the generating particle, satisfies the condition
\begin{equation}
\genp(\R) \cap \Gamma = \emptyset\ .
\label{no_int_con_gamma}
\end{equation}
We introduce a notion of equivalence on $\lama$.

\begin{definition}
We say that $u$ and $v \in\lama$ are equivalent and we write $u\sim v$ if the
corresponding maps of the generating particle $\genp$ and $\genpv$ are
homotopic in $\R^3\setminus\Gamma$, that is if there exists a continuous map
$h:\R\times [0,1] \to \R^3$, $T$-periodic in the first variable, and
such that
\begin{itemize}
\item[{\rm (i)}] $h(\R,s)\cap \Gamma = \emptyset,\ \forall s \in
[0,1]$\,,
\item[{\rm (ii)}] $h(\cdot,0) = \genp,\ h(\cdot,1)=\genpv$ .
\end{itemize}
\label{def:usimv}
\end{definition}

Let $\lamoa \subset \lama$ be the subset of all the maps $u$ satisfying the
following condition:

\smallbreak\noindent ({\cal C})\quad {\em $\genp$ is not homotopic to any map
$\genpv\in \accaunoR3$ of the form:
\[
\genpv(t) = \euno' + \delta\left[\cos \left(2\pi k\frac{t}{T}\right) \edue' +
\sin \left(2\pi k\frac{t}{T}\right) \etre'\right]
\]
where $\ej', \, j=1,2,3$ is an orthonormal basis, with $\euno'$ parallel to one
of the axes in $\Gamma$, $0 <\delta<<1$, $k \in \N\cup\{0\}$.}  

\smallbreak \noindent Let $\K(u) \in \lamoasim$ $\,$ denote the equivalence
class of $u \in \lamoa$. We have:
\begin{proposition} 
Each cone $\K\in \lamoasim$ satisfies (\ref{bordercone}) and $\A|_{\K}$ is
coercive.
\label{prop:AKu_coerc}
\end{proposition}
\begin{proof} $\typloop \in \partial\K$ if and only 
if there exists a time $\bar t$ and $R\in \RR\setminus\{I\}$ such that
$\genp(\overline t) \in r(R)$. Then from (\ref{partgen}),
writing $R\to j_R$ for the inverse of $j\to R_j$, we have: $\genp(\bar
t) = u_{j_R}(\bar t)= R \genp(\bar t)$ and therefore a collision. This
establishes (\ref{bordercone}).  To show coercivity we observe that condition
$({\rm C})$ above implies the existence of a constant $c_{\K}$ such that
\[
\max_{t_1,t_2\in[0,T]} |\genp(t_1) - \genp(t_2)| \geq c_{\K}
\min_{t\in[0,T]} |\genp(t)|,\ \forall \typloop\in\K\ .
\]
Therefore, if $t_m$ satisfies $|\genp(t_m)|=\min_{t\in[0,T]}|\genp(t)|$, we
have
\[
|\genp(t)| \leq |\genp(t_m)| + |\genp(t) -\genp(t_m)| \leq
 \left(1/c_{\K}+ 1\right) \max_{t_1,t_2\in[0,T]}
 |\genp(t_1)-\genp(t_2)|\ .
\] 
%
%
This concludes the proof.

\rightline{\small$\square$}
\end{proof}

On the basis of Proposition~\ref{prop:AKu_coerc} we shall show (see
Theorem~\ref{platorb_exist}) that for several $\K\in\lamoasim$ there
exists $\minloop\in\K$ corresponding to a smooth periodic motion of
the Newtonian $N$--body problem. For the precise statement of
Theorem~\ref{platorb_exist} and for the detailed analysis of the
existence of collisions we make use of two different ways of
characterizing the topology of the maps in a given
$\K\in\lamoasim$. Indeed we associate to $\K$ two different
topological invariants in the form of periodic sequences of integers.
The following Subsection is devoted to the definition of these
invariants.

\subsection{Characterizing the topology of $\K$}
\label{subsec:chartop}

For ${\cal R}\in\{{\cal T}, {\cal O}, {\cal I}\}$ we denote with $\tilde{\cal
R}\in \{\tilde{\cal T}, \tilde{\cal O}, \tilde{\cal I}\}$ the associated
reflection group.  A fundamental domain \cite{grove_benson} for $\tilde{\cal
R}$ can be identified with the open convex cone $\fundD\subset \R^3$ generated
by the positive $\asseuno$ semiaxis, and the semiaxes determined by $OM$ and
by $OV$.\footnote{The axes $\assej, j=1,2,3$, are defined as in
Section~\ref{sec:intro}. If ${\cal R}={\cal O}\; ({\cal R} ={\cal I})$, the
definition can be based indifferently on $\Cplato$ or $\Oplato$ (on $\Dplato$
or $\Iplato$).}  Then
\[
\tilde{R} \fundD\cap \fundD = \emptyset \quad \forall \tilde{R}\in\tilde{\cal
R}\setminus\{I\}\,,
\hskip 0,5cm
\bigcup_{\tilde{R}\in\tilde{\cal R}} \tilde{R}\overline{\fundD} = \R^3\ .
\]
This means that $\R^3$ is divided into exactly $|\tilde{\cal R}|$
non--overlapping {\em chambers} ($|\tilde{\cal R}| = 2|{\cal R}| =
24,48,120$ for ${\cal R} = {\cal T}, {\cal O}, {\cal I}$) each of
which is an isometric copy of $\fundD$.  Let ${\cal D} = \{ D
\subset\R^3 : D =\tilde{R}\fundD, \tilde{R}\in\tilde{\cal R}\}$. We
also let ${\fundS}_i, i=1,2,3$ be the (open) faces of $\fundD$ and let
${\cal S}=\{S\subset\R^3: S = \tilde{R} {\fundS}_i,
\tilde{R}\in\tilde{\cal R}, i=1,2,3\}$ be the set of the faces of all
the elements of ${\cal D}$. For each $S\in{\cal S}$ we define
$\tilde{R}_S \in\tilde{\cal R}$ as the reflection with respect to the
plane of $S$. We also note that each $S\in{\cal S}$ uniquely
determines a pair $D_S^i\in{\cal D}, i=1,2$ such that
$\overline{S}=\overline{D_S^1}\cap \overline{D_S^2}$ and $D_S^2 =
\tilde{R}_S D_S^1$.

We consider the set of sequences $\sigma = \{D_k \}_{k\in\Z} \subset
{\cal D}$ that satisfy
\begin{itemize}
\item[(I)] $\sigma$ is periodic: $\exists\; K\in\N$ such that
$D_{k+K} = D_k, k\in\Z$;
\item[(II)] $D_{k+1}$ is the mirror image of $D_k$ with respect to one of
the faces of $D_k$ and $D_{k+1}\neq D_{k-1}$;
\item[(III)] $\bigcap_{k\in\Z}\overline{D_k} = \{0\}$. 
\end{itemize}
We identify sequences that coincide up to translations, that is
$\sigma = \{D_k\}, \sigma' = \{D_k'\}$ are identified whenever there exists
$K\in\N$ such that
\begin{equation}
D_k = D'_{k+K}, \ k\in\Z\ .
\label{identif}
\end{equation}
Each sequence $\sigma$ satisfying (I), (II), (III) can be regarded as periodic
of period $nK_\sigma$, with $n\in\N$ and $K_\sigma$ the minimal period. We
write $(\sigma, n)$ to indicate that we are considering $\sigma$ with the
particular period $nK_\sigma$ and regard $(\sigma, n_1)$, $(\sigma, n_2)$,
with $n_1\neq n_2$, as different objects.

Our first algebraic characterization of the topology of maps in $\K$
is described in the following Proposition:
\begin{proposition}
Each pair $(\sigma,n)$ with $\sigma = \{D_k\}_{k\in\Z}$ satisfying {\rm (I)},
{\rm (II)}, {\rm (III)} and $n\in\N$ uniquely determines a cone
$\K\in\lamoasim$. Viceversa each $\K\in\lamoasim$ uniquely determines $n\in\N$
and (up to translation) a sequence $\sigma$ satisfying {\rm (I)}, {\rm (II)},
{\rm (III)}.
\label{first_alg_char}
\end{proposition}
\begin{proof}
To each pair $(\sigma,n)$ we associate in a canonical way a map $\usn
\in \lamoa$ as follows.  Let $\tau_k,k\in\Z$ be the spherical
triangle, intersection of $D_k$ with the unit sphere
$\unitsphere\subset\R^3$, and let $c_k\in\tau_k$ be the center of
$\tau_k$.  Set $\ell_k = |c_{k+1}-c_k|$, $\ell_0=0$, $L =
\sum_{j=1}^{nK_\sigma}\ell_j$ and define $\usn$ by
\begin{equation}
\begin{array}{ll}
\usngen(t) &= \displaystyle\Bigl(\sum_{j=0}^k\ell_j
-\frac{t}{T}L\Bigr)\frac{c_k}{\ell_k} + \Bigl(\frac{t}{T}L -
\sum_{j=0}^{k-1}\ell_j\Bigr)\frac{c_{k+1}}{\ell_k}\,,\cr
&\displaystyle\sum_{j=0}^{k-1}\frac{\ell_j}{L} \leq \frac{t}{T} \leq
\sum_{j=0}^k\frac{\ell_j}{L}\,,\hskip0.2cm
k=1,\ldots,nK_\sigma\ .\cr
\end{array}
\label{genpsigman}
\end{equation}
By construction $\usngen$ satisfies (\ref{no_int_con_gamma}) and
condition (C) by (III). It follows that $\usn\in\lamoa$ and therefore
each pair $(\sigma,n)$ uniquely determines an element $\K(\usn)\in
\lamoasim$. To prove that the viceversa is also true we shall show
that for each $u\in\lamoa$ there is a unique pair $(\sigma_u,n_u)$
such that
\[
u \sim {\rm u}^{(\sigma_u,n_u)}\ .
\]
By definition of $\lamoa$, for each $\typloop \in\lamoa$ there is $d_u>0$ such
that
\begin{equation}
d(\genp(\R),\Gamma) = d_u\ .
\label{dist_from_gamma}
\end{equation}
Therefore the set $\Theta = \{t\in\R: \genp(t)\in S, S\in{\cal S}\}$ is
closed. We note in addition that (\ref{dist_from_gamma})
yields
\[
\A(\typloop) <+\infty\,,\qquad\forall\; \typloop\in\lamoa\ .
\]
From this and (\ref{dist_from_gamma}) we conclude that, if $(t_1,t_2)$ is a
connected component of the open set $\R\setminus\Theta$ with the property that
\begin{equation}
\genp(t_1)\in S, \genp(t_2)\in S', \mbox{ for some } S,S'\in{\cal S}\,, \mbox{
  with } S\neq S',
\label{genpt1t2}
\end{equation}
then $t_2-t_1 > \delta_\typloop$ for some $\delta_u>0$. This
inequality shows that in any time interval of size $T$ there are only
a finite number, say $M$, of connected components of
$\R\setminus\Theta$ that satisfy (\ref{genpt1t2}).

\noindent Let $J$ be the set of the connected components of
$\R\setminus\Theta$ that satisfy (\ref{genpt1t2}) and consider two consecutive
intervals $(t_1,t_2)$, $(t_1', t_2') \in J$. The assumption that between $t_2$
and $t_1'$ there is no other interval belonging to $J$ implies the existence
of $S\in{\cal S}$ such that
\begin{equation}
\genp(t_2), \genp(t_1') \in S\,,
\hskip 1cm
\genp([t_2,t_1']) \subset D_S^- \cup S \cup D_S^+\,,
\label{ut2t1p}
\end{equation}
where $D_S^\pm\in\{D_S^1, D_S^2\}$ are determined by the conditions
\[
\genp((t_1,t_2)) \subset D_S^-\,,\qquad D_S^- \neq D_S^+\ .
\]
Define $\checkgenp$ by setting, for each interval $(t_1,t_2)\in J$,
\begin{equation}
\checkgenp(t) = \left\{
\begin{array}{ll}
\tilde{R}_S \genp(t) &\qquad t\in(t_2,t_1'), u_1(t)\in D_S^+\,,\cr
\genp(t) &\qquad\mbox{otherwise }\,,
\end{array}
\right.
\label{ucheck}
\end{equation}
as in the example in Figure~\ref{fig:check}.
\begin{figure}[h]
\psfragscanon
\psfrag{S}{$S$}
\psfrag{Ds-}{$D_S^-$}
\psfrag{Ds+}{$D_S^+$}
\psfrag{u1t}{$\genp(t)$}
\psfrag{u1ct}{$\checkgenp(t)$}
\psfrag{u1t1}{$\genp(t_1)$}
\psfrag{u1t2}{$\genp(t_2)$}
\psfrag{u1t1p}{$\genp(t_1')$}
\psfrag{u1t2p}{$\genp(t_2')$}
\centerline{\epsfig{figure=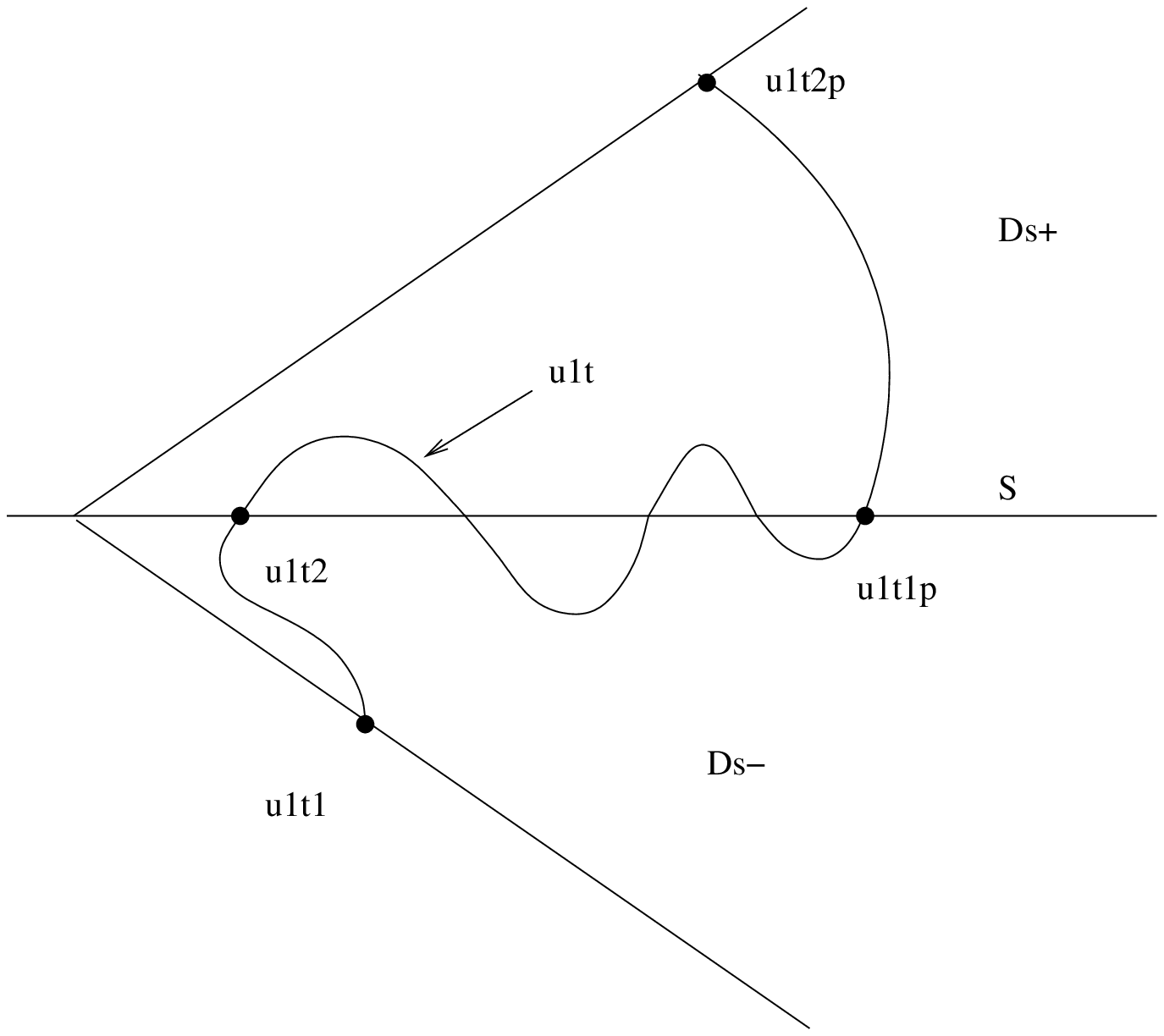,width=6cm}
\hskip 0.3cm
\epsfig{figure=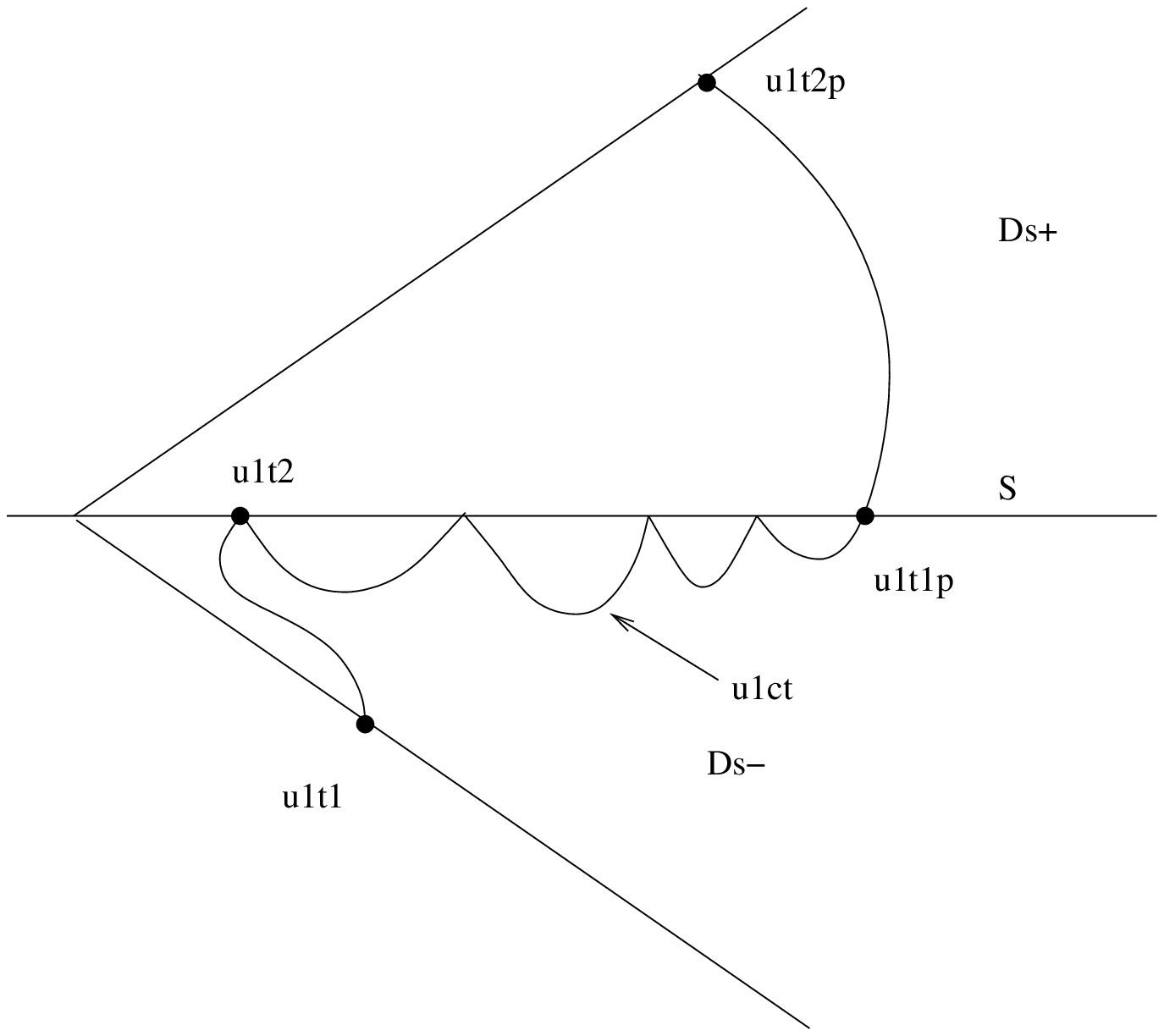,width=6cm}}
\psfragscanoff
\caption{Reflecting $\genp$ into $\checkgenp$.}
\label{fig:check}
\end{figure}
The convexity of $D_S^- \cup S \cup D_S^+$ implies that the map
\[
h(t,s) = (1-s)\genp(t) + s\checkgenp(t)\,, \quad  t\in\R, s\in[0,1]
\]
satisfies conditions (i), (ii) of Definition~\ref{def:usimv}, and therefore
$\check{\typloop}\sim \typloop$.

\noindent Write $J$ in the form $J = \{(t_1^k,t_2^k),\
t_1^k<t_1^{k+1},\ k\in\Z\}$, let $S_{k+1}\in{\cal S}$ be the face
associated to the interval $(t_2^k, t_1^{k+1})$, as in (\ref{ut2t1p}),
and let $D_k\in{\cal D}$ be the set $D_S^-$ corresponding to $S=S_{k+1}$,
that is
\begin{equation}
\genp\bigl((t_1^k,t_2^k)\bigr) \subset D_k\ .
\label{genpDk}
\end{equation}
From (\ref{ucheck}) and (\ref{genpDk}) it follows that 
\begin{equation}
\checkgenp\bigl([t_1^k,t_1^{k+1}]\bigr)\subset
\overline{D_k}\setminus\Gamma\,,\ k\in\Z
\label{ucheckinDk}
\end{equation}
and moreover that
\begin{equation}
\checkgenp\bigl((t_1^k,t_1^k+\delta_u)\bigr)\subset D_k\,,\ k\in\Z\
.
\label{ucheckdelta}
\end{equation}
From now on we shall drop the subscript 1 and write simply $t^k$ instead of
$t_1^k$. The periodicity of $\checkgenp$ implies
\[
[t^k + jT, t^{k+1}+jT] = [t^{k+jM}, t^{k+1+jM}], \ \ D_k = D_{k+jM}, \ \
j\in\Z\ .
\]
We can assume that 
\begin{equation}
D_k \neq D_{k+1}\ .
\label{DkDkp1}
\end{equation}
Indeed if $D_h = D_{h+1}$ for some $h\in\Z$, then (\ref{ucheckinDk}) implies
\[
\checkgenp\bigl([t^h+jT,t^{h+2}+jT]\bigr)\subset
\overline{D_{h+jM}}\setminus\Gamma
\]
and therefore, if we erase the subsequences $\{t^{h+1+jM}\}_{j\in\Z}$,
$\{S_{h+1+jM}\}_{j\in\Z}$, $\{D_{h+1+jM}\}_{j\in\Z}$ from the
sequences $\{t^k\}_{k\in\Z}$, $\{S_k\}_{k\in\Z}$, $\{D_k\}_{k\in\Z}$,
then, after relabeling, we still have that (\ref{ucheckinDk}),
(\ref{ucheckdelta}) hold. A finite number of steps of this kind
establishes (\ref{DkDkp1}).

By homotopy we can transform $\checkgenp$ into a map $\hatgenp$ that
satisfies (\ref{ucheckinDk}), (\ref{ucheckdelta}), \ref{DkDkp1}) and moreover
has the property that
\begin{equation}
D_{k-1} \neq D_{k+1}\ .
\label{Dkm1Dkp1}
\end{equation}
To see this we observe that $D_{h-1}=D_{h+1}$ together with
(\ref{ucheckinDk}), (\ref{DkDkp1}) imply $S_h = S_{h+1}$ and
\begin{equation}
\checkgenp(t^h), \checkgenp(t^{h+1})\in S_h\,,\hskip 1cm
\checkgenp\bigl([t^h,t^{h+1}]\bigr)
\subset\overline{D_h}\setminus\Gamma\ .
\label{ucheckinSk}
\end{equation}

\begin{figure}[h]
\psfragscanon
\psfrag{Sh}{$S_h$}
\psfrag{Dh}{$D_h$}
\psfrag{Dh-1Dh+1}{$D_{h-1} = D_{h+1}$}
\psfrag{uc1t}{$\check{u}_1(t)$}
\psfrag{uh1t}{$\hatgenp(t)$}
\psfrag{uc1th}{$\checkgenp(t^h)$}
\psfrag{uc1th+1}{$\checkgenp(t^{h+1})$}
\centerline{\epsfig{figure=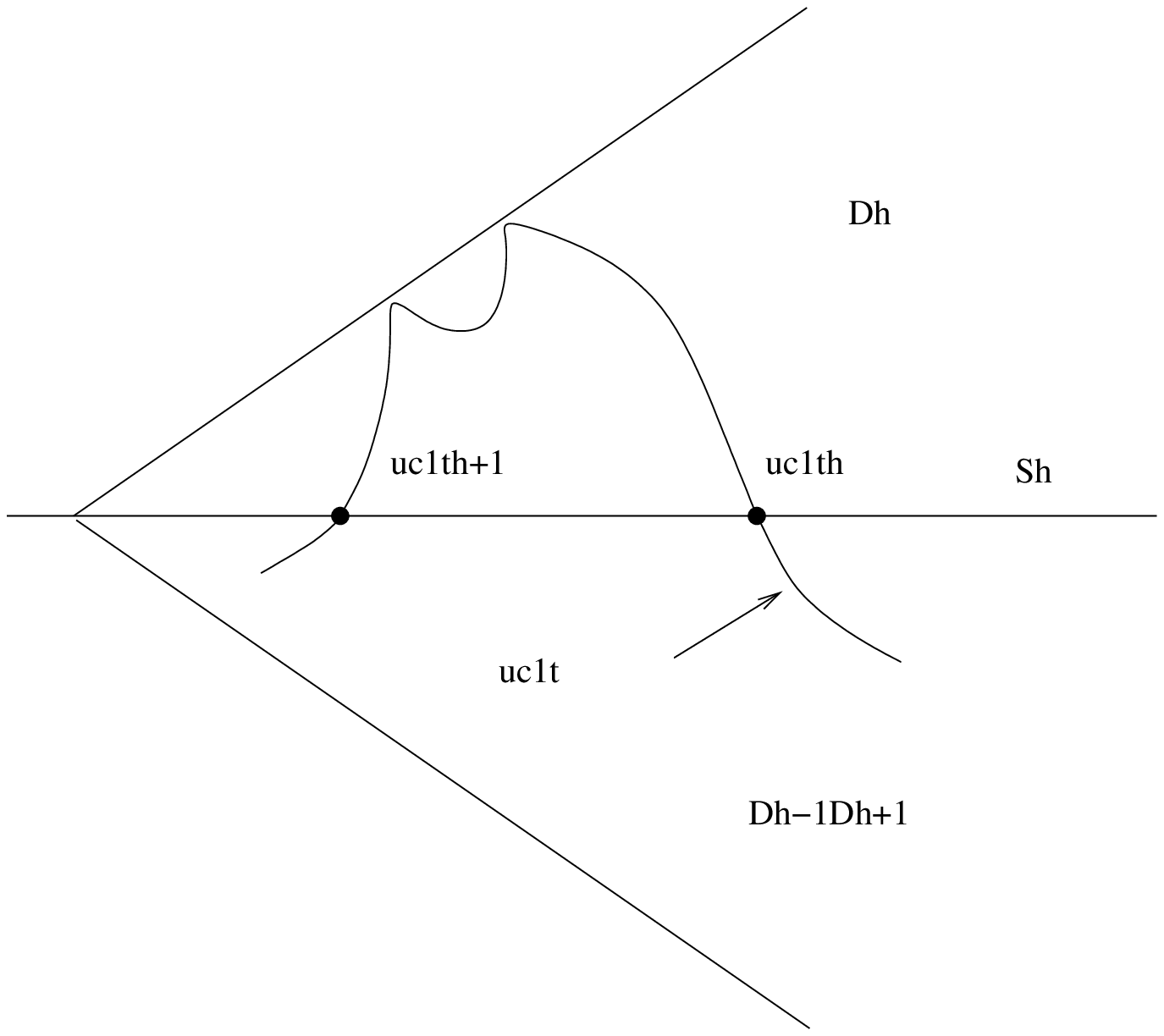,width=6cm}
\hskip 0.3cm
\epsfig{figure=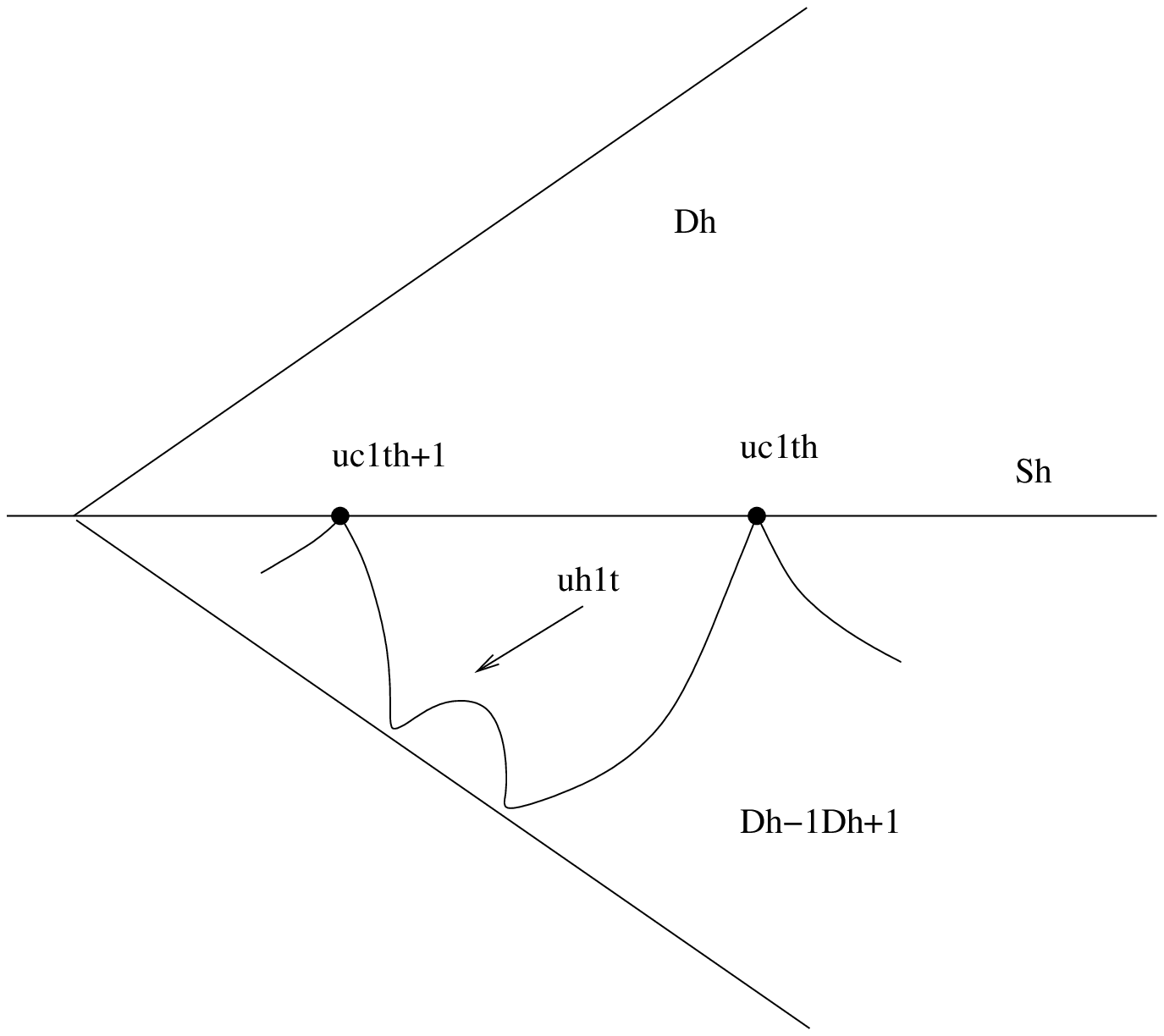,width=6cm}}
\psfragscanoff
\caption{Reflecting $\checkgenp$ into $\hatgenp$.}
\label{fig:hat}
\end{figure}
From (\ref{ucheckinSk}) and the same argument used above to define
$\check{\typloop}$ it follows that the map defined by
\[
\hatgenp(t) = 
\left\{
\begin{array}{ll}
\tilde{R}_{S_h} \checkgenp(t)\,,\qquad
 &t\in(t^{h+jM},t^{h+1+jM})\,,\; j\in\Z\,,\cr
\checkgenp(t)\,, &\mbox{otherwise}\cr
\end{array}
\right.
\]
(see Figure~\ref{fig:hat}) is homotopic to $\checkgenp$ and it has the
property that
\[
\hatgenp\bigl([t^{h-1},t^{h+2}]\bigr)\subset
\overline{D_{h-1}}\setminus\Gamma = \overline{D_{h+1}}\setminus\Gamma\ .
\]
This shows that, after erasing $t^{h+jM}$ and $t^{h+1+jM}$, $j\in\Z$, from
$\{t^k\}_{k\in\Z}$ ($S_{h+jM}, S_{h+1+jM}$, $j\in\Z$, from $\{S_k\}_{k\in\Z}$,
$D_{h+jM}, D_{h+1+jM}, j\in\Z$, from $\{D_k\}_{k\in\Z}$), and after relabeling,
the relations (\ref{ucheckinDk}), (\ref{ucheckdelta}), (\ref{DkDkp1}) still
hold and $M$ is reduced by 2 units. Therefore after a finite number of steps
we obtain a map $\hatgenp$ homotopic to $\genp$ and such that the
corresponding sequences $\{t^k\}_{k\in\Z}$, $\{S_k\}_{k\in\Z}$,
$\{D_k\}_{k\in\Z}$ satisfy (\ref{ucheckinDk})-- (\ref{Dkm1Dkp1}).  In
particular the sequence $\sigma_u = \{D_k\}_{k\in\Z}$ satisfies (I) and (II)
and we claim that $\hatgenp$ and therefore $\genp$ is homotopic to the map
$\usunugen$ where $n_u$ is defined by $n_u = M/K_u$, $K_u$ being the minimal
period of $\sigma_u$.  To prove the claim it suffices to observe that the
transformation $s \to h(\cdot,s)$, defined by
\[
h(\lambda,s) = (1-s)\hatgenp(t^k+\delta_u+\lambda(t^{k+1}-t^k)) +
s(c_k + \lambda(c_{k+1}-c_k)),\ \ \lambda,s\in[0,1]\,,
\]
continuously deforms, without touching $\Gamma$, each arc
$\hatgenp([t^k+\delta_u,t^{k+1}+\delta_u])\subset(\overline{D_k}\cup
\overline{D_{k+1}})\setminus\Gamma$ into the segment $[c_k, c_{k+1}]$
joining the centers of the spherical triangles $\tau_k=D_k\cap
\unitsphere$, $\tau_{k+1} = D_{k+1}\cap\unitsphere$. To conclude the
proof we only need to show that $\sigma_u$ satisfies (III). This
follows from the fact that, should (III) be violated, then $\usunugen$,
and in turn $\genp$, would not satisfy condition (C).

\rightline{\small$\square$}
\end{proof}

\noindent The above discussion establishes a one to one correspondence between
$\lamoasim$ and the set of the pairs $(\sigma,n)$ with $\sigma$ satisfying
(I), (II), (III) (with the identification (\ref{identif})).

\begin{definition}
A pair $(\sigma,n)$ is said to be {\em simple} if $\sigma$ does not
contain a string $D_k,\dots,D_{k+H}$ such that
\begin{itemize}
\item[a)] $\bigcap_{j=0}^H\overline{D_{k+j}} = r(R)$, for some $R\in{\cal
R}\setminus\{I\}$;
\item[b)] $H=2|{\cal C}|$, where ${\cal C}\subset {\cal R}$ is the
maximal cyclic group of the rotations around $r(R)$.
\end{itemize}
\label{def:simple}
\end{definition}
We say that $u\in\lamoa$ is simple if the corresponding $(\sigma_u,
n_u)$ is simple. On the basis of Definition~\ref{def:simple}, simple
$u$ are the ones such that $\genp$ does not coil around any of the
axes of the rotations in ${\cal R}$. We remark that the subset of
$\lamoasim$ of the cones corresponding to simple $u$ is infinite.


\medbreak Next we introduce our second algebraic characterization of the
topology of $\K\in\lamoasim$. We identify the faces $\fundS_1$ and $\fundS_2$
of $\fundD$, the fundamental domain introduced above, with the faces generated
by $\xi_1, OM$ and $OM, OV$ respectively.  Let $\Pi_i, i=1,2,3$ be the planes
of $\fundS_i$ and let $\tilde{R}_i\in\tilde{\cal R}$ be the reflection with
respect to $\Pi_i$. Since $\Pi_1$ and $\Pi_2$ are orthogonal there is a unique
$q\in\fundS_3\cap \unitsphere$ such that $q, \tilde{R}_1q, \tilde{R}_2q,
\tilde{R}_2\tilde{R}_1 q$ are the vertexes of a square. Since
$q=\tilde{R}_3q$, the orbit of $q$ under $\tilde{\cal R}$ contains only
$|\tilde{\cal R}|/2 = |{\cal R}|$ distinct points and coincides with the orbit
$\{Rq\}_{R\in{\cal R}}$ of $q$ under ${\cal R}$. The convex hull
$\mathcal{Q}_{\cal R}$ of $\{Rq\}_{R\in{\cal R}}$ is an Archimedean polyhedron
\cite{archimede} which is naturally associated to the rotation group ${\cal
R}$. Since the action of ${\cal R}$ on the square defined by
$q,\tilde{R}_1q,\tilde{R}_2q,\tilde{R}_3q$ generates $|\tilde{\cal R}|/4$
distinct squares, $\mathcal{Q}_{\cal R}$ has exactly $|\tilde{\cal R}|$ equal
sides and therefore all its faces ${\cal F}$ are regular polygons with axes
coinciding with $r(R)$ for some $R\in{\cal R}\setminus\{I\}$.

\begin{figure}[ht]
\centerline{
$\mathcal{Q}_{\cal T}$, (3434) \epsfig{figure=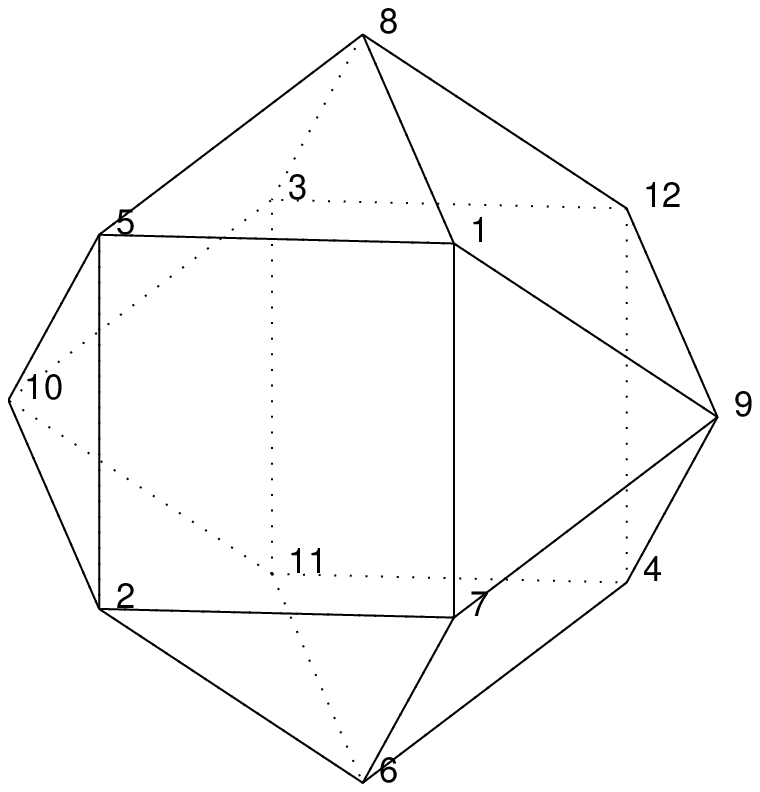,width=4.6cm}\hskip 0.3cm
$\mathcal{Q}_{\cal O}$, (3444) \epsfig{figure=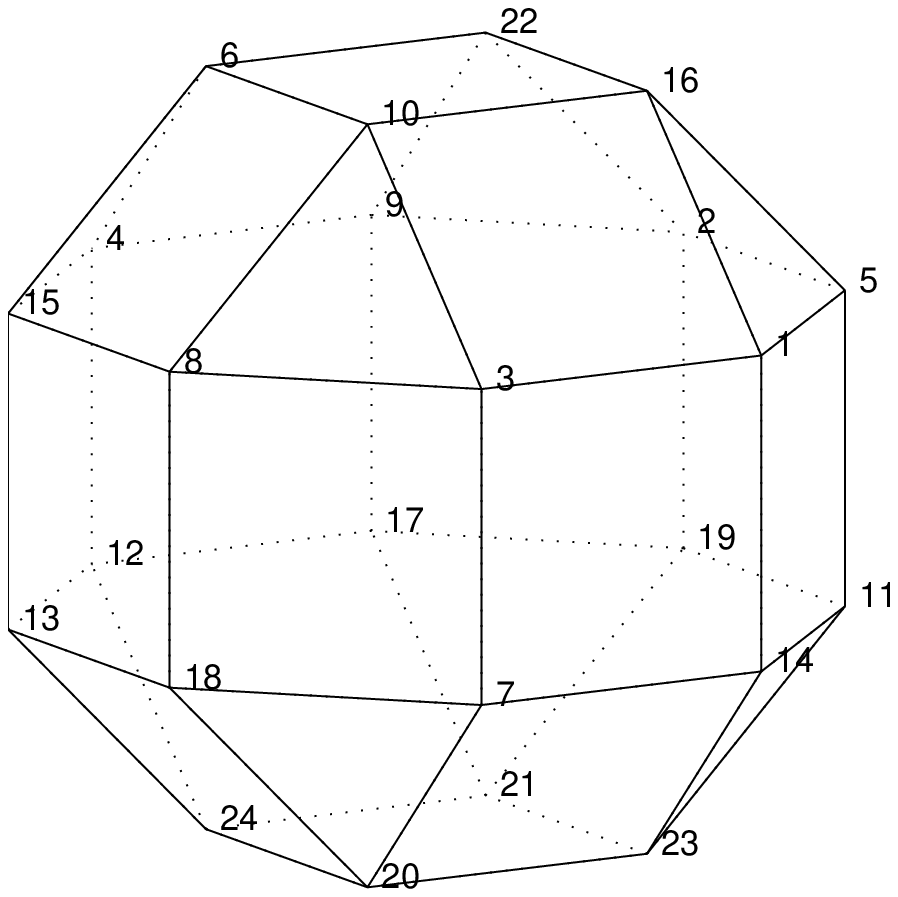,width=5cm}}
\vskip 0.1cm\centerline{
$\mathcal{Q}_{\cal I}$, (3454) \epsfig{figure=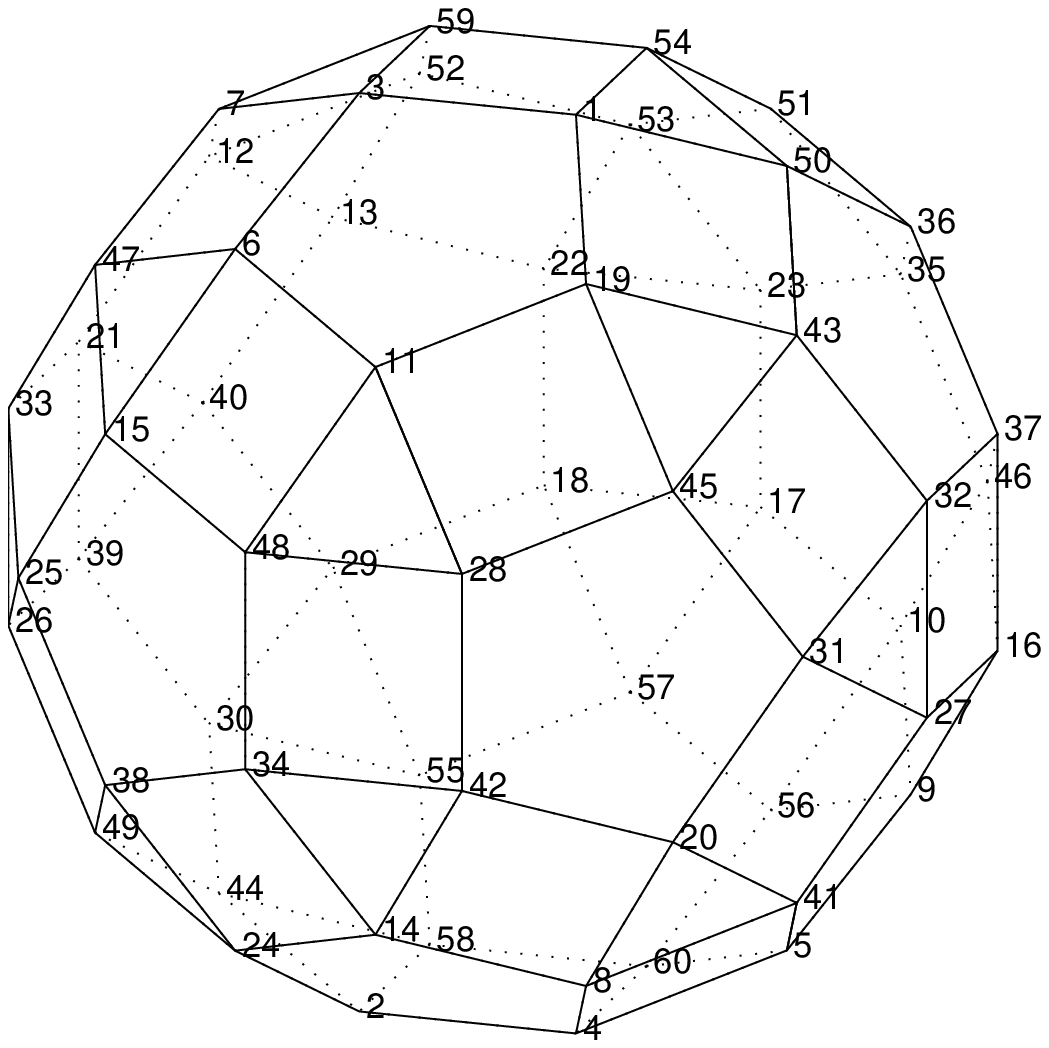,width=6cm}}
\caption{The three Archimedean polyhedra $\mathcal{Q}_{\cal T}, \mathcal{Q}_{\cal O},
\mathcal{Q}_{\cal I}$; the vertexes have been numbered for later reference.}
\label{fig:archimede}
\end{figure}

In Figure~\ref{fig:archimede} we show the three Archimedean Polyhedra
$\mathcal{Q}_{\cal T}, \mathcal{Q}_{\cal O}, \mathcal{Q}_{\cal I}$
corresponding to ${\cal R} = {\cal T}, {\cal O}, {\cal I}$
respectively.  The vertex configurations \cite{archimede} of
$\mathcal{Q}_{\cal T}, \mathcal{Q}_{\cal O}, \mathcal{Q}_{\cal I}$ are
$(3434), (3444), (3454)$.  By construction ${\cal L}_{\cal R}$, the
union of the edges of $\mathcal{Q}_{\cal R}$, avoids $\Gamma$. This
property allows us to introduce another algebraic characterization of
$\K\in\lamoasim$.

\begin{proposition}
Each $\K\in\lamoasim$ uniquely determines a number $n\in\N$ and (up to
translations) a periodic sequence $\nu=\{\nu_k\}_{k\in\Z}$ of vertexes
of $\mathcal{Q}_{\cal R}$ such that
\begin{itemize}
\item[{\rm[i]}] for each $k\in\Z$ the segment $[\nu_k,\nu_{k+1}]$ coincides
with one of the edges of $\mathcal{Q}_{\cal R}$;

\item[{\rm[ii]}] $\nu\not\subset\overline{\cal F}$, for all the faces
${\cal F}$ of $\mathcal{Q}_{\cal R}$.

\end{itemize}
Viceversa each pair $(\nu,n)$, $\nu$ a periodic sequence of vertexes
of $\mathcal{Q}_{\cal R}$ that satisfies {\rm[i]}, {\rm[ii]}, uniquely
determines a cone $\K\in\lamoasim$.
\label{prop_nu}
\end{proposition}

\begin{proof}
Let $c$ be the center of the spherical triangle
$\tau=\fundD\cap\unitsphere$ and set
\[
c_{\tilde{R}}(s) = (1-s)\tilde{R}c + s\tilde{R}q, \
s\in[0,1], \ \tilde{R}\in\tilde{\cal R}\ .
\]
Given $u\in\K$ let $(\sigma_u,n_u)$, $\sigma_u = \{D_k\}_{k\in\Z}$,
$n_u\in\N$, be the pair associated to $u$ by Proposition~\ref{first_alg_char}
and let $c_k(s)$ be determined by the condition
\[
c_k(s) = c_{\tilde{R}}(s)\,;\quad \tilde{R}\fundD = D_k\ .
\]
Moreover let $\usunusgen$ be the map defined by (\ref{genpsigman})
when $c_k, c_{k+1}$ are replaced by $c_k(s), c_{k+1}(s)$ and $\ell_k$
by $\ell_k(s) = |c_{k+1}(s) - c_k(s)|$. The map $h(t,s) = \usunusgen
(t), (t,s) \in\R\times[0,1]$ defines a homotopy that transforms
$\usunugen = \usunugenzero$ into $\usunugenuno$. By definition
$\usunugenuno$ ranges in ${\cal L}_{\cal R}$ and describes with
constant speed a closed path on ${\cal L}_{\cal R}$.  This determines
a sequence $\nu=\{\nu_j\}_{j\in\Z}$ consisting of the vertexes of
$\mathcal{Q}_{\cal R}$ visited one after the other by
$\usunugenuno$. If $\nu_j = \usunugenuno (t_j)$ we have
\[
\nu_j = \usunugenuno (t_j+T) = \nu_{j+K}\,,
\]
for some $K\in\N$. The integer $n$ is determined by $n=K/K_\nu$ with
$K_\nu$ the minimal period of $\nu$.  The sequence $\nu$ satisfies
[i], and also [ii] since otherwise $\usunuuno$ and therefore
$\typloop$ will not satisfy condition (C).

To prove the last statement of the Proposition, given $(\nu,n)$, $\nu$
a periodic sequence of vertexes of $\mathcal{Q}_{\cal R}$ satisfying
[i], [ii], we define the $T$-periodic map $\vnun$ by setting
\begin{equation}
\left\{
\begin{array}{l}
\gamma^{(\nu,n)} = \displaystyle\prod_{j=1}^{nK_\nu} \gamma_{\nu,j}\,,
\quad\gamma_{\nu,j}(s) = (1-s)\nu_j + s\nu_{j+1}\,,\quad s\in[0,1]\cr
{\rm v}_1^{(\nu,n)}(t) = \gamma^{(\nu,n)}(t/T)\ .\cr
\end{array}
\right.
\label{vnun}
\end{equation}
By definition the map $\vnun$ is $T$--periodic and $\vnungen(\R)\cap\Gamma =
\emptyset$. Moreover $\vnungen$ satisfies condition (C) by [ii]. Therefore we
have $\vnun\in\lamoa$.

\rightline{\small$\square$}
\end{proof}

\subsection{Results on the existence and the geometric structure 
of periodic motions in $\K$} 
\label{subsec:geostruct}

Based on the algebraic characterization of $\K\in\lamoasim$ discussed
in the previous Section we can now state the following
\begin{theorem}
For each one of the sequences $\nu$ listed below there exists a
$T$--periodic solution of the classical Newtonian $N$--body problem
which is equivalent to $\vnuuno$ in the sense of
Definition~\ref{def:usimv}.  (The sequences are given with reference
to the enumeration of the vertexes of $\mathcal{Q}_{\cal R}$ in
Figure~\ref{fig:archimede}).
\vskip 0.2cm
\leftline{${\cal R} = {\cal T}$}
\begin{eqnarray*}
 &&\nu^1 = [1, 8, 3, 12, 9, 7, 1]\,,\hskip 12cm\\
 &&\nu^2 = [1, 8, 3, 11, 6, 7, 1]\,,\\
 &&\nu^3 = [9, 1, 8, 12, 4, 6, 7, 9, 12, 3, 11, 4, 9]\,,\\
\end{eqnarray*}
\vskip 0.2cm
\leftline{${\cal R} = {\cal O}$}
\begin{eqnarray*}
 &&\nu^1 = [16, 5, 11, 14, 23, 20, 18, 8, 3, 10, 16]\,,\hskip 8cm \\
 &&\nu^2 = [1, 16, 10, 3, 7, 20, 23, 14, 1]\,,  \\ 
 &&\nu^3 = [1, 14, 23, 11, 5, 16, 1]\,, \\ 
 &&\nu^4 = [1, 5, 16, 1, 3, 7, 18, 20, 7, 14, 1]\,,\\ 
 &&\nu^5 = [1, 16, 10, 8, 18, 7, 3, 10, 6, 15, 8, 3, 1]\,, \\
 &&\nu^6 = [11, 5, 2, 22, 16, 10, 6, 15, 8, 18, 13, 24,
 20, 23, 21, 19, 11]\,,
\end{eqnarray*}
\vskip 0.2cm
\leftline{${\cal R} = {\cal I}$}
\begin{eqnarray*}
 &&\nu^1 = [11, 48, 34, 14, 42, 28, 11]\,,\\
 &&\nu^2 = [11, 48, 34, 42, 28, 11, 6, 15, 48, 28, 45, 19, 11]\,,\\
 &&\nu^3 = [15, 48, 34, 42, 28, 45, 31, 32, 43, 50, 36, 51, 54,
59, 52, 12, 7, 47, 33, 25, 15]\ .
\end{eqnarray*}
\label{platorb_exist}
\end{theorem}

\begin{figure}[ht]
\centerline{
$\nu^1$ \epsfig{figure=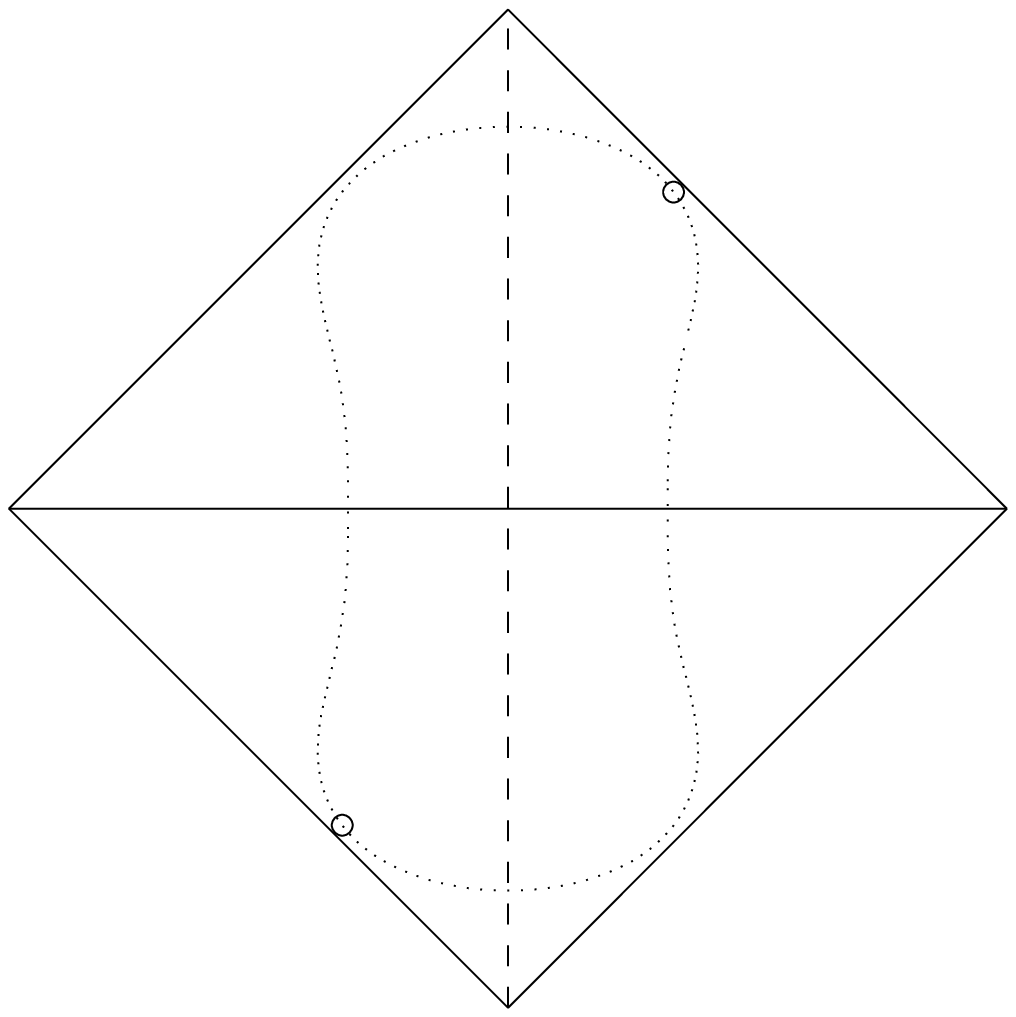,width=3.5cm}
\hskip 0.3cm
$\nu^2$ \epsfig{figure=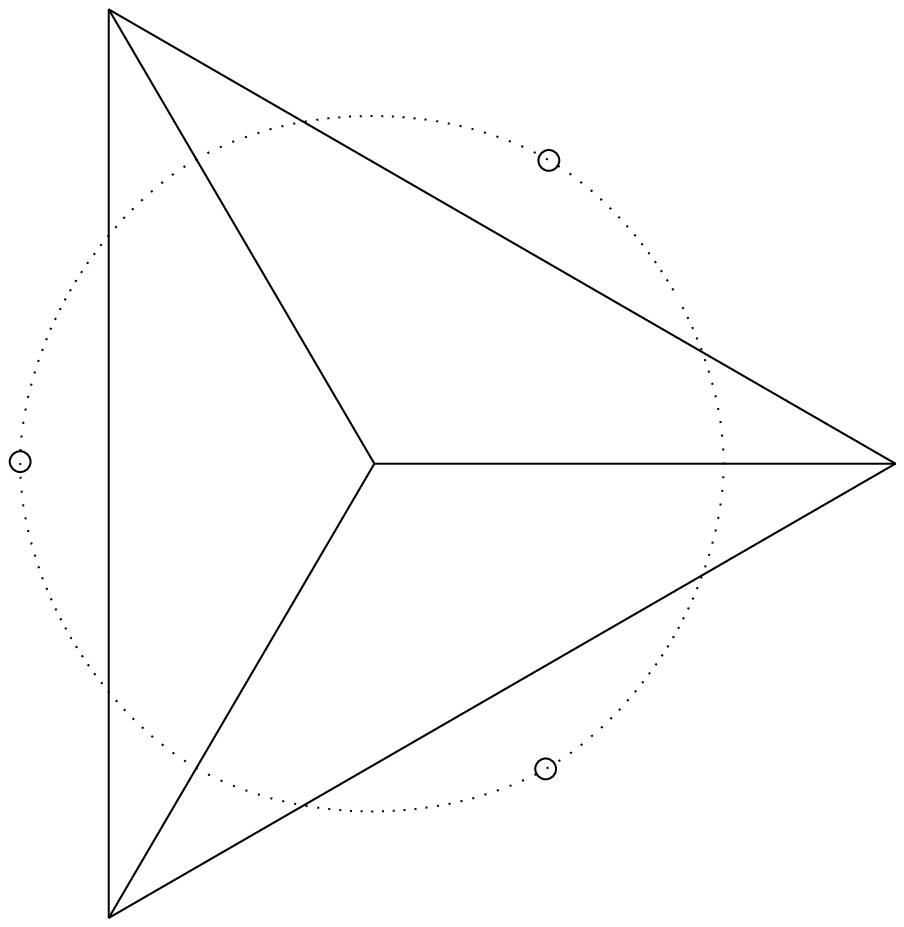,width=3.2cm}
\hskip 0.3cm
$\nu^3$ \epsfig{figure=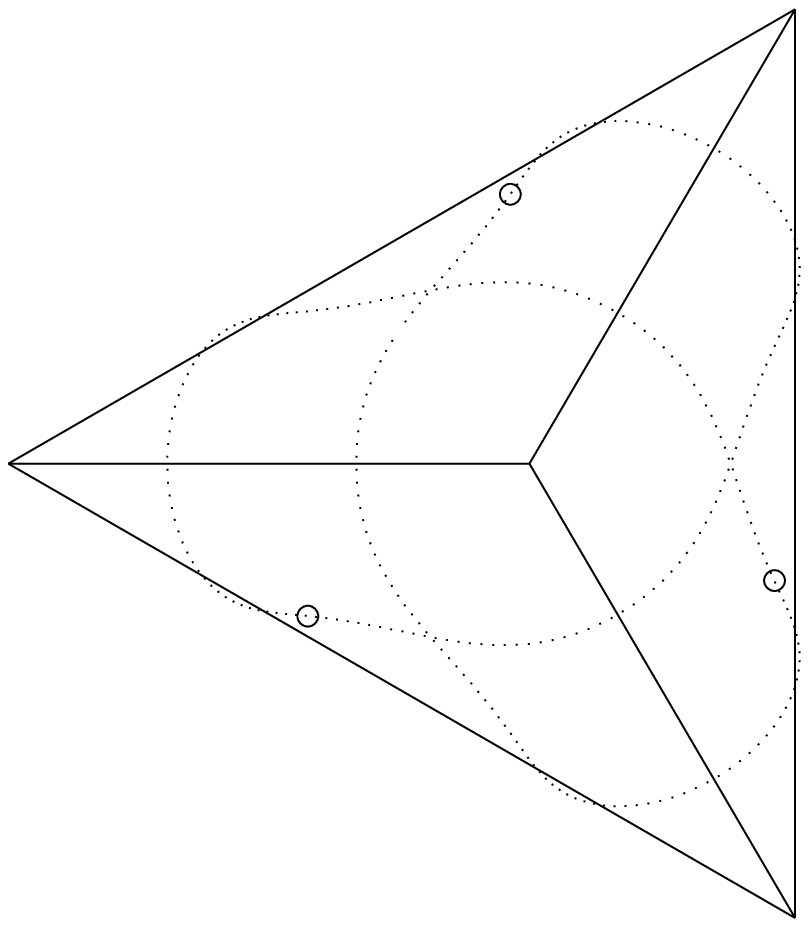,width=3.2cm}
}
\caption{Typical loops in $\K^\nu$ corresponding to $\nu^1$, $\nu^2$, $\nu^3$
for ${\cal R} = {\cal T}$.}
\label{fig:sketch_tetra}
\end{figure}
\begin{figure}[ht]
\centerline{
$\nu^1$ \epsfig{figure=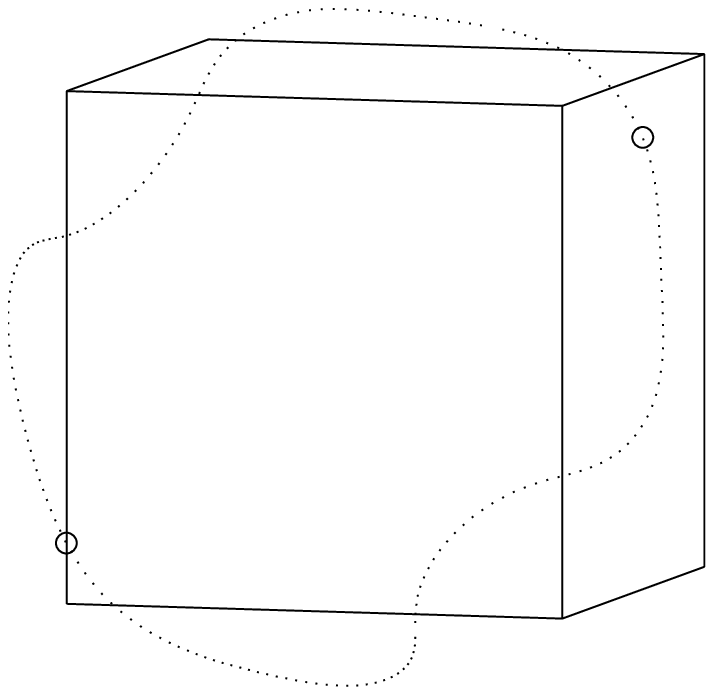,width=3.2cm}
\hskip 0.3cm
$\nu^5$ \epsfig{figure=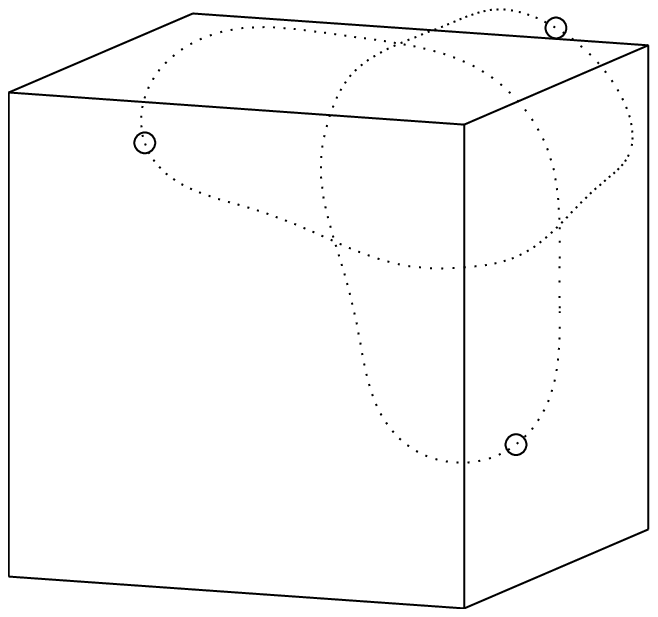,width=3.1cm}
\hskip 0.3cm
$\nu^6$ \epsfig{figure=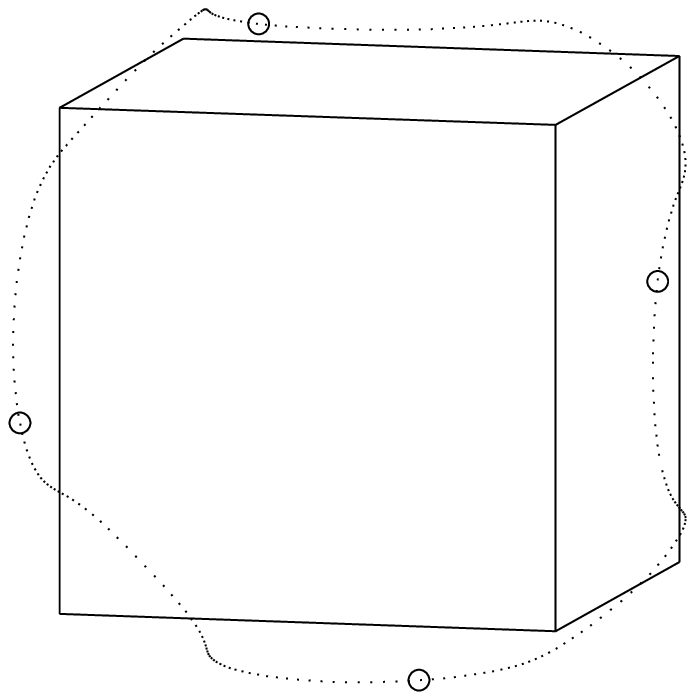,width=3.1cm}
}
\caption{Typical loops in $\K^\nu$ corresponding to $\nu^1$, $\nu^5$, $\nu^6$
for ${\cal R} = {\cal O}$.}
\label{fig:sketch_hexa}
\end{figure}
\begin{figure}[ht]
\centerline{
$\nu^1$ \epsfig{figure=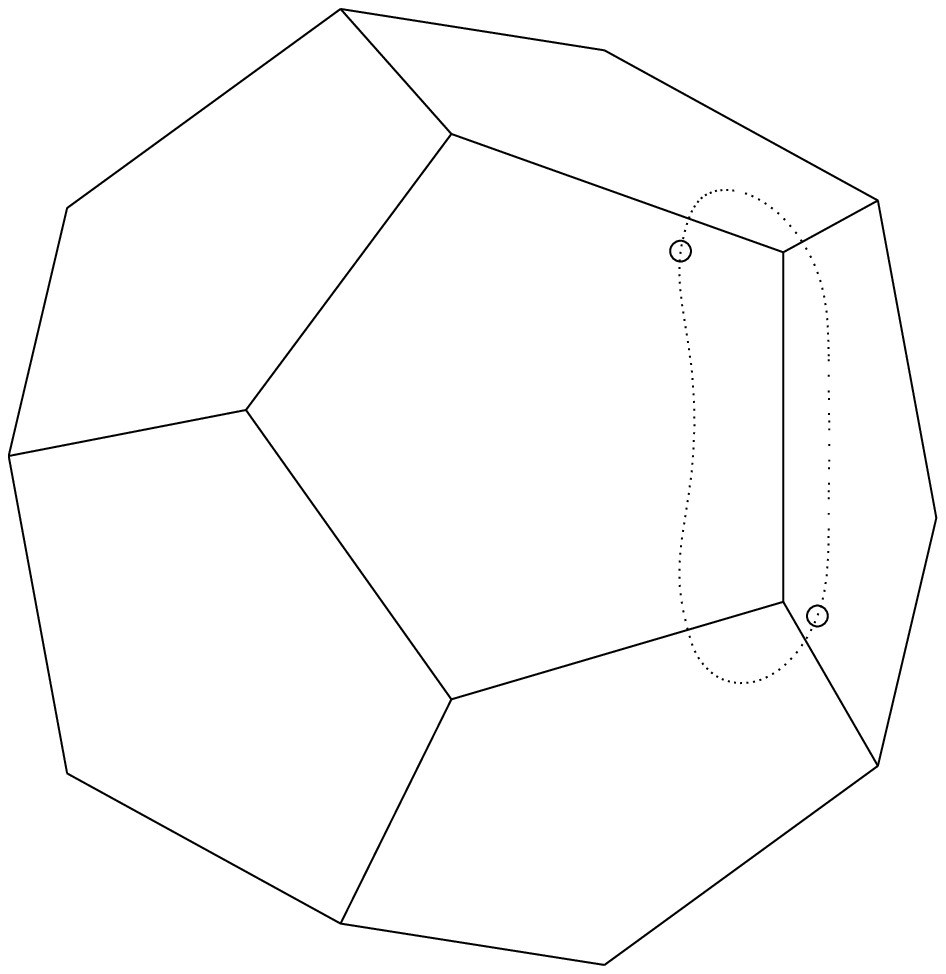,width=3.2cm}
\hskip 0.3cm
$\nu^2$ \epsfig{figure=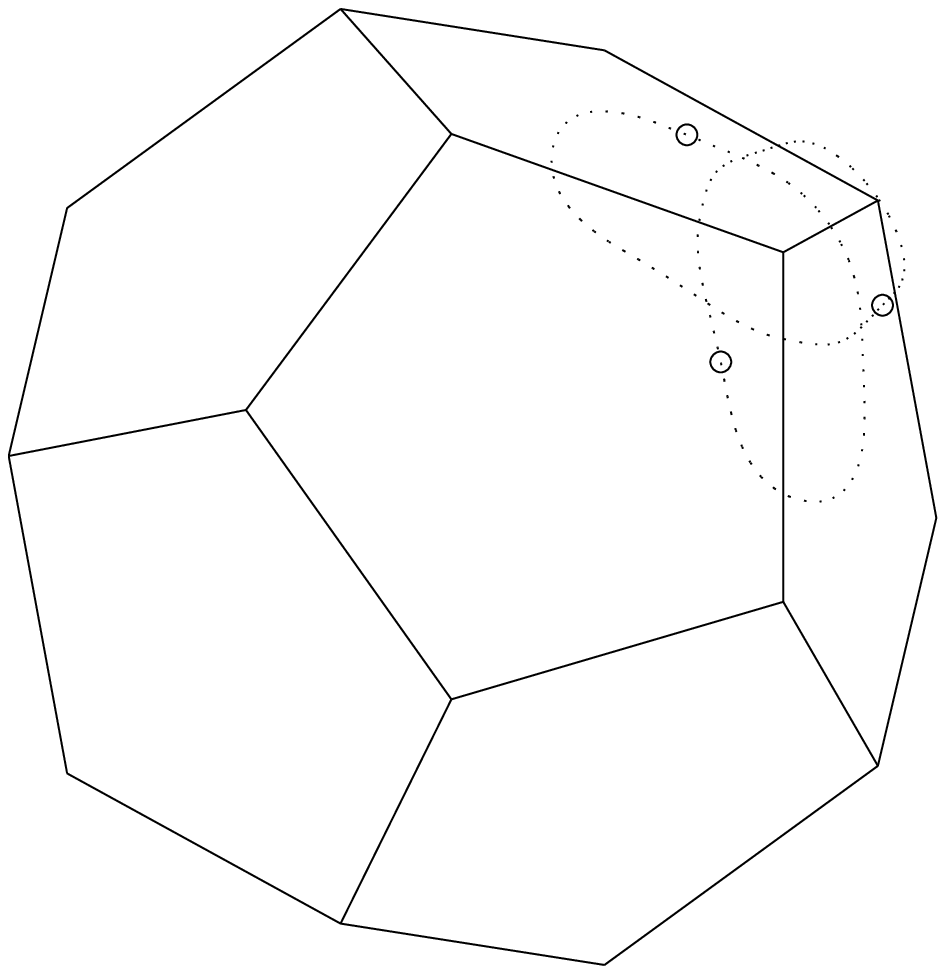,width=3.2cm}
\hskip 0.3cm
$\nu^3$ \epsfig{figure=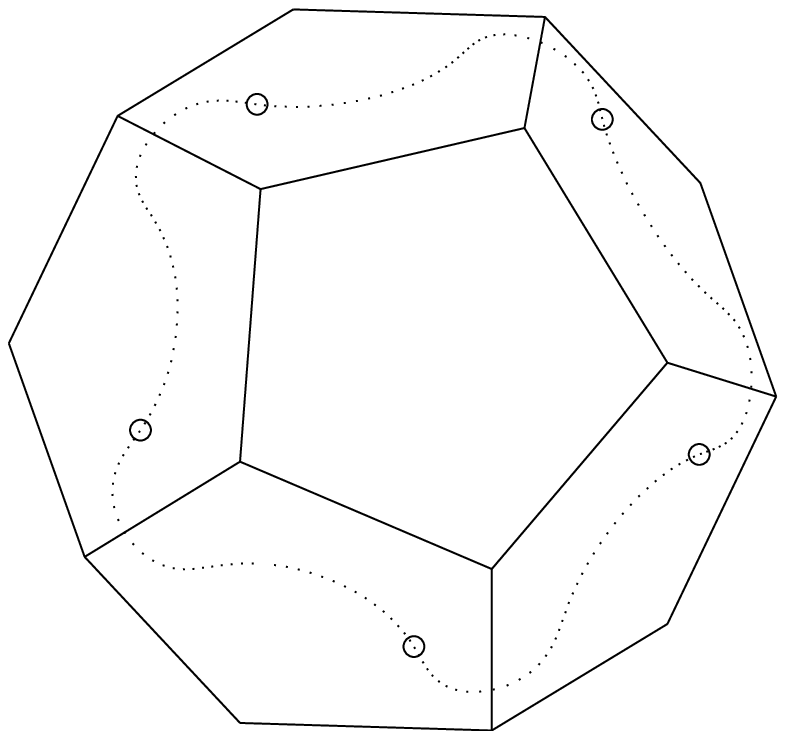,width=3.4cm}
}
\caption{Typical loops in $\K^\nu$ corresponding to $\nu^1$, $\nu^2$, $\nu^3$
  for ${\cal R} = {\cal I}$.}
\label{fig:sketch_dodeca}
\end{figure}

For later reference we observe that if $\nu$ is one of the particular
sequences listed in Theorem~\ref{platorb_exist}, then we can associate
to the corresponding cone 
\[
\K^\nu = \K(\vnuuno)
\] 
(where $\vnuuno$ is the map defined in (\ref{vnun})) a plane $\Pi$ with the
associated reflection $\tilde{R}_\Pi\in \tilde{\cal R}\setminus \{I\}$, a
number $M\in\{2,3,4,5\}$ and a rotation $R\in{\cal R}\setminus\{I\}$ of angle
$2\pi/M$ such that the symmetry conditions
\begin{equation}
\left\{
\begin{array}{l}
\genp(t) = \tilde{R}_\Pi \genp(-t)\cr
\genp(t+T/M) = R\genp(t)\cr
\end{array}
\right.
\label{symm_cond}
\end{equation}
are compatible with membership in $\K^\nu$. This is a straightforward
consequence of the fact that the map $\vnuunogen$ itself satisfies
these conditions for suitable $\tilde{R}_\Pi, M, R$. For each $\nu$ the
corresponding value of $M$ is given in column 3 of Table~\ref{tab:Ktildeni}.
We denote by $\tilde{\K}^\nu \subset\K^\nu$ the subset of $\K^\nu$ of the maps
that satisfy (\ref{symm_cond}). 
Figures~\ref{fig:sketch_tetra}, \ref{fig:sketch_hexa},
\ref{fig:sketch_dodeca} illustrate the structure of the map $\genp$
for a typical element $\typloop\in\tilde{\K}^\nu$.

\noindent In preparation for the proof of
Theorem~\ref{platorb_exist} we note the following
\begin{proposition}
Given $\typloop\in\lamoa$ there exist $\hat{u}\sim u$, $\delta_u>0$ and
sequences $\{t^k\}, \{S_k\}, \{D_k\}$, $k\in\Z$ such that
\begin{itemize}
\item[(i)] $\hatgenp(t^k) \in S_k\,,$
\item[(ii)]$\hatgenp((t^k, t^k+\delta_u))\subset D_k\,,\ \hatgenp([t^k,
t^{k+1}])\subset \overline{D_k}\setminus\Gamma\,,$
\item[(iii)] 
\begin{equation}
S_{k+1}\neq S_k\,,
\label{Skp1Sk}
\end{equation}
\item[(iv)] $\A(\tloophat) \leq \A(\typloop)$\ .
\end{itemize}
\label{prop:action_decrease}
\end{proposition}

\begin{proof}
We identify $\tloophat$ with the map constructed in the proof of
Proposition~\ref{first_alg_char}. Then (i), (ii), (iii) hold
trivially.  To prove (iv) we observe that, given $S\in\{\cal S\}$ there
is a bijection ${\cal R}\setminus\{I\}\ni R\to R^S\in{\cal
R}\setminus\{I\}$ such that (see Figure~\ref{fig:reflection})
\begin{equation}
\tilde{R}_S R x = R^S\tilde{R}_S x, \ x\in\R^3\ .
\label{eq:RSRx}
\end{equation}
\begin{figure}[!h]
\psfragscanon
\psfrag{S}{$S$}
\psfrag{x}{$x$}
\psfrag{Rx}{$Rx$}
\psfrag{Rtilx}{$\tilde{R}_S x$}
\psfrag{RSRtilx=RtilRx}{$R^S\tilde{R}_S x =\tilde{R}_S R x$}
\psfrag{r(R)}{$r(R)$}
\psfrag{r(RS)}{$r(R^S)$}
\centerline{\epsfig{figure=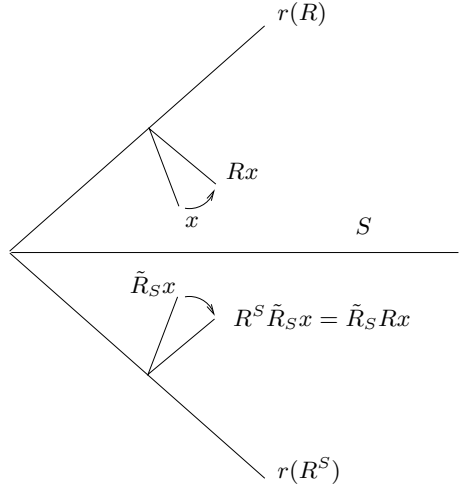,width=6cm}}
\psfragscanoff
\caption{Skematic illustration of $R$, $\tilde{R}_S$, $R^S$ and
equation~(\ref{eq:RSRx})}
\label{fig:reflection}
\end{figure}

\noindent From this it follows
\[
|R^S\tilde{R}_S x - \tilde{R}_S x| = |\tilde{R}_S(Rx - x)| = |Rx - x|
\]
and therefore we have, for every $t\in\R$ and for each $S\in{\cal S}$,
\begin{eqnarray*}
&&\sum_{R\in{\cal R}\setminus\{I\}} \frac{1}{|R \tilde{R}_S\genp(t) - 
\tilde{R}_S\genp(t)|} = 
\sum_{R\in{\cal R}\setminus\{I\}} \frac{1}{|R^S\tilde{R}_S\genp(t) - 
\tilde{R}_S\genp(t)|} = \\
&&=\sum_{R\in{\cal R}\setminus\{I\}} \frac{1}{|R\genp(t) - \genp(t)|}\ .
\end{eqnarray*}
This shows that the reflections used in the proof of
Proposition~\ref{first_alg_char} to construct $\checkgenp$ and then
$\hatgenp$ do not change the potential term of the action. This and
the fact that also the kinetic part of the action is unchanged by a
reflection proves (iv).

\rightline{\small$\square$}
\end{proof}

Before proceeding we observe that, on the basis of the coercivity of $\A|_\K$
proved in Propositions~\ref{conecoerc}, \ref{border&coerc},
\ref{prop:AKu_coerc}, standard arguments from Calculus of Variations
\cite{dacorogna}, \cite{giusti}, \cite{ventPhD} yield the existence of a
minimizer $\minloop\in\accauno\cap\Kclo$, $\Kclo$ the closure of $\K$ in the
$C^0$--topology. We also remark that, in all cases where it can be proved that
$\minloop\not\in\partial\K$, we can invoke the {\em principle of symmetric
criticality} \cite{palais} to deduce that $\minloop$ is actually a critical
point of the unconstrained action functional.

An interesting consequence of Proposition~\ref{prop:action_decrease}
is the following
\begin{theorem}
Given $\typloop\in\lamoa$, assume $\minloop\in\K(\typloop)$ is a
collision free minimizer of $\A|_{\K(\typloop)}$.  Then $\Theta_* :=
\{t\in[0,T): \minloopgen(t)\in S, S\in{\cal S}\}$ is a finite set and
\begin{itemize}
\item[(i)] $\#\Theta_* = n_u K_u$ ($n_u, K_u$ as in
Proposition~\ref{first_alg_char}) and therefore $\#\Theta_*$ is the
minimum compatible with the topological structure of $\typloop$;
\item[(ii)] if $\minloopgen(\bar{t})\in S$ for some $S\in{\cal S}$ and some
$\bar{t}\in[0,T)$, then $\dminloopgen(\bar{t})$ is transversal to $S$.
\end{itemize}
\label{teo:collfreemin}
\end{theorem}
\begin{proof}
Since $\minloop$ is a collision free minimizer, by elliptic regularity it is
a smooth function, therefore, if $\bar{t}\in[0,T)$ is such that
$\minloopgen(\bar{t})\in S$ and $\dminloopgen(\bar{t})$ is not transversal to
$S$, then necessarily $\dminloopgen(\bar{t})$ is parallel to $S$.  This
implies
\begin{equation}
\left\{
\begin{array}{l}
\tilde{R}_S R \minloopgen(\bar t) = R^S\tilde{R}_S \minloopgen(\bar t) =
R^S \minloopgen(\bar t)\cr
\tilde{R}_S R \dminloopgen(\bar t) = R^S\tilde{R}_S
\dminloopgen(\bar t) = R^S \dminloopgen(\bar t)\cr
\end{array}
\right.
\,,
\hskip 0.5cm R\in{\cal R}\ .
\label{RSRjgenp}
\end{equation}

\begin{figure}[h]
\psfragscanon
\psfrag{S}{$S$}
\psfrag{Rjminloop} {$(R\minloopgen(\bar{t}),R\dminloopgen(\bar{t}))$}
\psfrag{RjSminloop}{$(R^S\minloopgen(\bar{t}),R^S\dminloopgen(\bar{t}))$}
\psfrag{minloop}{$(\minloopgen(\bar{t}),\dminloopgen(\bar{t}))$}
\centerline{\epsfig{figure=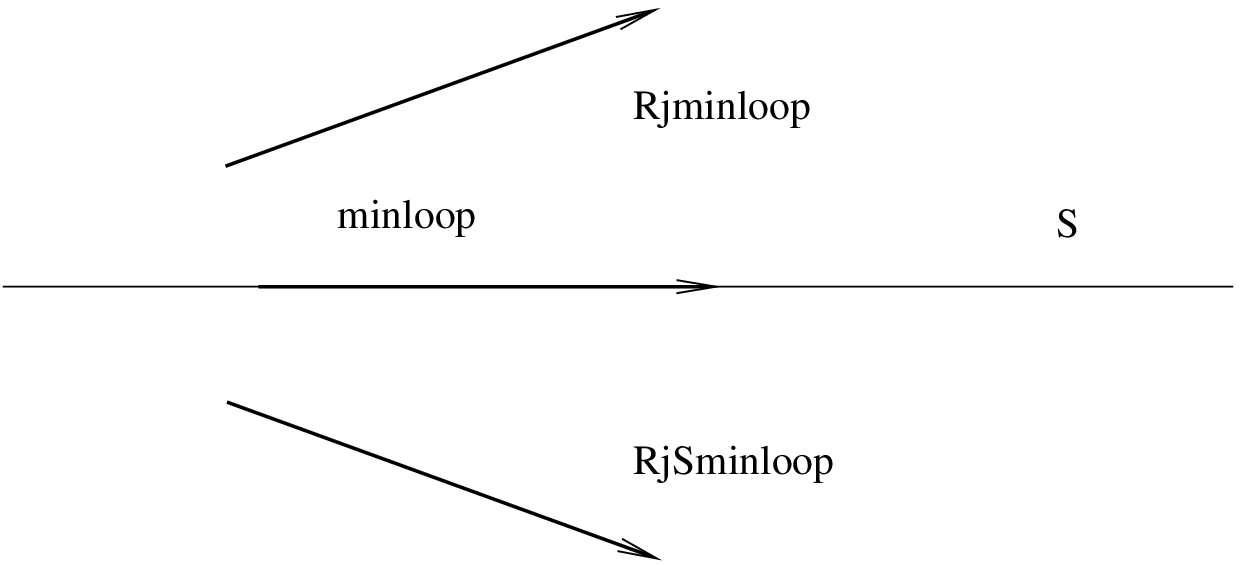,width=8cm}}
\psfragscanoff
\caption{The reflection $\tilde{R}_S$ changes the set
  $\{(R\minloopgen(\bar{t}),R\dminloopgen(\bar{t})), R\in{\cal R}\}$
  into itself.}
\label{fig:reflect}
\end{figure}

Equations (\ref{RSRjgenp}) say that at time $t=\bar t$ the set of
positions and velocities of the $N=|{\cal R}|$ particles of the system
is changed into itself by the reflection $\tilde{R}_S$. From the
symmetry of the equations of motion it follows that the same is true
for all $t$. In particular this implies that $\minloopgen(t)$ belongs
to the plane of $S$ for all $t$.  This is clearly incompatible with
membership in $\K(\typloop)$. Therefore (ii) is established.  Let
$\hatminloop$ be the map associated to $\minloop$ as in
Proposition~\ref{prop:action_decrease}, then $\hatminloop$ is a
minimizer and if (i) does not hold then $\hatminloop\neq\minloop$.  In
particular there is a $\bar t\in\R$ and $S\in{\cal S}$ such
that $\hatminloopgen(\bar t) = \minloopgen(\bar t) \in S$
and
\[
\dhatminloopgen (\bar t^-) = \dminloopgen(\bar t) \neq
\tilde{R}_S\dminloopgen(\bar t) = \dhatminloopgen (\bar t^+)
\]
since by (ii) $\dminloopgen(\bar t)$ is transversal to $S$.
This is in contradiction with (ii). 

\rightline{\small$\square$}
\end{proof}

\begin{remark}
Let $\K\in\lamoasim$ and $(\sigma, n)$ be the corresponding pair.  An
important consequence of Theorem~\ref{teo:collfreemin} is that in the search
of classical $T$-- periodic solutions, minimizing on $\K$ is equivalent to
minimize on the subset of $\K$ of the loops $u$ such that the map $t \to
\genp(t)$ visits periodically one after the other all the chambers $D_k$ in
the sequence $\{D_1,\dots,D_{n K_\sigma}\}$ entering $D_k$ from $S_k$ and
exiting $D_k$ from $S_{k+1}$ without touching the third face
$S^k\not\in\{S_k,S_{k+1}\}$ of $D_k$.
\label{rem:min_num_chambers}
\end{remark}
%

\begin{figure}[ht]
\psfragscanon
\psfrag{rk-1}{$r_{k-1}$}
\psfrag{rk}{$r_k$}
\psfrag{rk+1}{$r_{k+1}$}
\psfrag{rk+2}{$r_{k+2}$}
\psfrag{Sk-1}{$S_{k-1}$}
\psfrag{Sk}{$S_k$}
\psfrag{Sk+1}{$S_{k+1}$}
\psfrag{Sk+2}{$S_{k+2}$}
\psfrag{Sk+3}{$S_{k+3}$}
\psfrag{S^k-1}{$S^{k-1}$}
\psfrag{S^k}{$S^k$}
\psfrag{S^k+1}{$S^{k+1}$}
\psfrag{S^k+2}{$S^{k+2}$}
\centerline{\epsfig{figure=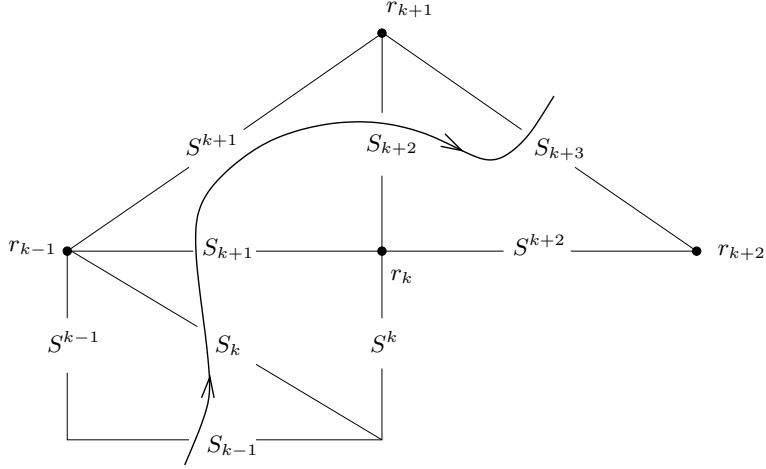,width=10cm}}
\psfragscanoff
\caption{The sequence $k\to r_k = (\overline{S_{k+1}} \cap \overline{S^k})
\setminus\{0\}$.}
\label{fig:rk_seq}
\end{figure}

Remark~\ref{rem:min_num_chambers} has important consequences on the kind of
partial collisions that have to be excluded in the proof that a minimizer
$\minloop\in\overline{\K}$ is collision free.  First of all, if a minimizer
$\minloop$ has a partial collision at time $t_c$, then there is $k$ such that
$\minloopgen(t_c)\neq 0$ belongs to one of the semiaxes on the boundary of
$D^k$. This largely reduces the set of partial collisions that $\minloop\in\K$
may exhibit.  Since each time we add one of the sets $\overline{D_k}$ to the
previous one $\overline{D_{k-1}}$ we introduce the semiaxis
$(\overline{S_{k+1}}\cap\overline{S^k})\setminus\{0\}$; the semiaxes on the
boundary of $D_k$, $k\in\Z$ all appear in the sequence $k\to r_k:=
(\overline{S_{k+1}}\cap \overline{S^k})\setminus\{0\}$.  Therefore the
sequence $k\to r_k$ characterizes the semiaxes where partial collisions can
occur for the particular cone $\K$ under consideration (see
Figure~\ref{fig:rk_seq}).
Other consequences of the preceding discussion on the possible partial
collisions and on their geometric structure will be considered in
Section~\ref{sec:coll} below.
\begin{remark}
Theorem~\ref{teo:collfreemin} and Remark~\ref{rem:min_num_chambers} apply
also to minimizers $\minloop\in\tilde{\K}^\nu$ of $\A|_{\tilde{\K}^\nu}$, with
$\tilde{\K}^\nu\subset\K^\nu$ the cone introduced above
(cfr. (\ref{symm_cond}).
\end{remark}

\section{Collisions}
\label{sec:coll}

In this Section we always denote by $\K$ either $\K_4$, or $\KPi$, or
$\tilde{\K}^\nu$ with $\nu$ as in Theorem~\ref{platorb_exist}. We show
that minimizers $\minloop\in\K$ are collision free.
We start by excluding total collisions. In Section~\ref{sec:part_coll}
we shall discuss the case of partial collisions.

\subsection{Total collisions}
\label{sec:tot_coll}
Our strategy to show that actually in all cases a minimizer
$\minloop\in\overline{\K}$ does not have total collisions is based on
{\em level estimates}, that is we show
\begin{itemize}
\item[(a)] the assumption that $\minloop$ has a total collision implies a
bound of the form
\[
\A(\minloop) \geq a>0\,;
\]
\item[(b)] there exists $v\in\K$ such that 
\[
\A(v) <a\ .
\]
\end{itemize}
This approach is quite natural in our context. Indeed a total
collision implies $\minloop\in \partial\K$ and therefore any attempt
to perturb $\minloop$ into a competing function $v$ such that
$\A(v)<\A(\minloop)$ to show that $\minloop$ is free of total
collisions runs against the difficulty of respecting the topological
constraints that characterize membership in $\K$.

We begin with the cone $\K_4$.
\begin{proposition}
Assume $u\in \overline{\K_4}$
has a total collision. Then
\[
\mathcal{A}( u) \ge a_4 = \frac{18}{2^{\frac{1}{3}}}
\pi^{\frac{2}{3}}T ^{\frac{1}{3}}
\]
\label{prop:A_totcoll}
\end{proposition}
\begin{proof}
From the definition of $\K_4$ it follows that, if $\genp = (u_{11}, u_{12},
u_{13})$,
\begin{eqnarray}
\mathcal{A}(u) &=& \int_0^T \biggl( 2 |\dot u_1|^2 +
\frac{1}{\sqrt{u_{11}^2 + u_{12}^2}} + \frac{1}{\sqrt{u_{11}^2 +
u_{13}^2}} + \frac{1}{\sqrt{u_{12}^2 + u_{13}^2}}\biggr) dt
\nonumber\\
&=& \int_0^T \biggl(\frac{1}{2} [2 (\dot u_{11}^2 + \dot u_{12}^2)]+
\frac{1}{\sqrt{u_{11}^2 + u_{12}^2}} \biggr)dt + \int_0^T
\biggl(\frac{1}{2} [2 (\dot u_{12}^2 + \dot u_{13}^2)]+ \nonumber\\
&+& \frac{1}{\sqrt{u_{12}^2 + u_{13}^2}} \biggr)dt + \int_0^T
\biggl(\frac{1}{2} [2 (\dot u_{11}^2 + \dot u_{13}^2)]+
\frac{1}{\sqrt{u_{11}^2 + u_{13}^2}} \biggr)dt \nonumber\\
&\stackrel{def}{=}& a(u_{11},u_{12},T) + a(u_{12},u_{13},T)
+a(u_{11},u_{13},T)\ .
\label{sum_3_actions}
\end{eqnarray}
If $m_1, m_2, K$ are positive constants such that
\[
\frac{m_1 m_2}{m_1+ m_2} = 2\,,\hskip 1cm Km_1 m_2 = 1\,,
\]
and $\xi:\R\rightarrow \R^2$ is a periodic map of period $\mathtt{T}$,
we can write
\[
a(\xi_1,\xi_2,\mathtt{T}) = \int_0^{\mathtt{T}}\biggl(\frac{1}{2} \frac{m_1
m_2}{m_1+ m_2} |\dot \xi|^2 + \frac{Km_1 m_2}{|\xi|}\biggr)\;dt\ .
\]
Gordon's Theorem~\cite{gordon77} implies that if $\xi(t)$ vanishes at some $t
\in [0, \mathtt{T})$, then
\begin{equation}
a(\xi_1,\xi_2,\mathtt{T}) \geq 3\biggl(\frac{K^2 \pi^2}{2(m_1 +
m_2)}\biggr)^{\frac{1}{3}} m_1 m_2 \mathtt{T}^{\frac{1}{3}} = 3
\pi^{\frac{2}{3}} \mathtt{T}^{\frac{1}{3}} \ .
\label{action_twobody}
\end{equation}
If $u\in \overline{\K_4}$ has a total collision then $u(t)$ vanishes
at least at two times $t_1, t_2$, with $t_2-t_1=T/2$, in the interval
$[0, T)$. From this, (\ref{sum_3_actions}) and (\ref{action_twobody})
we get, for $\mathtt{T}=T/2$,
\[
\mathcal{A}( u) \ge 3\cdot 2\cdot 3 \pi^{\frac{2}{3}}
\mathtt{T}^{\frac{1}{3}}= \frac{18}{2^{\frac{1}{3}}}
\pi^{\frac{2}{3}}T ^{\frac{1}{3}}\ .
\]

\rightline{\small$\square$}
\end{proof}

\begin{proposition}
There exists $\sloop\in\K_4$ such that
\[
\A(v)< \frac{18}{3^{\frac{1}{3}}} \pi^{\frac{2}{3}}T ^{\frac{1}{3}}\ .
\]
\label{prop:nototcoll_4b}
\end{proposition}
\begin{proof}
 We define $\sloop$ by describing the motion $\sloopgen$ of $\partgen$, the
generating particle. $\partgen$ moves with constant speed on a closed curve
which is the union of four half circumferences $C_1^\pm, C_2^\pm$ of radius
$\rho > 0$. $C_1^\pm$ has center on the axis $\assetre$ and lies on the plane
$\assetre = \pm \rho$. $C_2^\pm$ has center on the axis $\assedue$ and lies on
the plane $\assedue = \pm \rho$.  The kinetic part $A_K(\sloop)$ of
$\A(\sloop)$ is given by
\[
A_K(\sloop) = 4 \cdot \frac{T}{2} {\left(\frac{4 \pi \rho}{T}\right)}^2 =
32\frac{\pi^2\rho^2}{T}\ .
\]
From the definition of $\sloop$ it follows that $|\sloop_i - \sloop_j|\ge
2\rho$ whenever $i\ne j$. Then the potential part $A_U(\sloop)$ of
$\A(\sloop)$ satisfies
\[
A_U(\sloop) =\int_0^T \sum_{i< j}\frac{1}{|\sloop_i - \sloop_j|} <
6\frac{T}{2\rho}\,,
\]
therefore we have
\begin{equation}
\A(\sloop) <32\frac{\pi^2\rho^2}{T} + 3\frac{T}{\rho}\,, \qquad \forall \rho >
0\ .
\label{stima_per_A}
\end{equation}
If we choose $\rho =  {\left(\frac{3 T^2}{64\pi^2}\right)}^{\frac{1}{3}}$, to
minimize the r.h.s. of (\ref{stima_per_A}), we obtain
\[
\A(\sloop) < \frac{18}{3^{\frac{1}{3}}}\pi^{\frac{2}{3}} T^{\frac{1}{3}}\ .
\]

\rightline{\small$\square$}
\end{proof}
Propositions (\ref{prop:A_totcoll}) and (\ref{prop:nototcoll_4b})
imply that a minimizer $\minloop\in \Kclo$ can not have total
collisions.  
We now consider the cases $\K=\KPi$ or $\K=\tilde{\K}^\nu$.

\begin{proposition}
Let $\Omega\subset\R^s$ be an open connected set, $L_k:\Omega\times
\R^s\rightarrow \R$, $0\leq k\leq M$ be smooth Lagrangian functions and
$Q_k\subset \{u: u\in H^1_T(\R, \R^s), u(\R)\subset\Omega\}$.  Let
$\A_k:Q_k\to\R$ be the action functional
\begin{displaymath}
\A_k(q)=\int_0^T L_k(q,\dot q)\;dt\,, \hskip 1cm 0\leq k\leq
M\ .
\end{displaymath}
Assume
\begin{itemize}
\item[{\rm (i)}]
$\sum_{k=1}^M L_k(x,y)\leq L_0(x,y), \hskip 0.2cm\forall(x,y)\in
\Omega\times \R^s$,

\item[{\rm (ii)}] $Q_0\subset Q_k, \hskip 0.2cm 1\leq k\leq M$.
\end{itemize}
Then
\begin{displaymath}
\sum_{k=1}^M \inf_{q\in Q_k} \A_k(q)\leq \inf_{q\in Q_0} \A_0(q).
\end{displaymath}
\label{prop4.5}
\end{proposition}

\begin{proof} (i) and (ii) imply 
\begin{displaymath}
\sum_{k=1}^M \inf_{q\in Q_k} \A_k(q) \leq
\sum_{k=1}^M \inf_{q\in Q_0} \A_k(q) \leq
\inf_{q\in Q_0}\sum_{k=1}^M \A_k(q) \leq
\inf_{q\in Q_0} \A_0(q)\ .
\end{displaymath}

\rightline{\small$\square$}
\end{proof}

\begin{proposition}
Assume $u\in\lamoa$ has a total collision. Then
\begin{equation}
\A(u) \geq 3
N\left(\frac{\pi^2(N-1)^2}{32}\right)^{1/3}T^{1/3}\,,\qquad\mbox{with
}N=|{\cal R}|\ .
\label{lowerbound}
\end{equation}
\label{prop:lower}
\end{proposition}
\begin{proof}
To prove (\ref{lowerbound}) we apply Proposition~\ref{prop4.5} with
$M=1$ and
\[
L_0(x,y)=\frac{N}{2}\left(|y|^2 + \sum_{R\in{\cal R}\setminus\{I\}}
\frac{1}{|(R-I)x|}\right),
\]
\[
Q_0 = \{\genp\in H^1_T(\R,\R^3\setminus\Gamma) : \lim_{t\to\{0,T\}}
\genp(t) = 0\},
\]
\[
L_1(x,y)=\frac{N}{2}\left(\frac{(y\cdot x)^2}{x\cdot x} + \frac{(N-1)}{2|x|}
\right),
\]
\[
Q_1 = \{\genp \in H^1_T(\R,\R^3\setminus\{0\}) : \lim_{t\to\{0,T\}}
\genp(t) = 0\}\ .
\]
Set $\rho=|\genp|$. Then we have
\begin{equation}
\A_1(u) = \frac{N}{2}\int_0^T\left( {\dot{\rho}}^2 +
\frac{N-1}{2\rho}\right)\,dt \stackrel{def}{=} N a_1(\rho)\ .
\label{actionuno}
\end{equation}
Therefore minimizing $\A_1$ on $Q_1$ is equivalent to the minimization of
$a_1$ on $H^1_T(\R, (0,+\infty))$. From (\ref{actionuno}) $a_1$ is the
action of a system of two masses $m_1=m_2=1/2$ interacting with Newtonian
potential of gravitational constant $2(N-1)$.  From this and Gordon's
Theorem it follows
\[
a_1(\rho) \geq 3\left(\frac{\pi^2 (N-1)^2}{32}\right)^{1/3}\, T^{1/3}\ .
\]
This and (\ref{actionuno}) imply (\ref{lowerbound}).

\rightline{\small$\square$}
\end{proof}

The estimate (\ref{lowerbound}) in Proposition~\ref{prop:lower} is
based on the assumption, fulfilled by any $\typloop\in\lamoa$, that at each
time all the particles have the same distance $\rho$ from the origin $O$, and
on the obvious bound $|\typloop_i-\typloop_j|\leq 2\rho$, which implies
\[
\sum_{R\in{\cal R}\setminus\{I\}}\frac{1}{|(R-I)\genp|} 
\geq (N-1)\frac{1}{2\rho}\ .
\]

\small\noindent Simple geometric observations allow for sharper estimates. For
instance for ${\cal R}={\cal O}$ one can observe that at each time $t$ the
generating particle $\partgen$ is the vertex of 3 squares and 4 equilateral
triangles, the other vertexes of which are all particles distinct from each
other and from $\partgen$. Since the maximum of the side of a square with
vertexes on a sphere with radius $\rho$ is $\sqrt{2}\rho$ and the maximum of
the side of an equilateral triangle with vertexes on a sphere with radius
$\rho$ is $\sqrt{3}\rho$, we conclude that at each time $t$ there are 6
particles at distance $\leq\sqrt{2}\rho$ from $\partgen$ and 8 particles at
distance $\leq \sqrt{3}\rho$ from $\partgen$. From these observations it
follows that in the case ${\cal R}={\cal O}$ we have
\begin{equation}
\sum_{R\in{\cal R}\setminus\{I\}}\frac{1}{|(R-I)\genp|} 
\geq \left(\frac{9}{2} +
\frac{6}{\sqrt{2}} + \frac{8}{\sqrt{3}} \right)\frac{1}{\rho} > 1.161\cdot
\frac{23}{2\rho} \ .
\label{estforC}
\end{equation}
Similarly, if ${\cal R}={\cal T}$, $\partgen$ is the vertex of 4 equilateral
triangles, therefore we have
\begin{equation}
\sum_{R\in{\cal R}\setminus\{I\}}\frac{1}{|(R-I)\genp|} 
\geq \left(\frac{3}{2} +
\frac{8}{\sqrt{3}} \right)\frac{1}{\rho} > 1.112\cdot \frac{11}{2\rho} \ .
\label{estforT}
\end{equation}
Finally, if ${\cal R}={\cal I}$, $\partgen$ is a vertex of 6 regular pentagons
and 10 equilateral triangles. For each pentagon there are 2 particles whose
distance from $\partgen$ is equal to the side of the pentagon and 2 particles
at distance equal to $2\cos(\pi/10)$ times the radius of the
circumcircle. Therefore we have the estimate
\begin{equation}
\sum_{R\in{\cal R}\setminus\{I\}}\frac{1}{|(R-I)\genp|} \geq
\left(\frac{15}{2} + \frac{12}{2\sin(\pi/5)} +
\frac{12}{2\cos(\pi/10)} + \frac{20}{\sqrt{3}}\right)\frac{1}{\rho} >
1.205\cdot\frac{59}{2\rho}\ .
\label{estforI}
\end{equation}

\normalsize
We denote by $a_{\cal R}$ the right--hand side of (\ref{lowerbound})
and by $a'_{\cal R}$ the analogous lower bound obtained by using the
improved estimates (\ref{estforC}), (\ref{estforT}), (\ref{estforI}).
In Table~\ref{tab:aR} we list the values of $a_{\cal R}, a'_{\cal R}$
(approximated by truncation).
\begin{table}[h]
\caption{Lower bounds for the action in case of total collision.}
\label{tab:aR}
\begin{center}
\begin{tabular}{clclc}
\hline\noalign{\smallskip}
${\cal R}$ &\vline &$a_{\cal R}/T^{1/3}$ &\vline &$a'_{\cal R}/T^{1/3}$ \\ 
\hline
${\cal T}$ &\vline &\ 120.3042  &\vline &\ 129.1665  \\
\hline
${\cal C}$ &\vline &\ 393.4301  &\vline &\ 434.8151  \\
\hline
${\cal I}$ &\vline &1843.1348   &\vline &2087.7547 \\
\hline
\end{tabular}
\end{center}
\end{table}
The lower bounds given in Table~\ref{tab:aR} applies to any
$\typloop\in\overline{\K}$, which is known to have one total collision
per period. If $\typloop$ has $M>1$ total collisions per period and the
time intervals between subsequent collisions are all equal, then from
(\ref{lowerbound}) we derive
\begin{equation}
\A(\typloop) \geq M a_{\cal R}\frac{1}{M^{1/3}} = M^{2/3} a_{\cal R}
\label{lb_with_aR}
\end{equation}
and the same is true with $a'_{\cal R}$ in place of $a_{\cal R}$.  This
observation applies in particular to the cones $\KPi$. For $\typloop\in\KPi$,
the definition of $\KPi$ (cfr. $\condA$, $\condB$, $\condC$) implies
\[
|\typloop(t+T/H)| = |\typloop(t)|, \ \ \forall t\in\R\,,
\]
therefore if $\typloop$ has a total collision at time $t_c$ it also has a
total collision at time $t_c+T/H$ and we can apply (\ref{lb_with_aR}) with
$M=H$. In Table~\ref{tab:KPi} (lines 2, 3) we list the values of
$H^{2/3}a_{\cal R}/T^{1/3}, H^{2/3}a'_{\cal R}/T^{1/3}$.

\begin{table}[t]
\caption{Values of the action for $\typloop\in\KPi$.}
\label{tab:KPi}
\begin{center}
\begin{tabular}{c|c|c|c|c|c}
\hline\noalign{\smallskip}
  &$\Tplato$ &$\Cplato$ &$\Oplato$ &$\Dplato$ &$\Iplato$\\
\hline
& & & & & \\
 $H$ &3 &4 &3 &5 &3  \\
\hline 
& & & & & \\
 $H^{2/3} a_{\cal R}/T^{1/3}$ &$250.2428$ &$991.3818$ &$818.3676$ &$5389.3588$ &$3833.8749$ \\
\hline
& & & & & \\
 $H^{2/3} a'_{\cal R}/T^{1/3}$ &$268.6772$ &$1095.6654$ &$904.4519$ &$6104.6318$ &$4342.7048$ \\
\hline\hline
& & & & & \\
$\A(\sloop)/T^{1/3}, i=1$ &$220.2007$ &$734.9502$ &$589.9526$ &$2866.6116$ &$2027.2544$ \\
\hline
& & & & & \\
$\A(\sloop)/T^{1/3}, i=2$ &$168.0446$ &$553.1633$ &$589.9526$ &$2181.2066$ &$2452.2053$ \\
\hline
& & & & & \\
$\A(\sloop)/T^{1/3}, i=3$ &$266.7542$ &$896.4157$ &$819.8050$ &$3477.7486$ &$3208.5266$ \\
\hline
\end{tabular}
\end{center}
\end{table}
The inequality (\ref{lb_with_aR}) can also be applied to
$\typloop\in\tilde{\K}^\nu$. In Table~\ref{tab:Ktildeni}, for each
$\nu$ considered in Theorem~\ref{platorb_exist} we list the
corresponding values of $M\in\{2,3,4,5\}$ and the lower bounds given
by (\ref{lb_with_aR}) .

\begin{table}[t]
\caption{Values of the action for $u\in\tilde{\K}^\nu$.}
\label{tab:Ktildeni}
\begin{center}
\begin{tabular}{c|c|c|c|c|c}
\hline\noalign{\smallskip}
${\cal R}$ &$\nu$ &$M$ &$M^{2/3}a_{\cal R}/T^{1/3}$ &$M^{2/3}a'_{\cal R}/T^{1/3}$ &$\A(\sloop)/T^{1/3}$\\
\hline
${\cal T}$ &1     &2  &$190.9710$  &$205.0391$ &$168.0445$  \\
           &2     &3  &$250.2428$  &$268.6772$ &$168.0445$  \\
           &3     &3  &$250.2428$  &$268.6772$ &$266.7542$  \\
\hline
${\cal O}$ &1     &2  &$624.5314$  &$690.2260$  &$647.2635$ \\
           &2     &2  &$624.5314$  &$690.2260$  &$553.1632$ \\
           &3     &2  &$624.5314$  &$690.2260$  &$462.9895$ \\
           &4     &2  &$624.5314$  &$690.2260$  &$647.2635$ \\
           &5     &3  &$818.3676$  &$904.4519$  &$724.8489$ \\
           &6     &4  &$991.3818$  &$1095.6654$  &$859.5748$ \\
\hline
${\cal I}$&1  &2 &$2925.7941$  &$3314.1040$  &$1556.2362$ \\
          &2  &3 &$3833.8749$  &$4342.7048$  &$2463.1128$\\
          &3  &5 &$5389.3588$  &$6104.6318$  &$3447.1168$ \\
\end{tabular}
\end{center}
\end{table}

\begin{table}[t]
\caption{Minimal sequences for the cones $\KPi$ considered in
Theorem~\ref{teoplato}.}
\label{tab:KPi_seq}
\begin{center}
\begin{tabular}{c|c|l}
\hline\noalign{\smallskip}
$P$       &$\KPi$ &\hskip 2cm$\nu$\\
\hline
          &1 &$[2, 7, 1, 9, 12, 8, 3, 10, 5, 2]$ \\
$\Tplato$ &2 &$[2, 10, 3, 12, 9, 7, 2]$ \\
          &3 &$[2, 10, 5, 8, 3, 12, 8, 1, 9, 7, 1, 5, 2]$\\
\hline
          &1 &$[5, 1, 16, 10, 3, 8, 18, 7, 20, 23, 14, 11, 5]$\\
$\Cplato$ &2 &$[5, 16, 10, 8, 18, 20, 23, 11, 5]$\\
          &3 &$[1, 5, 16, 1, 3, 10, 8, 3, 7, 18, 20, 7, 14, 23, 11, 14, 1]$\\
\hline
          &1 &$[3, 1, 16, 10, 6, 15, 8, 18, 7, 3]$\\
$\Oplato$ &2 &$[1, 16, 22, 6, 15, 13, 18, 7, 14, 1]$\\
          &3 &$[3, 10, 16, 22, 6, 10, 8, 15, 13, 18, 8, 3, 7, 14, 1, 3]$\\
\hline
          &1 &$[1, 54, 59, 3, 7, 47, 6, 15, 48, 11, 28, 45, 19, 43, 50, 1]$\\
$\Dplato$ &2 &$[54, 59, 7, 47, 15, 48, 28, 45, 43, 50, 54]$\\
          &3 &$[54, 1, 3, 59, 7, 3, 6, 47, 15, 6, 11, 48, 28, 11, 19, 45, 
               43, 19, 1, 50, 54]$\\
\hline
          &1 &$[28, 45, 19, 11, 6, 15, 48, 34, 42, 28]$ \\
$\Iplato$ &2 &$[45, 19, 1, 3, 6, 15, 25, 38, 34, 42, 20, 31, 45]$\\
          &3 &$[45, 28, 11, 19, 1, 3, 6, 11, 48, 15, 25, 38, 34, 48, 28, 42, 20, 31, 45]$\\
\hline
\end{tabular}
\end{center}
\end{table}

To complete the proof that minimizers $\minloop\in\overline{\K},
\K=\KPi$ or $\K=\tilde{\K}^\nu$, $\nu$ as in
Theorem~\ref{platorb_exist}, are free of total collisions we now show
that in all these cases there exists $\sloop\in\K$ which is
collision free and has a value $\A(\sloop)$ of the action below the
lower bounds discussed above. For $\K=\tilde{\K}^\nu$ we choose
$\sloop=\lambda\vnuuno$, where $\vnuuno$ is defined by (\ref{vnun})
and $\lambda>0$ will be chosen later. For the cones $\KPi$ the
sequence $\nu$ is not uniquely determined. On the other hand, if
$\typloop\in\KPi$ then the map $\genp$ must necessarily visit, and in
a well determined order, certain domains $D\in{\cal D}$. This determines a
minimal sequence $\nu$ compatible with membership in $\KPi$.  This
minimal sequence $\nu$ is the one we use for defining our test
function $\sloop = \lambda\vnuuno \in\KPi$. With reference to the
numbering of the vertexes of $\mathcal{Q}_{\cal R}$ in
Figure~\ref{fig:archimede}, we list in Table~\ref{tab:KPi_seq} the
minimal $\nu$ corresponding to each $\KPi$.

Let $A_K = A_K^{(\nu,n)}, A_U = A_U^{(\nu,n)}$, the kinetic and the
potential part of the action:
\begin{equation}
A_K = \frac{N}{2}\int_0^T |\dvnungen|^2\;dt\,,
\hskip 0.5cm
A_U = \frac{N}{2}\int_0^T \sum_{R\in{\cal R}\setminus\{I\}}
\frac{1}{|(R-I)\vnungen|}\;dt\ .
\end{equation}
We choose $\lambda = (\frac{A_U}{2A_K})^{1/3}$, that gives to
$\A(\lambda\vnun) = \lambda^2 A_K + \frac{1}{\lambda}A_U$
its minimum value
\begin{equation}
\A(\sloop) = 3\Bigl(\frac{A_K A_U^2}{4}\Bigr)^{1/3}\ .
\label{Aminval_lp}
\end{equation}
The main reason for considering test functions of the form $\sloop = \lambda
\vnun$ is that the piecewise affine character of $\vnun$
implies that $A_U$ is the sum of $N-1$ elementary integrals and by consequence
$\A(\sloop)$ has an explicit analytic expression. Moreover, on the basis of
simple observations, the computation of $A_U$, and in turn the computation of
$\A(\sloop)$, can be reduced to a purely algebraic fact.  Going back to the
construction of the polyhedron $\mathcal{Q}_{\cal R}$ outlined in
Section~\ref{sec:platocones_two} we set $q_i = \tilde{R}_i q, i=1,2$
(cfr. Section~\ref{sec:platocones_two} for the definition of $q$ and
$\tilde{R}_i$) and we deduce from the discussion in
Section~\ref{sec:platocones_two} that the set ${\cal L}_{\cal R}$ of the sides
of $\mathcal{Q}_{\cal R}$ is the union of the two orbits $\{R[q,
q_i]\}_{R\in{\cal R}}, i=1,2$, of the segments $[q, q_i], i=1,2$. It follows
that we can associate to each $j\in\Z$ a uniquely determined pair
$(R_j,i_j)\in{\cal R}\times\{1,2\}$ such that $[\nu_{j-1}, \nu_j] = R_j[q,
q_{i_j}]$. For each given $R'\in{\cal R}$ set
\[
\upsilon_i(R') = \int_0^1\sum_{R\in{\cal R}\setminus\{I\}}
\frac{ds}{\bigl|(R-I)R'[(1-s)q+s q_i]\bigr|}\,, \qquad i=1,2\ .
\]
Since
\[
\bigl|(R-I)R'[(1-s)q+s q_i]\bigr| = \bigl|((R')^{-1}RR' -I)[(1-s)q+s
q_i]\bigr|
\]
and the map $R\mapsto (R')^{-1}RR'$ is an isomorphism of ${\cal R}$ onto
itself, we have
\[
\upsilon_i(R') = \upsilon_i(I) \stackrel{def}{=} \upsilon_i\ .
\]
From this and the fact that $\sloopgen^{(\nu,n)}$ travels each side
$[\nu_{j-1}, \nu_j]$ in a time interval of size $T/(nK_\nu)$, it
follows
\begin{eqnarray}
A_U &=& \frac{N}{2}\frac{T}{n K_\nu}\sum_{j=1}^{nK_\nu} \int_0^1\sum_{R\in{\cal
R}\setminus\{I\}} \frac{ds}{\bigl|(R-I)R_j[(1-s)q+s q_{i_j}]\bigr|} =\nonumber
\\
&=&\frac{N}{2}\frac{T}{nK_\nu}(N_1\upsilon_1 + N_2\upsilon_2)\,,
\label{action_pot}
\end{eqnarray}
where $N_i$ is the number of sides $[\nu_{j-1},\nu_j]$ in the orbit of
$[q, q_i]$.

There is a simple geometric criterium to decide whether $[\nu_{j-1},\nu_j]$ is
in the orbit of $[q, q_i]$.  Besides the $|\tilde{\cal R}|/4$ squares, the
other faces of $\mathcal{Q}_{\cal R}$ are the images under some $R\in{\cal R}$
of the polygons ${\cal F}_i$ with vertexes $\{Rq\}_{R\in C_i}, i=1,2$, where
$C_1, C_2$ are the cyclic groups of the rotations with axis $\asseuno$ and
$OV$ respectively.  The side $[\nu_{j-1}, \nu_j]$ is in the orbit of $[q,
q_i]$ if it is one of the sides of $R{\cal F}_i$ for some $R\in{\cal R}$.

In the case ${\cal R}={\cal T}$, ${\cal F}_1$ and ${\cal F}_2$ are
both equilateral triangles and it is straightforward to check that we
have $\upsilon_1 = \upsilon_2 \stackrel{def}{=}
\upsilon$, and (\ref{action_pot}) becomes
\footnote{We write the 
explicit expression of $\upsilon$ for the case ${\cal R}={\cal T}$:
\begin{eqnarray*}
\upsilon &=& -\ln \left( \sqrt {2}-1 \right) \sqrt {2}-2\,\ln \left(
 2-\sqrt {3} \right) -2\,\ln \left( 3 \right) + \\
&+& 2/3\,\sqrt {3}\ln \left( 3 \right) +2\,\ln \left( 2+\sqrt {3}
 \right) -2/3\,\ln \left( 2- \sqrt {3} \right) \sqrt {3}+ \\
&+& 2/3\,\ln
 \left( 2+\sqrt {3} \right) \sqrt {3}-\ln \left( \sqrt {2}-1 \right)\ .
\end{eqnarray*} 
We omit the analogous, but longer, expressions of $\upsilon_1,
\upsilon_2$ for ${\cal R} = {\cal O}, {\cal I}$.}
\[
A_U = 6T\upsilon\ .
\]
From (\ref{Aminval_lp}), (\ref{action_pot}) and
\[
A_K = \frac{N}{2}\frac{\ell^2n^2K_\nu^2}{T}\,,
\]
where $\ell$ is the length of a side of $\mathcal{Q}_{\cal R}$, we
finally obtain
\begin{equation}
\A(\sloop) = \frac{3}{2\cdot 4^{1/3}} N \ell^{2/3}(N_1\upsilon_1 +
N_2\upsilon_2)^{2/3} T^{1/3}\ .
\label{Asloop_val}
\end{equation}
In the last column of Table~\ref{tab:Ktildeni} we list the values of
$\A(\sloop)$ given by (\ref{Asloop_val}) for the cases considered in
Theorem~\ref{platorb_exist}. We recall that the function $\vnungen$, defined
in (\ref{vnun}), will automatically satisfy condition (\ref{symm_cond}). That
is $\vnun\in\tilde{\K}^\nu$. This and the fact that the values of
$\A(\sloop)$ given in the last column of Table~\ref{tab:Ktildeni} are smaller
than the corresponding values given in lines 3 and 4 proves that a minimizer
$\minloop\in\tilde{\K}^\nu$ is free of total collisions.

In the last 3 lines of Table~\ref{tab:KPi}, for each $(P,i)$, we list
the values of $\A(\sloop)$ corresponding to the minimal sequences in
Table~\ref{tab:KPi_seq}. Again we remark that $\sloop$ enjoys all the
symmetries and the topological constraints required for membership in
$\KPi$.  From this and the fact that the values $\A(\sloop)$ in the
line corresponding to the pair $(P,i)$ are strictly less than the
values in lines 2 or 3 shows that minimizers $\minloop\in\KPi$ are
free of total collisions. This concludes our analysis of total
collisions.  It remains to exclude the occurrence of partial
collisions: this is done in the following Subsection.

\subsection{Partial collisions}
\label{sec:part_coll}

Our strategy to exclude partial collision consists of two steps.  Let us
assume that a minimizer $\minloop\in\overline{\K}$ has a partial collision:
first we show that the collision is isolated; then we prove that we can
construct a local perturbation $\sloop$ with $\A(\sloop) < \A(\minloop)$. In
this construction we can not rely on techniques of the type used in
\cite{marchal01}, \cite{FT04}, based on Marchal's idea of averaging the action
on a set of perturbations $\sloop_\theta$, depending on a parameter
$\theta$. In fact the condition $\sloop_\theta\in \overline{\K}$ is a kind of
unilateral constraint and may be violated for some value of $\theta$.  Besides
in certain cases, for instance for the cones $\KPi$, the use of this technique
is not allowed due to the presence of reflection symmetries.  We base our
discussion of partial collisions on the fact that, as we discuss below, all of
them can be regarded as binary collisions and we can take advantage of the
knowledge of the geometric--kinematic structure of such collisions.

\begin{lemma}
Let $\minloop\in\overline{\K}$ be a minimizer of the action. Assume
that $\minloop$ has a partial collision at time $t_c$, then the
collision is isolated.
\label{lem:isolated}
\end{lemma}
\begin{proof}
1. In \cite{marchal01}, \cite{chenciner02}, \cite{FT04} it is shown that,
if $M\leq N$ particles of the system all collide together at time
$t_c$ and there is $\delta>0$ such that for $t\in(t_c-\delta,
t_c+\delta)$ there is no collision involving only a proper subset of
the $M$ particles colliding at time $t_c$, then $t_c$ is not an
accumulation point of collisions of the $M$ particles.

\noindent
2. Let $G$ be the group of order $|G|=4$ generated by the rotations $R_j$ of
$\pi$ around the axes $\assej, j=1,2,3$, if $\K=\K_4$, or $G = {\cal
  R}\in\{{\cal T},{\cal O}, {\cal I}\}$ otherwise.  As we have already
observed in Proposition~\ref{border&coerc} at the collision time $t_c$ the
generating particle $\minloopgen(t_c)$ lies on one of the axes,
say $r$, of some rotation $R\in G\setminus\{I\}$.  Since we deal with a
partial collision we have that $\minloopgen(t_c)\in r\setminus\{0\}$ and the
collision involves all the $|C|$ particles associated to the maximal subgroup
$C\subset G$ of rotations of axis $r$; there are $|G|/|C|$ contemporary
partial collisions of clusters of $|C|$ bodies.  Assume that there is another
partial collision at time $t'_c$ such that $\minloopgen(t'_c)\in
r'\setminus\{0\}$ for $r'\neq r$ the axis of another rotation in $G$. Then
\[
|t_c' - t_c| \geq
 \frac{N}{2\A(\minloop)} d^2(\minloopgen(t_c),r')
 \stackrel{def}{=}\delta\,;
\]
indeed we have that 
\begin{eqnarray*}
d(\minloopgen(t_c),r') &\leq& |\minloopgen(t'_c) -
\minloopgen(t_c)| \leq
|t'_c-t_c|^{1/2}
\left[\int_{t_c}^{t'_c}
|\dminloopgen(t)|^2\;dt\right]^{1/2} \leq \\
&\leq& |t'_c-t_c|^{1/2} \left(\frac{2\A(\minloop)}{N}\right)^{1/2}\ .
\end{eqnarray*}
It follows that all collisions of $\minloopgen$ on the interval
$(t_c-\delta, t_c+\delta)$ take place on the axis $r$ and involve
exactly the $|C|$ particles associated to $C$.
From this and 1. the Lemma follows.

\rightline{\small$\square$}
\end{proof}

If $\minloop\in\overline{\K}$ is a minimizer and $(t_1,t_2)$ is an
interval of regularity, then $\minloopgen$ is a solution of Newton's
equation
\begin{equation}
\ddot{w} = \sum_{R\in G\setminus\{I\}}\frac{(R - I)w}{
|(R - I)w|^3}\,, \qquad t\in (t_1,t_2) .
\label{K_eq_motion}
\end{equation}
If $r$ is the axis of some rotation in $G$ ($G$ as in
Lemma~\ref{lem:isolated}) and $C\subset G\setminus\{I\}$ is the maximal
subgroup of the rotations with axis $r$ we can rewrite (\ref{K_eq_motion}) in
the form
\begin{equation}
\ddot{w} = \sum_{R\in C\setminus\{I\}}\frac{(R - I)w}{ |(R - I)w|^3} +
\sum_{R\in G\setminus C}\frac{(R - I)w}{ |(R - I)w|^3}\ .
\label{K_motion_decomp}
\end{equation}
If we call $R_\pi$ the rotation of $\pi$ around $r$ and set 
\[
\alpha = \sum_{j =1}^{|C| -1}\frac{1}{\sin\bigl(\frac{j \pi}{|C|}\bigr)}
\]
then we have
\[
\sum_{R\in C\setminus\{I\}}\frac{(R - I)w}{ |(R - I)w|^3} =
\alpha\frac{(R_\pi - I)w}{ |(R_\pi - I)w|^3}\,,
\]
that shows (\ref{K_motion_decomp}) is of the general form
\begin{equation}
\ddot{w} = \alpha \frac{(R_\pi - I)w}{ |(R_\pi - I)w|^3} + V_1(w) \,,
\label{K_motion_general}
\end{equation}
where $V_1(w)$ is a smooth function defined in an open set $\Omega\subset\R^3$
that contains $r\setminus\{0\}$. The form (\ref{K_motion_general}) of Newton's
equation is well suited for the analysis of partial collisions occurring on
$r$ and implies that all partial collisions a minimizer
$\minloop\in\overline{\K}$ may present can be regarded as binary collisions.

By a similar computation the first integral of energy can be written
in the form
\begin{equation}
|\dot{w}|^2 - \alpha \frac{1}{|(R_\pi - I)w|} - V(w)= h \ .
\label{K_energy_general}
\end{equation}
For the case at hand 
\[
V_1(w) =\sum_{R\in G\setminus C}\frac{(R - I)w}{ |(R - I)w|^3},
\qquad
V(w) =\sum_{R\in G\setminus C}\frac{1}{ |(R - I)w|}\ .
\]
From these expressions it follows that if $r'\neq r$ is the axis of
some rotation in $G\setminus\{I\}$ and $\tilde{R}$ is the reflection
with respect to the plane determined by $r, r'$, then $V_1, V$ satisfy
the symmetry conditions
\begin{equation}
V_1(\tilde{R}w) = \tilde{R}V_1(w)\,,
\hskip 1cm
V(\tilde{R}w) = V(w)\ .
\label{symcond_Rtilde}
\end{equation}
In the following Proposition we list a number of properties of {\em
ejection solutions} to (\ref{K_motion_general}), that is solutions such that
\[
\lim_{t\to t_c^+} w(t) = w(t_c)\in r\setminus\{0\}\ .
\]
By shifting the origin of the coordinates and of the time we can assume
$w(t_c)=0$, $t_c=0$. We denote by $\er$ a unit vector parallel to $r$.
\begin{proposition}
Let $w:(0,\bar t)\rightarrow \R^3$ be a maximal solution of
(\ref{K_motion_general}). Assume that
\begin{equation}
\lim_{t\rightarrow 0^+}w(t) = 0\ .
\label{eq:P2.14}
\end{equation}
Then \\
\noindent
(i) there exists $\olda \in \R$ and a unit vector ${\mathsf n}$,
orthogonal to $r$, such that
\begin{equation}
\lim_{t\rightarrow 0^+}\frac{\dot w(t) + R_\pi\dot w(t)}{2} = \olda \er, 
\label{lim_pdot}
\end{equation}
\begin{equation}
\lim_{t\rightarrow 0^+}\frac{w(t) - R_\pi w(t)}{|w(t) - R_\pi w(t)|} =
\lim_{t\rightarrow 0^+}\frac{w(t)}{|w(t) |} = {\mathsf n}\ .
\label{eq:P2.16}
\end{equation}

\noindent(ii) The rescaled function $w^\lambda:[0,1]\to\R^3$ defined by
$w^\lambda(0)=0, w^\lambda(\tau) = \lambda^{2/3} w(\tau/\lambda),
\lambda>1/\bar t$, satisfies
\begin{equation}
\begin{array}{l}
\lim_{\lambda\to+\infty}|w^\lambda(\tau) - s^\alpha(\tau){\mathsf n}|
= 0 \mbox{ uniformly in } [0,1]\,,\cr
\lim_{\lambda\to+\infty}|\dot{w}^\lambda(\tau) -
\dot{s}^\alpha(\tau){\mathsf n}| = 0 \mbox{ uniformly in } [\delta,1],
0<\delta<1\,,\cr
\label{wlambda_eqs}
\end{array}
\end{equation}
where
\[
s^\alpha(\tau) = \frac{3^{2/3}}{2}\alpha^{1/3} \tau^{2/3}\,,\quad
\tau\in[0,+\infty)
\]
is the parabolic ejection motion, that is the solution of $\dot{s} =
(\alpha/2)^{1/2} s^{-1/2}$ that satisfies $\lim_{\tau\to 0^+}s(\tau) = 0$.
\label{prop:limcoll}
\end{proposition}
\begin{proof}
The change of variables 
\[
p = \frac{w + R_\pi w}{2},\qquad q = \frac{w - R_\pi w}{2}
\]
transforms (\ref{K_motion_general}), (\ref{K_energy_general}) into
\begin{equation}
\left\{
\begin{array}{l}
\ddot p = \frac{1}{2}(I + R_\pi) V_1(p + q)\cr 
\ddot q = -\displaystyle\frac{\alpha}{4} \frac{q}{|q|^3} + \frac{1}{2} (I -
R_\pi) V_1(p + q)\cr
\end{array}
\right.
\ ,
\label{pqddot}
\end{equation}
\begin{eqnarray}
|\dot p|^2 + |\dot q|^2 =  \frac{\alpha}{2|q|}  + V(p+q) + h\ .
\label{eq:P2.24}
\end{eqnarray}
Fix a number $d > 0$ such that $B_d = \{|w| < d\}\subset\subset\Omega$ and let
$(0,t_d)$ be the maximal interval in which the solution $(p(t),q(t))$ of
(\ref{pqddot}) remains in $B_d$. For $t\in(0,t_d)$ the boundedness of $V_1$
and (\ref{pqddot})$_1$ implies $|\ddot p| \le C_1$, for some constant $C_1>0$.
This, the assumption (\ref{eq:P2.14}) and the fact that, by definition, $p$ is
parallel to $r$ yield the existence of $\olda \in \R$ such that
\begin{equation}
\lim_{t\rightarrow 0^+}\dot p(t) = \olda \er\,,
\label{eq:P2.25}
\end{equation}
which proves (\ref{lim_pdot}) and implies
\begin{equation}
|\dot p(t)| \le C_2 \,, \qquad t \in (0, \min\{1, t_d\})
\label{eq:P2.26}
\end{equation}
where $C_2$ is a positive constant that depends only on $C_1$ and $\olda$.  If
we set $\rho = |q|$ and $\oldnu = \frac{q}{|q|}$, (\ref{pqddot})$_2$ and
(\ref{eq:P2.24}) become
\begin{equation}
\ddot \rho \oldnu + \frac{1}{\rho}\frac{d (\rho^2\doldnu)}{dt} =
-\frac{\alpha}{4 \rho^2} \oldnu + \frac{1}{2}(I - R_\pi)V_1,
\label{eq:P2.27}
\end{equation}
\begin{equation}
\dot \rho^2 + \rho^2|\doldnu|^2 = \frac{\alpha}{2\rho} + V + h - |\dot
p|^2.
\label{eq:P2.28}
\end{equation}
In the remaining part of the proof and in the following
Propositions~\ref{prop:collection_est} and \ref{prop:unique} $\Cquad$
will denote a positive constant that may depends only on $\olda$ and $h$.
The value of $\Cquad$ can change from line to line.

By projecting (\ref{eq:P2.27}) on $\oldnu$ and on its orthogonal
complement we get
\begin{equation}
\left\{
\begin{array}{l}
\ddot \rho = \rho |\doldnu|^2 - \frac{\alpha}{4 \rho^2} +\frac{1}{2}
((I - R_\pi) V_1)\cdot \oldnu \cr
\frac{d(\rho^2\doldnu)}{dt} + \rho^2 |\doldnu|^2 \oldnu=
\rho\frac{1}{2} ((I - R_\pi) V_1)^{\perp}\cr
\end{array}
\right.
\label{nuproj_and_compl}
\end{equation}
where the suffix $\perp$ denotes the projection on the plane orthogonal to
$\oldnu$. We claim that there is a right neighborhood of $t = 0$ where
\begin{equation}
\label{eq:P2.30}
\dot \rho > 0\ .
\end{equation}
If this is not the case, in any neighborhood of $t = 0$ there is $t_0$ such
that $\dot \rho(t_0) =0$. For $t= t_0$ (\ref{eq:P2.26}) and (\ref{eq:P2.28})
imply
\begin{equation}
\label{eq:P2.31}
\rho |\doldnu|^2 = \frac{\alpha}{2\rho^2} + \frac{1}{\rho} (V + h - 
|\dot p|^2) \ge \frac{\alpha}{2 \rho^2} - \frac{\Cquad}{\rho}, \quad (t =
t_0)\ .
\end{equation}
This inequality and
(\ref{nuproj_and_compl})$_1$ imply
\begin{equation}
\label{eq:P2.32}
\ddot \rho \ge \frac{\alpha}{4\rho^2} - \Cquad(1 + \frac{1}{\rho}) , \quad (t
= t_0)\ .
\end{equation}
Since the assumption (\ref{eq:P2.14}) implies $\phantom{aa}$
$\lim_{t\rightarrow0^+}\rho(t) = 0$, from (\ref{eq:P2.32}) we obtain that all
points $t_0$ in a small neighborhood of $ t =0$ where $\dot \rho(t_0)= 0 $ are
relative minima of $\rho$. This is clearly impossible and (\ref{eq:P2.30}) is
established.  Next we show
\begin{equation}
\lim_{t\rightarrow 0^+}\rho|\doldnu| = 0
\label{eq:P2.33}
\end{equation}
and therefore that by (\ref{eq:P2.25}), (\ref{eq:P2.28}) 
\begin{equation}
\label{eq:P2.34}
\lim_{t\rightarrow 0^+} \dot \rho^2 - \frac{\alpha}{2\rho} = V(0) + h
- \olda^2 = \Cquad\ .
\end{equation}
To show (\ref{eq:P2.33}) we first prove the weaker statement
\begin{equation}
\label{eq:P2.35}
\lim_{t\rightarrow 0^+} \rho^2 |\doldnu| = 0\ .
\end{equation}
Suppose that, on the contrary, there is $\delta > 0$ and a sequence
$\{t_j\}, \, t_j \rightarrow 0^+$ such that $\rho^2|\doldnu| \ge
\delta$ along this sequence. Then for $t = t_j$ (\ref{eq:P2.28}) implies
$(\frac{\delta}{\rho})^2 \le \frac{\alpha}{2 \rho} + V + h - | \dot p|^2$,
which is impossible for large $j$, and (\ref{eq:P2.35}) is
established. If we take the vector product of
$(\ref{nuproj_and_compl})_2$ by $\oldnu$, integrate on $(0, t)$ and use
(\ref{eq:P2.35}), we get
\begin{equation}
\rho^2 \doldnu \times \oldnu = \frac{1}{2}\int_0^t
\rho((I-R_\pi)V_1)^\perp \times \oldnu dt'.
\label{eq:P2.36.0}
\end{equation}
From this it follows that, provided $t$ is restricted to a small neighborhood
of $t = 0$ so that (\ref{eq:P2.30}) holds, we have
\begin{equation}
\label{eq:P2.36}
\rho |\doldnu|\le \frac{1}{2}\int_0^t |((I-R_\pi)V_1)^\perp |dt' \le
\Cquad t\,,
\end{equation}
and (\ref{eq:P2.33}) is established. From (\ref{nuproj_and_compl})$_1$ we
have, for $t$ in a neighborhood of $t = 0$,
\begin{equation}
\label{eq:P2.37}
 \ddot \rho + \frac{\alpha}{4\rho^2}\ge - \Cquad \Rightarrow \frac{d}{dt}(\dot
 \rho^2 -\frac{\alpha}{2\rho}) \ge -2 \Cquad\dot \rho\,,
\end{equation}
where we have also used (\ref{eq:P2.30}).  This and (\ref{eq:P2.34}) implies
that for $t$ in a neighborhood of $t = 0$
\begin{equation}
\label{eq:P2.38}
 \dot \rho^2  \ge  \frac{\alpha}{2\rho} -  \Cquad (1+\rho)\ .
\end{equation}
On the other hand (\ref{eq:P2.28}) implies
\begin{equation}
\label{eq:P2.39}
 \dot \rho^2  \le  \frac{\alpha}{2\rho} +  \Cquad\ .
\end{equation}
The inequalities (\ref{eq:P2.38}), (\ref{eq:P2.39}) imply that there exists
$t_0$ that depends only on $h, \olda$ such that (\ref{eq:P2.30}) holds for
$t\in (0, t_0)$ and moreover
\begin{equation}
\left\{
\begin{array}{ll}
\Cquad t^{\frac{2}{3}} \le \rho(t)\,,&\qquad \rho(t)\le \Cquad
t^{\frac{2}{3}}\cr
\Cquad t^{-\frac{1}{3}} \le \dot\rho(t)\,,&\qquad \dot\rho(t) \le \Cquad
t^{-\frac{1}{3}}, \, t\in (0, t_0]\cr
\end{array}
\right.
\ .
\label{eq:P2.40}
\end{equation}
For later reference we also observe that (\ref{eq:P2.38}), (\ref{eq:P2.39})
imply the asymptotic formulas
\begin{equation}
\rho(t) \propto s^\alpha(t) = \frac{3^{2/3}}{2}\alpha^{\frac{1}{3}}
t^{\frac{2}{3}}\,,
\hskip 1cm 
\dot \rho(t) \propto \dot{s}^\alpha(t) = 3^{-1/3}\alpha^{\frac{1}{3}}
t^{-\frac{1}{3}}\ .
\label{eq:P2.41}
\end{equation}
The inequality $(\ref{eq:P2.40})_1$ and (\ref{eq:P2.36.0}) yield
\begin{eqnarray}
\label{eq:P2.42}
\rho^2 |\doldnu|\le \frac{1}{2} \int_0^t \rho |(I - R_\pi)V_1|d t' \le
\frac{3}{10} \Cquad t^{\frac{5}{3}}, \, t \in (0, t_0]
\end{eqnarray}
which together with $(\ref{eq:P2.40})_1$ imply
\begin{equation}
\label{eq:P2.43}
|\doldnu| \le \Cquad t^{\frac{1}{3}}, \, t \in (0, t_0]\ .
\end{equation}
Therefore we deduce from (\ref{eq:P2.40}) that there exists $\rho_0$
depending only on $h, \olda$ such that
\begin{equation}
\label{eq:P2.44}
\Bigl|\frac{d\oldnu}{d \rho}\Bigr| = \frac{|\doldnu|}{|\dot \rho|}\le
 \Cquad \rho, \, \rho \in (0, \rho_0]\ .
\end{equation}
From the estimate (\ref{eq:P2.43}) it follows that there exists a unit vector
${\mathsf n}$ such that
\begin{equation}
\label{eq:P2.45}
{\mathsf n} = \lim_{t\rightarrow 0^+}\oldnu(t) = \lim_{t\rightarrow
0^+}\frac{w(t) - R_\pi w(t)}{|w(t) - R_\pi w(t)|}\ .
\end{equation}
Moreover by definition $\oldnu(t)\cdot \er = 0$ and therefore ${\mathsf
n}\cdot \er =0$. 
We set
\[
\er \times {\mathsf n} = \eort,
\]
\begin{equation}
\left\{
\begin{array}{l}
q = x {\mathsf n} + z \eort,\cr
p = y \er\cr
\end{array}
\right.
\label{eq:P2.46}
\end{equation}
which imply
\begin{equation}
\label{eq:P2.47}
w = p+q = x {\mathsf n} + y \er + z \eort\ .
\end{equation}
From (\ref{eq:P2.40}) we can take $\rho\in (0, \rho_0]$ as the independent
variable. From (\ref{eq:P2.46}) and $q = \rho \oldnu$ we get
\[
\frac{dx}{d\rho} = \oldnu\cdot {\mathsf n} + \rho \frac{d\oldnu}{d\rho}\cdot
{\mathsf n} = 1 + \Bigl(\oldnu - {\mathsf n} + \rho
\frac{d\oldnu}{d\rho}\Bigr)\cdot {\mathsf n}\ .
\]
This and (\ref{eq:P2.44}) imply
\begin{equation}
\Bigl|\frac{dx}{d\rho} - 1\Bigr| \le \Cquad \rho^2, \, \rho \in (0, \rho_0]\ .
\label{x_come_r}
\end{equation}
From this we derive, using again (\ref{eq:P2.44}),
\begin{eqnarray}
\label{eq:P2.48}
\Bigl|\frac{dz}{d\rho}\Bigr| &=& \Bigl|\frac{d}{d\rho}(\rho \oldnu - x
 {\mathsf n})\cdot \eort \Bigr| = 
\Bigl|\Bigl(\oldnu + \rho\frac{d\oldnu}{d\rho} -
 \frac{dx}{d\rho}{\mathsf n}\Bigr)\cdot\eort\Bigr| \nonumber\\
&&\hskip -1.8cm = \Bigl|\Bigl[\Bigl(1 - \frac{dx}{d\rho}\Bigr) {\mathsf n} +
\oldnu - {\mathsf n} + \rho \frac{d\oldnu}{d\rho}\Bigr]\cdot \eort\Bigr| \le
\Cquad \rho^2,\quad \rho\in(0,\rho_0]\,.
\end{eqnarray}
From (\ref{eq:P2.25}) and (\ref{eq:P2.40}) we get 
\begin{equation}
\label{eq:P2.49}
\Bigl|\frac{dy}{d\rho}\Bigr| \le \Cquad \rho^{1/2}, \, \rho \in (0, \rho_0]\ .
\end{equation}
These estimates and the fact that by (\ref{x_come_r}) we can take $x$ as the
independent variable in some interval $(0,x_0)$, with $x_0>0$ depending only
on $h, \olda$, imply
\begin{equation}
\label{eq:P2.50}
\lim_{t\rightarrow 0^+}\frac{w(t)}{|w(t)|} = \lim_{x\rightarrow
0^+}\frac{{\mathsf n} + \frac{y}{x}\er + \frac{z}{x}\eort}{(1 + \frac{y^2 +
z^2}{x^2})^{\frac{1}{2}}} = {\mathsf n}\,,
\end{equation}
which completes the proof of (\ref{eq:P2.16}). 

We denote by $'$ the derivative with respect to $x$ and observe that
(\ref{eq:P2.48}), (\ref{eq:P2.49}) imply 
\begin{equation}
|y'| \le \Cquad \sqrt{x}\,,
\hskip 1cm
|z'| \le  \Cquad x^2\,, \hskip 1cm x\in(0,x_0)\ .
\label{eq:P2.52}
\end{equation}
Relations (\ref{x_come_r}), (\ref{eq:P2.40}) and (\ref{eq:P2.52}) yield
\[
|x(t) - \rho(t)|\leq \Cquad t^2\,,
\hskip 0.5cm |y(t)| \leq \Cquad t\,,
\hskip 0.5cm |z(t)| \leq \Cquad t^2\ .
\] 
Using also
(\ref{eq:P2.41})$_1$, that implies
\[
1 = \lim_{t\to 0^+} \frac{\rho(t)}{s^\alpha(t)} =
\lim_{\lambda\to+\infty}
\frac{\rho(\tau/\lambda)}{s^\alpha(\tau/\lambda)} =
\lim_{\lambda\to+\infty}
\frac{\lambda^{2/3}\rho(\tau/\lambda)}{s^\alpha(\tau)}\,,
\]
the proof of (\ref{wlambda_eqs})$_1$ follows from the inequality
\begin{eqnarray*}
&&|w^\lambda(\tau)-s^\alpha(\tau){\mathsf n}| \leq
|\lambda^{2/3}\rho(\tau/\lambda) - s^\alpha(\tau)| +
\lambda^{2/3}|x(\tau/\lambda) - \rho(\tau/\lambda)| +\\
&& +\lambda^{2/3}(|y(\tau/\lambda)| + |z(\tau/\lambda)| )\ .
\end{eqnarray*}
The proof of (\ref{wlambda_eqs})$_2$ is similar.

\rightline{\small$\square$}
\end{proof}
For later use, in the following Proposition we collect some of the estimates
obtained in the proof of Proposition~\ref{prop:limcoll}.

\begin{proposition}
Let $w:(0,\bar t)\rightarrow \R^3$ be a maximal solution of
(\ref{K_motion_general}). Assume that
\[
\lim_{t\rightarrow 0^+}w(t) = 0
\]
and define $\olda, {\mathsf n}$ as in Proposition~\ref{prop:limcoll}.
Let $x,y,z,\rho$ be defined by
\[
w = x{\mathsf n} + y\er + z \er\times{\mathsf n}\,,\hskip 1cm
\rho = \frac{1}{2}|(R_\pi - I)w|\ .
\]
Then the following estimates hold for some positive constants $t_0$, $\rho_0$,
$x_0$, $\Cquad$, depending only on $\olda$ and $h$:
\begin{equation}
\left\{
\begin{array}{ll}
\Cquad t^{\frac{2}{3}} \le \rho(t)\,;
&\rho(t)\le \Cquad t^{\frac{2}{3}}\,, \cr
\Cquad t^{-\frac{1}{3}} \le \dot \rho(t)\,;
&\rho(t)\le \Cquad t^{-\frac{1}{3}}, \quad
t\in (0, t_0]\,, \cr
|\frac{dx}{d\rho} - 1| \le \Cquad \rho^2\,, &\rho \in (0, \rho_0]\,,\cr
|y'| \le \Cquad x^{1/2}\,; &|z'| \le  \Cquad x^2\,,\cr
\end{array}
\right.
\label{est_collect}
\end{equation}
where by $'$ we mean differentiation with respect to $x$.
\label{prop:collection_est}
\end{proposition}
Propositions~\ref{prop:limcoll}, \ref{prop:collection_est} are stated and
proved for ejection solutions. Analogous statements and proofs with obvious
modifications apply to collision solutions of (\ref{K_motion_general}), that
is solutions satisfying $\lim_{t\to 0^-}w(t)=0$.

Given unit vectors $\npm$ define $\omega = \omega^{\alpha,\npm}: \R \to
\R^3$ by setting
\[
\omega^{\alpha,\npm}(\pm t) = \npm s^\alpha(t), \ t\geq 0\,,
\qquad
s^\alpha(t) = \frac{3^{2/3}}{2}\alpha^{1/3}t^{2/3}\ .
\]
If $-1<\npiu\cdot \nmeno<1$ we let $\Theta_d$, with $0<\Theta_d<\pi$,
be the angle determined by $O$ and $\npm$ and let $\Theta_i$, with
$\pi<\Theta_i<2\pi$ be the complement of $\Theta_d$. Given $t^+,t^->0$
there exist unique Keplerian arcs $\omega_d: [-t^-,t^+]\to\R^3$ and
$\omega_i:[-t^-,t^+]\to\R^3$ that connect $\omega(-t^-)$ to
$\omega(t^+)$ in the time interval $[-t^-,t^+]$ and satisfy
\[
\left\{
\begin{array}{l}
\omega_d((-t^-,t^+)) \subset \Theta_d\cr
\omega_i((-t^-,t^+)) \subset \Theta_i\cr
\end{array}
\right.\ ;
\]
$\omega_d$ and $\omega_i$ are called the {\em direct} and {\em
indirect} Keplerian arc \cite{albouy}. In the boundary case
$\npiu\cdot \nmeno=-1$ both angles have measure $\pi$ and the
distinction between $\omega_d$ and $\omega_i$ does not make sense. In
the other boundary case $\npiu\cdot \nmeno = 1$ the indirect arc does
not exist and we can assume $\Theta_d = 0, \Theta_i = 2\pi$.

\begin{proposition}
The following inequalities hold:
\begin{itemize}
\item[(i)] $\A(\omega_d) < \A(\omega|_{[-t^-,t^+]}),\ \forall \npm$\,,
\item[(ii)] $\A(\omega_i) < \A(\omega|_{[-t^-,t^+]}),\ \forall \npm$ such
that $\npiu\cdot \nmeno<1$\,,
\end{itemize}
where
\[
\A(w) = \int_{-t^-}^{t^+} \Bigl(\frac{|\dot{w}|^2}{2} + \frac{\alpha}{4|w|}
\Bigr)\;dt\ .
\]
\label{prop:marchal}
\end{proposition}
\begin{proof} 
For the proof of this Proposition we refer to \cite{marchal01},
\cite{chenciner_priv}, and to \cite{terrvent}, where the
generalization to the case of potentials of the form $1/r^\gamma,
\gamma\in(0,2)$, is also considered. Below we present a proof for the
symmetric case $t^+=t^-=\tau$ which is the one we use in the
following.

(The case $t^+=t^-=\tau$) We consider a unit mass $\particle$, moving in the
plane under the attraction of a mass $\mu$ ($\mu=\frac{\alpha}{4}$) fixed in
the origin $O$, and we let $h$ be the energy and $J$ the constant of angular
momentum. We let $(r,\phi)$ be the polar coordinates of $\particle\neq O$. We
define $\rho>0$ and $\theta\in[0,\pi]$ by setting $\npm = (\cos\theta,
\pm\sin\theta)$, $\omega(\pm\tau) = \npm\rho$. Then $\tau$ and $\rho$ are
related by
\begin{equation}
\tau = \frac{2^{1/2}}{3}\frac{\rho^{3/2}}{\mu^{1/2}}\ .
\label{parab_motion_time}
\end{equation}
For each $\theta\in(0,\pi)$ there is a unique Keplerian arc
$\omega_\theta$ that connects $\omega(-\tau)$ with $\omega(\tau)$ in
the time interval $[-\tau,\tau]$ and intersects the polar axis at a
point $(\rho_0,0)$ with $\rho_0>0$. The arc $\omega_\theta$ is
therefore the direct arc if $\theta\in(0,\pi/2)$ and the indirect arc
if $\theta\in(\pi/2,\pi)$. We denote by $2A$ the action of the arc
$\omega_\theta$ and by
\begin{equation}
A_0 = 2^{3/2} (\rho\mu)^{1/2}
\label{parab_action}
\end{equation} 
the action of the parabolic ejection arc connecting $O$ with
$\omega(\tau)$.  To prove the Proposition is the same as to show
that the ratio $a:=A/A_0$ is strictly $< 1$ for all
$\theta\in[0,\pi)$. It is easily seen that
\[
\theta\lesseqqgtr \frac{2^{1/2}}{3} \qquad\Longrightarrow\qquad
\rho_0\gtreqqless \rho
\]
and moreover the eccentricity $e$ of $\omega_\theta$ satisfies
\begin{equation}
\left\{
\begin{array}{ll}
e\in [0,1),       &0<\theta\leq 2^{1/2}/3 \cr
e\in (0,+\infty), &2^{1/2}/3 <\theta \leq \pi/2\cr
e\in (0,-1/\cos\theta),\qquad &\pi/2<\theta\leq \pi \cr
\end{array}
\right.
\ .
\label{e_theta_cond}
\end{equation}
In particular it follows that $(\rho_0,0)$ is the apocenter if
$\theta<2^{1/2}/3$ and the pericenter if $\theta>2^{1/2}/3$. Therefore
the polar equation of $\omega_\theta$ reads
\begin{equation}
r = \frac{J^2/\mu}{1\mp e\cos\phi}\,,\quad -\theta\leq \phi\leq \theta
\label{polar_traject}
\end{equation}
where here and in the following $\mp = -$ if $\theta<2^{1/2}/3$ and
$\mp = +$ if $\theta>2^{1/2}/3$. The constants $h$ and $J$ are related
to $\rho_0$ and $e$ by
\begin{equation}
h = \frac{\mu}{2\rho_0}(-1\mp e)\,,
\label{energy_restpoint}
\end{equation}
\begin{equation}
J^2 = \rho_0\mu(1\mp e)\ .
\label{angmom_squared}
\end{equation}
The values of $\rho_0$ and $e$ are determined by the conditions
\begin{equation}
\rho = \frac{J^2/\mu}{1\mp e\cos\theta}\,,
\label{polar_restpoint}
\end{equation}
\begin{equation}
\frac{1}{J}\int_0^\theta r^2 \;d\phi =
\frac{J^3}{\mu^2}\int_0^\theta\frac{d\phi}{(1\mp e\cos\phi)^2} = \tau
\label{area_law}
\end{equation}
which express the fact that $(r,\phi)=(\rho,\theta)$ fulfills
(\ref{polar_traject}) and that the travel time from $(\rho_0,0)$ to
$(\rho,\theta)$ along $\omega_\theta$ coincides with $\tau$. From
(\ref{angmom_squared}), and (\ref{polar_restpoint}) it follows
\[
\frac{\rho}{\rho_0} = \frac{1\mp e}{1\mp e\cos\theta}\ .
\]
Using this, (\ref{parab_motion_time}) and (\ref{angmom_squared}), we
obtain from (\ref{area_law}) the equation
\begin{equation}
\int_0^\theta \frac{d\phi}{(1\mp e\cos\phi)^2}
=\frac{2^{1/2}}{3}\frac{1}{(1\mp e\cos\theta)^{3/2}}
\label{eq_for_e_of_theta}
\end{equation}
that determines $e=e(\theta)$. For the action ratio we have
\begin{equation}
a = \frac{1}{A_0}\Bigl(h\tau + 2\mu\int_0^\tau\frac{dt}{r}\Bigr) = 
\frac{1}{A_0}\Bigl(h\tau + 2J\int_0^\theta\frac{d\phi}{1\mp e\cos\phi}\Bigr)\ .
\label{action_ratio}
\end{equation}
If we set $I_n = \int_0^\theta\frac{d\phi}{(a+b\cos\phi)^n}$ we have
the identity
\[
(a^2-b^2)I_2 = aI_1 -
b\frac{\sin\theta}{a+b\cos\theta}\ \qquad(a^2\neq b^2)\ .
\]
From this with $a=1, b=\mp e$ and (\ref{eq_for_e_of_theta}) it follows
\[
\int_0^\theta\frac{d\phi}{1\mp e\cos\phi} =
\frac{2^{1/2}}{3}\frac{1-e^2}{(1\mp e\cos\theta)^{3/2}} \mp
\frac{e\sin\theta}{1\mp e\cos\theta}\ .
\]
If we introduce this expression
of $I_1$ into (\ref{action_ratio}) and use
\[
\frac{h\tau}{A_0} = -\frac{1-e^2}{12(1\mp e\cos\theta)}\,,
\qquad
\frac{2J}{A_0} = \frac{1}{2^{1/2}}(1\mp e\cos\theta)^{1/2}\,,
\]
that follow from (\ref{parab_motion_time}), (\ref{parab_action}),
(\ref{energy_restpoint})--(\ref{polar_restpoint}),
we finally get
\begin{equation}
a = \frac{1}{4}\frac{(1-e^2)}{(1\mp e\cos\theta)} \mp 
\frac{1}{2^{1/2}}\frac{e\sin\theta}{(1\mp e\cos\theta)^{1/2}}\ .
\label{a_di_e_theta}
\end{equation}
Note that this expression, derived under the assumption $0<\theta$,
is valid also for $\theta=0$.  To conclude the proof, instead of
studying directly the function $a(\theta)$ obtained by inserting the
solution $e=e(\theta)$ of (\ref{eq_for_e_of_theta}) into
(\ref{a_di_e_theta}), we show that $a$ has a unique maximum $a_M<1$ on
each line
\begin{equation}
e\cos\theta = {\rm const}\,,
\label{ecostheta}
\end{equation}
for $e, \theta$ satisfying (\ref{e_theta_cond}).  Differentiating
(\ref{ecostheta}) with respect to $e$ we get $\theta' =
\frac{1}{e\tan\theta}$.  From this and (\ref{a_di_e_theta}), it
follows that the derivative $a'$ of $a$ with respect to $e$ along the
lines $e\cos\theta = {\rm const}$ is given by
\begin{equation}
a' = \frac{1}{2^{1/2}\sin\theta(1\mp e\cos\theta)^{1/2}}
\Bigl(\mp 1 - \frac{e\sin\theta}{2^{1/2}(1\mp e\cos\theta)^{1/2}}
\Bigr)\ .
\label{aprimo}
\end{equation}
For $\theta\leq 2^{1/2}/3$ we have from (\ref{e_theta_cond}) that $e\in[0,1)$,
and (\ref{aprimo}) implies $a'\leq 0$, therefore $a$ takes its maximum $a_M =
\frac{1}{4}(1+e)<\frac{1}{2}$ for $\theta=0$. For $\theta>2^{1/2}/3$ we have
$a'=0$ on the curve $\ell$ defined by
\begin{equation}
1 = \frac{e\sin\theta}{2^{1/2}(1+e\cos\theta)^{1/2}}
\quad
\Longleftrightarrow
\quad
e = e_0(\eta) \stackrel{def}{=} (1+(1+\eta)^2)^{1/2} >1\,,
\label{ell_curve}
\end{equation}
with $\eta=e\cos\theta$. From (\ref{aprimo}) it follows that 
\[
a' \lesseqqgtr 0 \ \ \Longleftrightarrow \ \ e\gtreqqless e_0
\]
and therefore $a$ attains its maximum $a_M$ on $\ell$.  Inserting
(\ref{ell_curve}) into (\ref{a_di_e_theta}) yields
\[
a_M = 1 - \frac{1}{4}(1+e\cos\theta) < 1\ .
\]
This concludes the proof.

\rightline{\small$\square$}
\end{proof}

It is exactly the possibility of choosing between $\omega_d$ and
$\omega_i$ still reducing the action that, whenever $\minloop$ is
assumed to have a partial collision that implies
$\minloop\in\partial\K$, allows us to perturb $\minloop$ inside $\K$,
thus preserving the constraint of membership in $\K$. On the basis of
Proposition~\ref{prop:marchal} this can always be done if $\npiu\cdot
\nmeno<1$. The special case $\npiu\cdot \nmeno = 1$ is excluded since in this
case the indirect Keplerian arc does not exist at all and, if the
direct Keplerian arc does not allow the construction of a perturbation
$\sloop\in\K$, then Proposition~\ref{prop:marchal} can not be used.
The discussion of these situations is more delicate and it is based on
a uniqueness result that we prove below
(cfr. Proposition~\ref{prop:unique}).

\medbreak
We begin our analysis of partial collisions by
\begin{proposition}
Let $\minloop\in\overline{\K}$ be a minimizer of the action and
assume that $\minloop$ has a partial collision at time $t_c$. Let
$\npiu, \nmeno$ be the unit vectors associated to the collision of the
generating particle in the sense of Proposition~\ref{prop:limcoll}. Then 
\[
\npiu = \nmeno\ .
\]
\label{prop:npiunmeno}
\end{proposition}
\begin{proof}
We show that the assumption that $\minloop$ has a partial collision with
$\npiu\cdot\nmeno<1$ leads to the contradiction of the existence of a
perturbation $v\in\Kclo$ of $\minloop$ such that $\A(v)<\A(\minloop)$.  The
equivariance condition that characterizes $u\in\K$ implies that it suffices to
define $v$ only in a fundamental interval $I_\K$ that contains the collision
time $t_c$. Then $v$ is automatically extended to the whole $\R$ by
equivariance.

We call $r$ the axis where the collision of the generating particle takes
place and $C\subset G$ the maximal subgroup of the rotations with axis $r$.
Set $w(t) = \minloopgen(t_c+t) -\minloopgen(t_c)$.  For every fixed
$\lambda>1$ the restriction of $w$ to $\left[-\frac{1}{\lambda},
\frac{1}{\lambda}\right]$ is a minimizer of
\begin{eqnarray*}
\A^\lambda(\phi) &=& \lambda^{1/3} \frac{|G|}{2}
\int_{-1/\lambda}^{1/\lambda}\biggl(| \dot\phi|^2 + \sum_{R\in
C\setminus\{I\}} \frac{1}{|(R-I)\phi|} \biggr)\,dt +\\
&& \lambda^{1/3}\frac{|G|}{2} \int_{-1/\lambda}^{1/\lambda}\sum_{R\in
G\setminus C} \frac{1}{|(R-I)(\phi +\minloopgen(t_c))|}\, dt
\end{eqnarray*}
on the set of functions $\phi$ in
$H^1\left((-\frac{1}{\lambda},\frac{1}{\lambda}),\R^3 \right)$ that satisfy
\[
\phi\Bigl(\pm\frac{1}{\lambda}\Bigr) = w\Bigl(\pm\frac{1}{\lambda}\Bigr)\ . 
\]
The map
\[
f:H^1\Bigl(\Bigl(-\frac{1}{\lambda},\frac{1}{\lambda}\Bigr),\R^3 \Bigr) \to 
H^1\left((-1,1),\R^3 \right)
\]
defined by 
\begin{equation}
\left\{
\begin{array}{l}
f(\phi) = \psi\cr
\phi(t) = \lambda^{-2/3}\psi(\lambda t), \ \ t\in[-1/\lambda,1/\lambda]\cr
\end{array}
\right.
\label{scaling}
\end{equation}
is a bijection and we have
\begin{eqnarray}
&&\A^\lambda(\phi) = \hat{\A}^\lambda(\psi) \stackrel{def}{=}
\frac{|G|}{2} \int_{-1}^1 \Bigl(\Bigl|\frac{d\psi}{d\tau}\Bigr|^2 +
\sum_{R\in C\setminus\{I\}} \frac{1}{|(R-I)\psi|} \Bigr)\,d\tau
+\nonumber \\
&&\hskip -1cm +\frac{|G|}{2} \int_{-1}^1 \sum_{R\in G\setminus C}
\frac{1}{|(R-I)(\psi +\lambda^{2/3}\minloopgen(t_c))|}\, d\tau
\stackrel{def}{=} a(\psi) + a^\lambda(\psi)\,,
\label{rescaled_action}
\end{eqnarray}
where $a(\psi)$ and $a^\lambda(\psi)$ denote the two terms in the
definition of $\hat{\A}^\lambda(\psi)$.  Therefore from the minimality
of $w|_{[-1/\lambda,1/\lambda]}$ and (\ref{rescaled_action}) it
follows that the map $w^\lambda:[-1,1]\to\R^3$ defined by
\[
w^\lambda(\tau) = \lambda^{2/3} w\Bigl(\frac{\tau}{\lambda}\Bigr)
\]
is a minimizer of $\hat{A}^\lambda(\psi) = a(\psi) + a^\lambda(\psi)$.  From
Proposition~\ref{prop:limcoll}, for $\lambda\to+\infty$, $w^\lambda$ converges
uniformly in $[-1,1]$ to $w^\infty = \omega^{\alpha,\npm}$ with $\npm$
orthogonal to $r$.  Assume $\npiu\cdot \nmeno < 1$ and let $w^\pm$ be the
direct and indirect Keplerian arcs connecting $w^\infty(-1)$ to $w^\infty(+1)$
in the interval $[-1,1]$, and define
\[
\hat{w}^{\lambda,\pm} = w^\pm + \tilde{w}^\lambda\,,
\]
where we have set
\[
\tilde{w}^\lambda(\tau) = (w^\lambda(-1) - w^\infty(-1))(\frac{1-\tau}{2})
+ (\frac{1+\tau}{2})(w^\lambda(1) - w^\infty(1)),\hskip 0.2cm
\tau\in[-1,1].
\]
From (\ref{rescaled_action}), the boundedness of
$\hat{w}^{\lambda,\pm}$ and Lebesgue's dominate convergence theorem we have
\begin{equation}
\left\{
\begin{array}{l}
\lim_{\lambda\to+\infty}\hat{A}^\lambda(w^\lambda) = a(w^\infty)\cr
\lim_{\lambda\to+\infty}\hat{A}^\lambda(\hat{w}^{\lambda,\pm}) =
a(w^\pm)\cr
\end{array}
\right.\ .
\label{twolims}
\end{equation}
This and Proposition~\ref{prop:marchal} imply that, for $\lambda>>1$,
\begin{equation}
\A({\sloop}^{\lambda,\pm}) < \A(\minloop)\,,
\label{action_ineq}
\end{equation}
where ${\sloop}^{\lambda,\pm}$ is defined through
\[
{\sloopgen}^{\lambda,\pm}(t) =
\left\{
\begin{array}{l}
\minloopgen(t), \ t\in I_\K\setminus [t_c-1/\lambda,t_c+1/\lambda]\cr
\lambda^{-2/3}\hat{w}^{\lambda,\pm}(\lambda(t-t_c)) + \minloopgen(t_c), \
t\in[t_c-1/\lambda,t_c+1/\lambda]\cr
\end{array}
\right. \ .
\]
The inequality (\ref{action_ineq}) contradicts the minimality of $\minloop$
because the definition of $\K$ implies that either ${\sloop}^{\lambda,+}$ or
${\sloop}^{\lambda,-}$ belong to $\overline{\K}$.  In the above argument we
have tacitally assumed that $t_c$ is in the interior of $I_\K$.  If $t_c$ is
one of the boundary points of $I_\K$ the same argument applies verbatim with
the provision of replacing $t\in[t_c- 1/\lambda, t_c+ 1/\lambda]$ with
$t\in[t_c- 1/\lambda, t_c+ 1/\lambda]\cap I_\K$ in the definition of
$v_1^{\lambda,\pm}(t)$.

\rightline{\small$\square$}
\end{proof}

\begin{remark}
The idea of constructing local variations of the parabolic ejection--collision
orbit obtained by blowing up the rescaled collision solution
has already been used in \cite{ventPhD}, \cite{FT04}.
\end{remark}
\begin{remark}
Let $(\sigma,n)$ be a pair such that the sequence $\sigma$ is simple in the
sense of Definition~\ref{def:simple}.  The argument in the proof of
Proposition~\ref{prop:npiunmeno} can be applied to show that a minimizer
$\minloop\in\overline{\K^{(\sigma,n)}}$ can not have isolated partial
collisions such that the unit vectors $\npm$ associated to the collisions of
the generating particle, as in Proposition~\ref{prop:limcoll}, satisfy
$\npiu\cdot \nmeno<1$.  In fact the assumption that $\sigma$ is simple and the
fact that the curve $\gamma=\{w^\pm(\tau) + \minloopgen(t_c), \tau\in[-1,1]\}$
is a simple closed curve linked to the axis $r$ imply that either
${\sloop}^{\lambda,+}$ or ${\sloop}^{\lambda,-}$ belongs to $\overline{\K}$.
If $\sigma$ is not simple it may be very difficult, or even impossible, to use
a local argument, as in Proposition~\ref{prop:npiunmeno}, to show that a
minimizer $\minloop\in\overline{\K^{(\sigma,n)}}$ does not have partial
collisions. Indeed from Gordon's Theorem \cite{gordon77} we know that, in the
class of loops with index $\neq -1,0,1$ with respect to the origin, the
only minimizers of the planar Kepler problem with an attracting fixed mass at
the origin are collision--ejection loops.
\end{remark}
\begin{remark}
The arguments in the proof of Proposition~\ref{prop:npiunmeno} are based on
the possibility of choosing between $\sloop^{\lambda,+}$ and
$\sloop^{\lambda,-}$, and in turn on the possibility of choosing between the
direct and indirect Keplerian arcs used in the construction of
$\sloop^{\lambda,\pm}$. If $\npiu = \nmeno$, then the indirect arc does not
exist and the construction in the proof of Proposition~\ref{prop:npiunmeno}
yields a single perturbation $\sloop^\lambda$ of $\minloop$ by utilizing the
direct Keplerian arc.  In all cases where it results
$\sloop^\lambda\in\overline{\K}$, again one reaches the contradiction
$\A(\sloop^\lambda)<\A(\minloop)$ and the supposed partial collision is
excluded.
\label{rem:npiunmeno}
\end{remark}
From Remark~\ref{rem:npiunmeno} it follows that a minimizer
$\minloop\in\overline{\K}$ is collision free if it can be excluded that
$\minloop$ has an isolated partial collision such that
\begin{itemize}
\item[(i)] $\npiu = \nmeno = \mathsf{n}$;
\item[(ii)] $\sloop^\lambda \not\in\overline{\K}$,
\end{itemize}
where $\sloop^\lambda$ is the perturbation constructed as in the proof
of Proposition~\ref{prop:npiunmeno}, by means of the direct Keplerian
arc.  

\noindent
We refer to partial collisions that satisfy (i) and (ii) as
collisions of type $(\rightrightarrows)$
\begin{lemma}
Assume $\minloop\in\overline{\K}$ has a collision of type
$(\rightrightarrows)$ at $t = t_c$; let $r$ be the axis on which the
collision of the generating particle takes place and set
$\mathsf{n}=\npm$. Then
\[
\mathsf{n} \in\Pi\,,
\]
where $\Pi$ is a plane through $r$ that contains two distinct axes of
rotations in $G\setminus\{I\}$.
\label{lem:n_in_Pi}
\end{lemma}

\begin{proof}
We give different proofs for the cases $\K = \K_4, \KPi$ and
$\tilde{\K}^\nu$.

1) $\K = \K_4$. From (\ref{4body_cond1}) it follows $\genp(t+T/2) = R_1
  \genp(t)$. This implies that if a minimizer $\minloop\in\overline{\K}$ has a
  collision at time $t_c$, then it also has a collision at time $t_c+T/2$. The
  constraint (\ref{conecond4b}) that defines $\K$ involves the values of
  $\genp(t)$ only for $t=0,T/4$.  Therefore if $t_c\not\in \{0,T/4\} \mod
  T/2$, then any sufficiently small perturbation of $\minloop$ in a compact
  interval $[t_c-\delta,t_c+\delta]$ will remain in $\overline{\K}$ provided
  it satisfies conditions~(\ref{4body_cond1}), (\ref{4body_cond2}).  It
  follows that the collision is not of type $(\rightrightarrows)$.  By
  consequence if $\minloop\in\overline{\K}$ has a collision of type
  $(\rightrightarrows)$, then necessarily $t_c\in \{0,T/4\} \mod T/2$.  From
  this and (\ref{4body_cond1}) it follows
\[
\left\{
\begin{array}{l}
\npiu = S_3 \nmeno \qquad \mbox{ if } t_c = 0\,,\cr
\npiu = S_2 \nmeno \qquad \mbox{ if } t_c = T/4\,,\cr
\end{array}
\right.
\]
which together with $\mathsf{n}=\npm$ imply
\[
\left\{
\begin{array}{l}
\mathsf{n} \in {\rm span}\{\euno,\edue\} \qquad \mbox{ if } t_c = 0\,,\cr
\mathsf{n} \in {\rm span}\{\euno,\etre\} \qquad \mbox{ if } t_c = T/4\ .\cr
\end{array}
\right.
\]

2) $\K = \KPi$. By arguing as in 1) we conclude that if a minimizer
  $\minloop\in\overline{\K}$ has a collision of type
  ($\rightrightarrows)$ at time $t_c$, then $t_c\in \{0,T/(2H)\} \mod
  T/H$. Then from $\condB$, $\condC$ it follows 
\begin{equation}
\left\{
\begin{array}{l}
\npiu = S_3 \nmeno \qquad\ \ \mbox{ if } t_c = 0\,,\cr
\npiu = R S_3 \nmeno \qquad \mbox{ if } t_c = \frac{T}{2H}\,,\cr
\end{array}
\right.
\label{npiuS3}
\end{equation}
where $R$ is the rotation of $2\pi/H$ around $\asseuno$. The operator
$R S_3$ coincides with the reflection with respect to the plane
determined by $\asseuno$ and $V$. From this observation, from equation
(\ref{npiuS3}) and $\mathsf{n}=\npm$ it follows
\[
\left\{
\begin{array}{l}
\mathsf{n} \in {\rm span}\{\euno,\eM\} \qquad \mbox{if } t_c = 0\,,\cr
\mathsf{n} \in {\rm span}\{\euno,\eV\} \qquad \mbox{ if } t_c = T/(2H)\ .\cr
\end{array}
\right.
\]

3) $\K = \tilde{\K}^\nu$, $\nu$ as in Theorem~\ref{platorb_exist}.  From the
discussion at the end of Section~\ref{subsec:geostruct}, if $\minloop
\in\overline{\K}$ is a minimizer which has a partial collision, then a
collision of the generating particle takes necessarily place on one of the
semiaxes in the sequence $k \to r_k:=(\overline{S_{k+1}}\cap
\overline{S^k})\setminus\{0\}$.  For each $k$ we let $\tilde{k}>k$ be defined
by the condition that $S^{\tilde{k}}$ borders with $S^k$ along $r_k$ in the
sense that $r_k \subset \overline{S^{\tilde{k}}}\cap \overline{S^k}$.  The
formal definition of $\tilde{k}$ is: $\tilde{k} = k+h$, where $h\geq 1$ is
determined by the conditions (cfr. Figure~\ref{fig:Dk_seq})
\[
\left\{
\begin{array}{l}
r_k \subset \overline{S^k} \cap \overline{S^{k+h}}\cr
0<j<h \Rightarrow \overline{S^k} \cap \overline{S^{k+j}} = \emptyset\cr
\end{array}
\right.\ .
\]
\begin{figure}[ht]
\psfragscanon
\psfrag{Dk}{$D_k$}
\psfrag{Dkp1}{$D_{k+1}$}
\psfrag{Dkp2}{$D_{k+2}$}
\psfrag{Dkp3}{$D_{k+3}$}
\psfrag{Dkp4}{$D_{k+4}$}
\psfrag{Dkp5}{$D_{k+5}$}
\psfrag{Dkp6}{$D_{k+6}$}
\psfrag{Sk}{$S_k$}
\psfrag{Skp7}{$S_{k+7}$}
\psfrag{S^k}{$S^k$}
\psfrag{SkSkp6}{$S^{\tilde{k}} = S^{k+6}$}
\psfrag{th}{$\theta$}
\psfrag{rk}{$r_k$}
\psfrag{rkp1}{$r_{k+1}$}
\psfrag{rkp2}{$r_{k+2}$}
\psfrag{rkp3}{$r_{k+3}$}
\psfrag{rkp4}{$r_{k+4}$}
\psfrag{rkp5}{$r_{k+5}$}
\psfrag{rkp6}{$r_{k+6}$}
\psfrag{n+}{$\npiu$}
\psfrag{n-}{$\nmeno$}
\centerline{\epsfig{figure=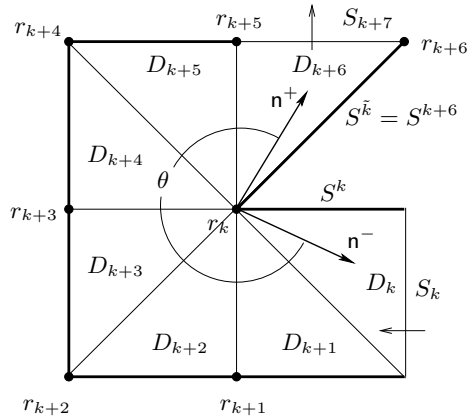,width=6cm}}
\psfragscanoff
\caption{The definition of $S^{\tilde{k}}$. The marked segments correspond to
  $S^j$, $j=k,\ldots,k+6=\tilde{k}$.}
\label{fig:Dk_seq}
\end{figure}

\noindent Let $E = \cup_{j=0}^h (S_{k+j}\cup D_{k+j}\cup S_{k+j+1})$.
Since $\minloopgen$ is the limit of a sequence of maps that enter $E$
through $S_k$ and exit $E$ through $S_{k+h+1}$, the unit vectors
$\nmeno$ and $\npiu$ applied at $\minloopgen(t_c)\in r_k$ point inside
$\overline{E}$.  We let $0\leq \theta\leq 2\pi$ be the angle swept by
$\nmeno$ when
\begin{itemize}
\item[(i)] $\nmeno$ is rotated around $r_k$ until it coincides with $\npiu$,

\item[(ii)] during the rotation $\nmeno$ always points inside $E$.
\end{itemize}
We observe that $\theta = 2\pi$ is equivalent to 
\[
S^k = S^{\tilde{k}}\,,\hskip 1cm \nmeno = \npiu = \mathsf{n}\in S^k\ .  
\]
This concludes the proof. Indeed, if $0\leq \theta < 2\pi$, by using either
the direct or the indirect Keplerian arc, as in the proof of
Proposition~\ref{prop:npiunmeno}, the collision can be excluded. More
precisely, to construct the perturbation of $\minloop$ we use the direct
Keplerian arc if $0\leq \theta\leq \pi$ and the indirect arc if instead
$\pi\leq \theta<2\pi$.

\rightline{\small$\square$}
\end{proof}

For the analysis of collisions of type $(\rightrightarrows)$ we need
a uniqueness result that we present next.  Assume
$\minloop\in\overline{\K}$ has an isolated partial collision at time
$t_c$, then, as we have discussed at the beginning of this Section, we
have $\minloopgen(t_c) \in r\setminus\{0\}$, $r$ being the axis of one
of the rotations in $G\setminus\{I\}$, and $w(t) = \minloopgen(t_c+t) -
\minloopgen(t_c)$ is a solution of (\ref{K_motion_general}).

\begin{proposition}
Let $w_i:(0,\bar t_i)\rightarrow \R^3, \bar{t}_i>0, i=1,2$ be two maximal
solutions of (\ref{K_motion_general}) such that
\[
\lim_{t\rightarrow 0^+}w_i(t) = 0\ .
\]
If $h_i$, $\olda_i$, ${\mathsf n}_i$  are the corresponding values
of the energy and the values of $\olda$ and ${\mathsf n}$ given
by Proposition~\ref{prop:limcoll}, then
\[
\left\{
\begin{array}{l}
h_1 = h_2\cr
\olda_1 = \olda_2\cr
{\mathsf n}_1 = {\mathsf n}_2\cr
\end{array}
\right.
\hskip 1cm\Longrightarrow
\hskip 1cm
\left\{
\begin{array}{l}
\bar{t}_1 = \bar{t}_2\cr
w_1 = w_2\cr
\end{array}
\right.
\ .
\]
\label{prop:unique}
\end{proposition}
\begin{proof}
We project (\ref{K_motion_general}) onto ${\mathsf n}, \er,
\eort$. 
Recalling that $p=y\er, q=x\mathsf{n} + z\eort$ and (\ref{pqddot})
we get
\begin{equation}
\left\{
\begin{array}{l}
\ddot x = -\frac{\alpha}{4 x^2( 1+ \frac{z^2}{x^2})^{\frac{3}{2}}} +
V_1 \cdot {\mathsf n}\cr
\ddot y = V_1 \cdot \er\cr
\ddot z = -\frac{\alpha z}{4 x^3( 1+ \frac{z^2}{x^2})^{\frac{3}{2}}} +
V_1 \cdot\eort\cr
\end{array}
\right.
\ .
\label{xyzddot}
\end{equation}
From Proposition~\ref{prop:collection_est} we can take $x$ as the
independent variable.
We rewrite (\ref{K_energy_general}) in the form
\begin{equation}
\label{eq:P2.53}
\dot x^2(1 + |y'|^2 + |z'|^2) = \frac{\alpha}{2 x( 1+
\frac{z^2}{x^2})^{\frac{1}{2}}} + V + h\,,
\end{equation}
where $'$ denotes differentiation with respect to $x$.
From this and $(\ref{xyzddot})_1$ we get
\begin{equation}
\left\{
\begin{array}{l}
\displaystyle\frac{1}{\dot x ^2} =\frac{2 x}{\alpha}\frac{(1 + |y'|^2
+ |z'|^2)( 1+ \frac{z^2}{x^2})^{\frac{1}{2}}}{ 1 + 2x( 1+
\frac{z^2}{x^2})^{\frac{1}{2}} \frac{V + h}{\alpha}} :=
\frac{2x}{\alpha}(1 + \mathcal W)\cr
\cr
\displaystyle\frac{\ddot x}{\dot x ^2} =-\frac{1}{2x}\frac{(1 + |y'|^2
+ |z'|^2)}{( 1+ \frac{z^2}{x^2})} \frac{\left( 1- 4x^2( 1+
\frac{z^2}{x^2})^{\frac{3}{2}} \frac{V_1\cdot
\mathsf{n}}{\alpha}\right)}{\left(1 + 2x( 1+
\frac{z^2}{x^2})^{\frac{1}{2}}\frac{V + h}{\alpha}\right)}:=
-\frac{1}{2x}(1 + \mathcal U)\cr
\end{array}
\right.
\ .
\label{eq:P2.54}
\end{equation}
Therefore, taking into account that for any function $f(x(t))$ we have
\begin{equation}
\dot f = f' \dot x\,, \hskip 1cm
\ddot f =( f''  + \frac{\ddot x}{\dot x^2} f')\dot x^2\,,
\label{eq:P2.55}
\end{equation}
we can rewrite system (\ref{xyzddot}) in the form
\begin{equation}
\label{eq:P2.56}
\left\{
\begin{array}{l}
y'' - \frac{1}{2 x}(1 + \mathcal U) y' = 2x(1 + \mathcal W)
\frac{V_1 \cdot \er}{\alpha}:= x \mathcal A\,,\cr
z'' - \frac{1}{2 x}(1 + \mathcal U) z' = 
\left[ -\frac{z}{2 x^2( 1+
\frac{z^2}{x^2})^{\frac{3}{2}}} + \frac{2x}{\alpha} V_1 \cdot
\eort\right](1 + \mathcal W) \cr
\hskip 2.8cm := - \frac{z}{2 x^2} (1 + \mathcal
V) + x \mathcal B\ . \cr
\end{array}
\right.
\end{equation}
We now change the independent variable and rewrite (\ref{eq:P2.56}) as a first
order system. We set $x = e^s$, $s \in(-\infty, s_0]$, where $s_0 < 0$ is
chosen later.  We introduce the new variables
$$
\eta = \frac{d y}{ds}, \qquad \zeta = \frac{dz}{ds}
$$
and observe that 
\[
y' = \frac{1}{x} \frac{dy}{ds} = \frac{\eta}{x}\,, \qquad
y'' = \frac{1}{x^2}\Bigl( \frac{d\eta}{ds} - \eta\Bigr) 
\]
and similarly for $z$. If we insert these expressions into
(\ref{eq:P2.56}) and multiply by $x^2 = e^{2s}$, we get the first order
system
\begin{equation}
\left\{
\begin{array}{ll}
\displaystyle\frac{dy}{ds} = \eta\,,
\hskip 0.5cm
&\displaystyle\frac{d\eta}{ds}= \Bigl(\frac{3}{2} + \frac{\mathcal U}{2}\Bigr)
\eta + e^{3s}\mathcal A \cr &\cr
\displaystyle\frac{dz}{ds} = \zeta\,,
\hskip 0.5cm
&\displaystyle\frac{d\zeta}{ds}= \Bigl(\frac{3}{2} + \frac{\mathcal
U}{2}\Bigr) \zeta - \Bigl(\frac{1}{2} + \frac{\mathcal V}{2}\Bigr)z +
e^{3s}\mathcal B\cr
\end{array}
\right.\ .
\label{eq:P2.58}
\end{equation}
We rewrite (\ref{eq:P2.58}) in the compact form
\begin{equation}
\frac{d \gamma}{ds} = M \gamma + \mathcal{N}(\gamma)
\label{dgammads}
\end{equation}
where $\gamma = (y, z, \eta, \zeta)^T$, $\mathcal{N}(\gamma) =(0,0,
\frac{\mathcal U}{2}\eta + e^{3s}\mathcal A, \frac{\mathcal U}{2} \zeta -
\frac{\mathcal V}{2}z + e^{3s}\mathcal B)^T$ and $M$ is the constant matrix
\[
M=
\left(
\begin{array}{cccc}
0&0&1&0\\
0&0&0&1\\
0&0& \frac{3}{2}&0\\
0&-\frac{1}{2}&0& \frac{3}{2}
\end{array}
\right)\ .
\]
From the analysis in the proof of Proposition~\ref{prop:limcoll} and the
equivalence between systems (\ref{xyzddot}) and (\ref{dgammads}) it follows
that to each solution $w$ of (\ref{K_motion_general}) that satisfies
(\ref{eq:P2.14}) there corresponds a solution $\gamma_w$ of (\ref{dgammads})
and $w\neq \tilde{w} \Rightarrow \gamma_w \neq \gamma_{\tilde{w}}$.  From the
estimates in Proposition~\ref{prop:collection_est} it is straightforward to
check that $\gamma_w$ satisfies
\begin{equation}
|\gamma_w(s)| \leq \Cquad e^{\frac{3}{2}s}\,,\qquad s\in(-\infty,s_0]\ .
\label{stima_gammaw}
\end{equation}
A key point in the proof of the claimed uniqueness is the fact that the
constant $\Cquad$, as discussed in the proof of
Proposition~\ref{prop:limcoll}, depends only on $h, \olda$ and therefore is the
same for all solutions of (\ref{K_motion_general}), (\ref{eq:P2.14}) that can
be associated to given $h, \olda, \mathsf{n}$.

The matrix $M$ has eigenvalues $\lambda_i$ and eigenvectors $\rho_i$ as
follows
\begin{equation}
\begin{array}{ll}
\lambda_1 = 0\,,    \qquad  &\rho_1 =(1,0,0,0)^T\nonumber\cr
\lambda_2 = \frac{1}{2}\,,  &\rho_2 =(0,1,0,\frac{1}{2})^T\nonumber\cr
\lambda_3 = 1\,,            &\rho_3 =(0,1,0,1)^T\nonumber\cr
\lambda_4 = \frac{3}{2}\,,  &\rho_4 =(1,0,\frac{3}{2},0)^T \cr
\end{array}\hskip 1cm ;
\label{eq:P2.61}
\end{equation}
these properties of the matrix $M$ imply the estimate
\begin{equation}
|e^{Ms}| \le C e^{3/2s}    
\label{stima_eMs}
\end{equation}
for some constant $C>0$ and $s\in[0,+\infty)$.

Let $P$ be the matrix $[\rho_1,\rho_2, \rho_3,\rho_4]$. Given a
constant $\delta \in \R$ define $\hat \delta\in \R^4$ by setting
\[
P\hat \delta = \delta \rho_4
\] 
and let
\begin{equation}
\label{eq:P2.63}
\gamma_\delta(s) = e^{Ms}P \hat \delta = e^{\frac{3}{2}s} \delta \rho_4\ .
\end{equation}
Clearly $\gamma_\delta$ is a solution of the homogeneous equation
$\frac{d\gamma}{ds}Ê = M \gamma$. Given $K > 0$ and $c\in
(0, \frac{1}{2}]$, consider the set $X$ of continuous maps
$\gamma : (-\infty, s_0]\rightarrow \R^4$ defined by
\begin{equation}
X =\{ \gamma :(-\infty, s_0]\rightarrow \R^4 : |(\gamma -\gamma_\delta)(s)| \le
K e^{(1 + c)s}\}\ .
\label{eq:P2.64}
\end{equation}
The set $X$ with the distance
\begin{equation}
\label{eq:P2.64.1}
d(\gamma, \tilde \gamma)=\max_{(-\infty, s_0]}|\gamma(s) -
\tilde\gamma(s)| e^{-s}
\end{equation}
is a complete metric space.  Note that, for each fixed $\delta$ and provided
the constant $K$ is taken sufficiently large, the estimate
(\ref{stima_gammaw}) implies $\gamma_w\in X$ for all solutions $w$ of
(\ref{K_motion_general}), (\ref{eq:P2.14}) corresponding to given values of
$h, \olda, \mathsf{n}$.  Solutions to (\ref{K_motion_general}) corresponds
through (\ref{pqddot}), (\ref{xyzddot}), etc. to continuous solutions $\gamma
:(-\infty, s_0]\rightarrow \R^4$ of the nonlinear integral equation
\begin{equation}
\gamma(s)=\gamma_\delta(s) + \int_{-\infty}^s e^{M(s - r)}{\mathcal
N}(\gamma(r))dr.
\label{integral_eq}
\end{equation}
To conclude the proof of Proposition~\ref{prop:unique} we shall show that

(I) If $-s_0> 0$ is sufficiently large, then
\begin{equation}
\label{eq:P2.66}
(T\gamma)(s)=\gamma_\delta(s) + \int_{-\infty}^s e^{M(s - r)}{\mathcal
N}(\gamma(r))dr
\end{equation}
defines a contraction on $X$ and therefore (\ref{integral_eq}) has a
unique solution for each $\delta \in \R$.

(II) the choice of $\delta$ is uniquely determined by the value of $\olda$ in
(\ref{lim_pdot}).

\noindent
To prove (I) we need to use the estimates in
Proposition~\ref{prop:collection_est} to derive corresponding estimates for
the functions ${\cal A},{\cal B},{\cal U},{\cal V}$ appearing in the
expression of ${\cal N}$ in (\ref{dgammads}). 
From (\ref{eq:P2.54}) and Proposition~\ref{prop:collection_est} we have
\begin{equation}
\left\{
\begin{array}{ll}
\mathcal U, \mathcal V &= O(|y'|^2 + |z'|^2 + \frac{|z|^2}{x^2} + x ) 
= \cr
&= O(e^{-2s}(|\eta|^2 + |\zeta|^2 + |z|^2 ) + e^s) = O(e^s)\,,\cr
\mathcal A, \mathcal B &= O(1)\ .\cr
\end{array}
\right.
\label{UVAB_est}
\end{equation}
Similarly we have the following estimates for the gradients:
\begin{equation}
\mathcal {U}_\gamma\,, \mathcal {V}_\gamma\,, 
\mathcal {A}_\gamma\,, \mathcal {B}_\gamma = O(e^{-\frac{s}{2}}) \ .
\label{grad_est}
\end{equation}

We have from (\ref{UVAB_est}), (\ref{grad_est}) and (\ref{stima_gammaw})
\begin{equation}
\begin{array}{ll}
|(\mathcal {U }\eta)(s)- (\tilde{\mathcal {U }} \tilde\eta)(s)| 
&\le |(\mathcal {U }(s)- \tilde{\mathcal {U }}(s))\eta(s)| +
|\tilde{\mathcal {U }}(s)(\eta(s) -\tilde\eta(s))| \cr
&\hskip -1cm
\le C e^{-\frac{s}{2}}|\eta(s)||\gamma(s) - \tilde{\gamma}(s)| + C
e^s|\eta(s) - \tilde{\eta}(s)| \cr
&\hskip -1cm
\le Ce^s|\gamma(s) - \tilde{\gamma}(s)| \le C e^{2s}d(\gamma, \tilde
\gamma),\cr
\end{array}
\label{eq:P2.68}
\end{equation}
\begin{equation}
e^{3s}|\mathcal {A }(s)- \tilde{\mathcal {A }}(s)|\le C
e^{\frac{5}{2}s}|\gamma(s) - \tilde{\gamma}(s)| \le C
e^{\frac{7}{2}s}d(\gamma,\tilde{\gamma})\,;
\label{eq:P2.69}
\end{equation}
here and in the remaining part of the proof $C$ is a generic constant that may
change value from line to line.

\noindent From these and similar estimates for the other terms appearing in
$\mathcal N$ we get
\begin{equation}
\left\{
\begin{array}{l}
|\mathcal N(\gamma(s))| \le Ce^{\frac{5}{2}s}\cr
|\mathcal N(\gamma(s)) -\tilde{\mathcal N}(\tilde{\gamma}(s))|  
\le Ce^{2 s}d(\gamma, \tilde \gamma) \cr
\end{array}
\right.\ .
\label{N_DeltaN}
\end{equation}
From (\ref{N_DeltaN}), (\ref{stima_eMs})
it follows
\begin{equation}
\label{eq:P2.72}
|(T\gamma)(s) - \gamma_\delta(s)| 
\le Ce^{\frac{5}{2}s}\,,
\end{equation}
\begin{equation}
\label{eq:P2.73}
|(T\gamma)(s) - (T\tilde\gamma)(s)| \le C d(\gamma, \tilde \gamma)
 \int_{-\infty}^{s} e^{2r} dr \le C e^{2s}d(\gamma, \tilde \gamma)\,,
\end{equation}
so that
\[
e^{-s}|(T\gamma)(s) - (T\tilde\gamma)(s)| \le Ce^{s_0} d(\gamma,
\tilde \gamma) \ \ \forall\; s\in(-\infty,s_0]\,,
\]
\[
\Longrightarrow \ \ d(T\gamma, T\tilde\gamma)\le Ce^{s_0}
d(\gamma, \tilde \gamma)\ .
\]
Therefore we see that, provided we take $-s_0>0$ sufficiently large, $T: X
\rightarrow X$ is a contraction. This concludes the proof of (I).

\medbreak\noindent To prove (II) we first observe that if $\overline
\gamma$ is the fixed point of $T$, the estimate (\ref{eq:P2.72}) implies
that
\begin{equation}
\label{eq:P2.75}
\lim_{s\rightarrow - \infty}|\overline \gamma(s) - \gamma_\delta(s)|
e^{-\frac{3}{2}s} = 0\ .
\end{equation}
We also observe that from (\ref{eq:P2.41}), (\ref{est_collect}) and $x =
e^s$ it follows
\begin{equation}
\left\{
\begin{array}{l}
t \propto \frac{2}{3} \sqrt{\frac{2}{\alpha}} x^{\frac{3}{2}} =
\frac{2}{3} \sqrt{\frac{2}{\alpha}} e^{\frac{3}{2}s},\cr
\frac{dt}{dx} \propto \sqrt{\frac{2}{\alpha}} x^{\frac{1}{2}} =
\sqrt{\frac{2}{\alpha}} e^{\frac{s}{2}}\cr
\end{array}
\right.\ .
\label{t_dtdx}
\end{equation}
From these asymptotic formulas and the definition of $\eta$ it follows
\begin{equation}
\eta = y' x = \dot y \frac{x}{\dot x} \propto \dot y
\sqrt{\frac{2}{\alpha}} e^{\frac{3}{2}s}\propto \dot y \frac{3}{2}t\ .
\label{eta_behavior}
\end{equation}
From (\ref{eq:P2.63}) and (\ref{t_dtdx})$_1$ we obtain
\begin{equation}
\gamma_\delta(t) \propto \frac{3}{2} \sqrt{\frac{\alpha}{2}}\delta\rho_4 t\ .
\label{gammab_behavior}
\end{equation}
Inserting (\ref{eta_behavior}), (\ref{gammab_behavior}) into
(\ref{eq:P2.75}) yields
\[
\left\{
\begin{array}{l}
\lim_{t \rightarrow 0^+} \frac{y(t) - \frac{3}{2}
\sqrt{\frac{\alpha}{2}} t \delta}{t} =0,\cr
\lim_{t \rightarrow 0^+} (\dot y(t) - \frac{3}{2}
\sqrt{\frac{\alpha}{2}} \delta) =0 \cr
\end{array}
\ .
\right.
\]
If we take $\delta = \frac{2}{3} \sqrt{\frac{2}{\alpha}}\olda$, then these
equations show that the unique solution determined by the fixed point
$\overline \gamma$ of $T$ satisfies (\ref{lim_pdot}). This concludes
the proof of Proposition~\ref{prop:unique}.

\rightline{\small$\square$}
\end{proof}

\noindent Proposition~\ref{prop:unique} is formulated for {\em ejection}
solutions. An analogous uniqueness result applies to {\em collision} solutions
and can be derived from Proposition~\ref{prop:unique} by the variable change
$t\rightarrow -t$.
\begin{corollary}
Let $w:(0, \bar{t}) \to \R^3$ ($(-\bar{t},0) \to \R^3$) be a maximal ejection
(collision) solution to (\ref{K_motion_general}) and let ${\mathsf n} =
\lim_{t\to 0^+} \frac{w(t)}{|w(t)|}$ ($= \lim_{t\to 0^-} \frac{w(t)}{|w(t)|}$)
be the unit vector in Proposition~\ref{prop:limcoll}.  Assume there exists an
axis $r'\neq r$ of some rotation in $G\setminus\{I\}$ such that 
${\mathsf n}$ is parallel to the plane $r r'$.  Then
\begin{equation}
w(t) \in {\rm span}\{\er,{\mathsf n}\},\ \forall t\in(0,\bar{t}) \ \; (\forall
t\in(-\bar{t},0))\ .
\label{eq:w_in_a_plane}
\end{equation}
\label{cor:w_in_a_plane}
\end{corollary}
\begin{proof}
From (\ref{symcond_Rtilde}) it follows that if $w$ is a solution to
(\ref{K_motion_general}), then $\tilde{R} w$ is also a solution to
(\ref{K_motion_general}) with the same values of $\olda$ and $h$. Then
Proposition~\ref{prop:limcoll} implies $\tilde{R}w = w$ and
(\ref{eq:w_in_a_plane}) follows.

\rightline{\small$\square$}
\end{proof}
The following consequence of Proposition~\ref{prop:unique} is not used in the
rest of the paper, but it is of independent interest:
\begin{corollary}
Let $w:(-\bar{t}^-, 0)\cup (0, \bar{t}^+)\rightarrow \R^3$ be a maximal {\it
ejection-collision} solution of
(\ref{K_motion_general}),(\ref{K_energy_general}). Assume that there exists a
unit vector ${\mathsf n}$ such that
\begin{itemize}
\item[(i)] $\lim_{t\rightarrow 0^+}\frac{w(t)}{|w(t)|} = \lim_{t\rightarrow
  0^-}\frac{w(t)}{|w(t)|} = {\mathsf n},$

\item[(ii)] $\lim_{t\rightarrow 0^\pm}(\dot{w}(t) + R\dot{w}(t)) = 0,$

\item[(iii)] the energy constant $h^-$ in the interval $(-\bar{t}^-,0)$
  coincides with the energy constant $h^+$ in the interval $(0,\bar{t}^+)$.
\end{itemize}
\noindent
Then $\bar{t}^- = \bar{t}^+ = \bar{t}$ and $w(t) = w(-t)$, $ t \in (- \bar{t},
\bar{t})$.
\label{cor:unique}
\end{corollary}
We are now in the position of proving
\begin{proposition}
Let $\minloop\in\overline{\K}$ be a minimizer of the action. Then
$\minloop\in\K$ and it is collision free.
\end{proposition}
\begin{proof}
By the discussion in Section~\ref{sec:tot_coll} $\minloop$ does not have total
collisions. From Proposition~\ref{prop:npiunmeno} and
Remark~\ref{rem:npiunmeno} $\minloop$ can not have partial collisions apart
from collisions of type $(\rightrightarrows)$.  From Lemma~\ref{lem:n_in_Pi}
and Corollary~\ref{cor:w_in_a_plane} it follows that
\[
\minloopgen(t)\in\Pi\,, \qquad \forall t \in (t_c-\bar{t}, t_c+\bar{t})\ .
\]
Since the only possible collisions for $\minloop$ are partial collisions of
type $(\rightrightarrows)$ at times $t=t_c-\bar{t}^-, t_c+\bar{t}^+$,
$\minloop$ has collisions of this kind. Therefore, invoking again
Corollary~\ref{cor:w_in_a_plane}, we obtain
\[
\minloopgen(t) \in\Pi\,, \qquad \forall t\in\R\ .
\]
This contradicts membership in $\overline{\K}$ and concludes the proof.

\rightline{\small$\square$}
\end{proof}

\begin{remark}
As we have seen there are infinitely many pairwise disjoint cones
$\K(u)\subset\lamoa$ corresponding to simple $\genp$. Moreover the
subset of such cones that satisfy the condition
\begin{equation}
S^{\tilde{k}} \neq S^k\,, \qquad k\in\Z
\label{starcond}
\end{equation}
is infinite. This is trivially seen by observing that if $u\in\lamoa$
is such that the sequence $\sigma_u$ satisfies condition
(\ref{starcond}).  it follows that all the cones $\K_n, n=1,..$
corresponding to the pairs $(\sigma_u,n), n=1,...$ satisfy condition
(\ref{starcond}).

It is easily seen that besides these rather trivial cases there are
many other situations where (\ref{starcond}) holds.

Are all these minimizers genuine periodic solutions of the classical Newtonian
$N$--body problem? Our conjecture is that when the linking of the trajectory
of the generating particle $\tau_1$ with $\Gamma$ becomes more and more
complex, from the point of view of minimization of the action, it may be
preferable to {\em crash} some of the complications into a total collision.
If this conjecture is correct then only a finite number of minimizers
corresponding to simple linking of $\tau_1$ with $\Gamma$ are genuine
periodic solutions of the $N$-body problem. 
\label{rem:finitecones}
\end{remark}




\newcommand{\TT}{{\mathcal{T}}}
\newcommand{\T}{{\mathcal{T}}}

\section{Complements, conjectures and numerical results}
\label{sec:conj}

Theorem~\ref{teo:collfreemin} yields valuable information on the geometric
structure of the periodic motion $\minloop\in\K$ considered in
Theorem~\ref{platorb_exist}. In particular we know that the orbit $\tau_1$ of
the generating particle is contained in the cone $C_\K\subset\R^3$, $C_\K =
\cup_{k\in\Z} S_k\cup D_k$ with $\{D_k\}_{k\in\Z}$ the sequence corresponding
to $\K$. Another consequence of Theorem~\ref{teo:collfreemin} is the existence
of motions that violate the coercivity condition (\ref{zeromean}) (cfr.
Theorem~\ref{coerteor}). We have indeed
\begin{theorem}
Assume $P\in\{\Cplato, \Iplato, \Dplato\}$. Then there exists a smooth
$T$--periodic solution $\motionloopPi\in\KPi$ of the classical
$N$--body problem such that the trajectory $\tau_1$ of the generating
particle is contained in the half space $\asseuno>0$. In particular we
have
\[
\int_0^T \motionloopPigen(t)\;dt >0\ .
\]
\label{teo:not_meanzero}
\end{theorem}
\begin{proof}
Let $\nu$ be the minimal sequence associated to $\KPi$ in
Section~\ref{sec:coll} (cfr. Table~\ref{tab:KPi_seq}).  Let $\motionloopPi$ be
the minimizer of the action restricted to the cone $\K(\vnuuno)\subset
\KPi$. All the arguments developed in Section~\ref{sec:coll} to prove that a
minimizer $\minloopPi$ of $\A|_{\KPi}$ is free of collisions apply verbatim to
show that also $\motionloopPi$ is free of collisions. Therefore
$\motionloopPi$ is a classical solution of the $N$--body problem.  To conclude
the proof we observe that if $P\in\{\Cplato, \Iplato, \Dplato\}$ and
$u=\motionloopPi$, the cone $C_{\K(u)}\subset\R^3, C_{\K(u)} = \cup_{k\in\Z}
S_k\cup D_k$ with $\{D_k\}_{k\in\Z}$ the sequence corresponding to $\K(u)$, is
actually contained in the half space $\asseuno>0$

\rightline{\small$\square$}
\end{proof}
The motion $\minloopPi$ is a minimizer of $\A$ on $\KPi$, which
properly contains $\K(\vnuuno)$; therefore $\minloopPi$ and
$\motionloopPi$ may well be different solutions of the $N$--body
problem in $\KPi$. However we conjecture that $\motionloopPi =
\minloopPi$, and this conjecture is supported by numerical
experiments.

Given a pair $(\sigma,n)$, with $\sigma$ satisfying (I), (II), (III)
in Section~\ref{subsec:chartop} and $n\in\N$, the associated $\usngen$
defined in (\ref{genpsigman}) possesses all the space and time
symmetries compatible with its topological structure. We can minimize
the action on $\K(\usn)$ or on the subset $\K^S(\usn)$ of the maps
$\typloop$ with the property that $\genp$ has the same symmetries as
$\usngen$. All numerical experiments we have made confirm the
conjecture that a global minimizer $\minloop\in\K(\usn)$ is actually
in $\K^S(\usn)$.

Assume $\minloop\in\K(\usuno)$ is a minimizer of $\A|_{\K(\usuno)}$.
Then the map $\minloop^n\in\K(\usn)$ defined by
\[
\minloop^n(t) := n^{-2/3} \minloop(nt), \ t\in[0,T)
\]
is a critical point of $\A|_{\K(\usn)}$. This follows from the fact
that the map $h_n: \K(\usuno) \to \K(\usn)$ defined by $(h_nu)(t) =
u^n(t) = n^{-2/3} u(nt)$ satisfies
\begin{equation}
\A(u^n) = n^{2/3}\A(u)\ .
\label{action_un}
\end{equation}
However we do not expect $\minloop^n$ to be a minimizer of
$\A|_{\K(\usn)}$ for every $n>1$.  Indeed we have
\begin{proposition}
Given $u\in\K(\usuno)$, if $n>n_0$, for some $n_0>1$, 
there exists $\widehat{u^n}\in\K(\usn)$ such that
\[
\A(u^n) - \A(\widehat{u^n}) > 0\ .
\]
\end{proposition}
\begin{proof}
We consider the case of $n$ even. The other case is analogous. Fix a number
$\tau\in(0,1)$ and set
\begin{eqnarray*}
T_1 &=& \Bigl(\frac{1}{2} - \frac{1}{n}\Bigr)(1+\tau)T\,,\\
T_2 &=& T_1 + \frac{T}{n}\,,\\
T_3 &=& T_2 + \Bigl(\frac{1}{2} - \frac{1}{n}\Bigr)(1-\tau)T \ .
\end{eqnarray*}
Define $\widehat{u^n}$ by setting
\begin{equation}
\widehat{\genp^n}(t) = 
\left\{
\begin{array}{ll}
\displaystyle 
\Bigl(\frac{1}{1+\tau}\Bigr)^{-2/3} \genp^n \Bigl(\frac{t}{1+\tau}\Bigr)\,, 
&0\leq t\leq T_1
\cr
&\cr
\displaystyle 
n\left[\frac{T_2-t}{T} \Bigl(\frac{1}{1+\tau}\Bigr)^{-2/3} +
\frac{t-T_1}{T} \Bigl(\frac{1}{1-\tau}\Bigr)^{-2/3}\right] \genp^n(t-T_1)\,,
&T_1\leq t\leq T_2
\cr&\cr
\displaystyle 
\Bigl(\frac{1}{1-\tau}\Bigr)^{-2/3} \genp^n
\Bigl(\frac{t-T_2}{1-\tau}\Bigr)\,,
&T_2\leq t\leq T_3
\cr&\cr
\displaystyle 
n\left[\frac{T-t}{T} \Bigl(\frac{1}{1-\tau}\Bigr)^{-2/3} +
\frac{t-T_3}{T} \Bigl(\frac{1}{1+\tau}\Bigr)^{-2/3}\right] \genp^n(t-T_3)\,,
&T_3\leq t\leq T\cr
\end{array}
\right.
\ .
\end{equation}
Then we have, with obvious notation, recalling also (\ref{action_un}),
\begin{equation}
\left\{
\begin{array}{l}
\A(\widehat{u^n}, (0,T_1)) = 
\displaystyle \frac{n-2}{2} \Bigl(\frac{1+\tau}{n}\Bigr)^{1/3}\A(u)
= \displaystyle \frac{n-2}{2n}(1+\tau)^{1/3}\A(u^n) \leq\cr
\hskip 2.3cm
\leq\displaystyle \Bigl[\frac{1}{2}(1+\tau)^{1/3} +\frac{C}{n}\Bigr]\A(u^n)\cr
\cr
\A(\widehat{u^n}, (T_2,T_3)) = 
\displaystyle \frac{n-2}{2} \Bigl(\frac{1-\tau}{n}\Bigr)^{1/3}\A(u)
= \displaystyle \frac{n-2}{2n}(1-\tau)^{1/3}\A(u^n) \leq\cr
\hskip 2.3cm\leq\displaystyle
\Bigl[\frac{1}{2}(1-\tau)^{1/3} +\frac{C}{n} \Bigr]\A(u^n)\cr
\cr
\A(\widehat{u^n}, (T_1,T_2)\cup(T_3,T)) \leq \displaystyle 
\frac{C}{n}\A(u^n)  \cr
\end{array}
\right.
\end{equation}
where $C>0$ is a constant that does not depend on $n$. It follows
\[
\A(u^n) - \A(\widehat{u^n}) \geq \left[\frac{1}{2}
(2 - (1+\tau)^{1/3} - (1-\tau)^{1/3}) - 3\frac{C}{n}\right]\A(u^n)\ .
\]
This inequality and $2>(1+\tau)^{1/3} + (1-\tau)^{1/3}$ conclude the
proof.

\rightline{\small$\square$}
\end{proof}

By analogy with the fact that, in the Kepler problem, if the index of the
orbit with respect to the center of attraction is different from $-1,0,1$, the
only minimizer is the collision--ejection motion \cite{gordon77}, we don't
expect that, if $\typloop\in\lamoa$ is not simple in the sense of
Definition~\ref{def:simple}, then a minimizer $\minloop\in\K(u)$ is collision
free. In Section~\ref{sec:coll}, Remark~\ref{rem:finitecones}, we have
conjectured that, in spite of the fact that there are infinitely many pairwise
disjoint cones $\K(u)\subset\lamoa$ corresponding to a simple $\typloop$, only
for a finite number of them the corresponding minimizers are indeed genuine
$T$--periodic solutions of the Newtonian $N$--body problem.
To put this conjecture in perspective we consider the $N$-body problem with
generalized potential $\frac{1}{r^\alpha}$, $\alpha>0$, and the corresponding
action functional
\begin{equation}
\label{eq:con3}
\A^\alpha(u) = \frac{N}{2}\int_0^T \Big(|\dot \genp(t)|^2 + \sum_{R\in
\mathcal R\setminus \{I\}}\frac{1}{|(R-I)\genp|^\alpha}\Big)\ .
\end{equation}
It is well known that on the basis of Sundman's estimates in the case of
strong forces ($\alpha\ge 2$) $\A^\alpha(u)< +\infty$ implies that $u$ is
collision free. Therefore, if $\mathcal{N}^\alpha$ is the number of the cones
$\K \in\lamoasim$ such that a minimizer $\minloop^\alpha
\in \K$ of $\A^\alpha|_\K$ is collision free, we have $\mathcal{N}^\alpha =
+\infty$ for $\alpha \in [2, +\infty)$. We conjecture that, on
the other hand, $\mathcal{N}^\alpha < +\infty$ for $\alpha \in
(0,2)$ and that $\mathcal{N}^\alpha$ is a non-decreasing function of
$\alpha$. It is also natural to conjecture that for each $(\sigma,n)$
there is a critical value $\alpha^{(\sigma,n)} > 0$ such that a
minimizer of $\A^\alpha|_{\K(\usn)}$ is collision free for $\alpha >
\alpha^{(\sigma,n)}$, while it has collisions for $\alpha\in
(0,\alpha^{(\sigma,n)}]$.

\noindent
Natural questions to ask are: what are the transformations a
minimizer $\minloop^\alpha \in \K(\usn)$ undergoes when $\alpha$
crosses the critical value $\alpha_{(\sigma,n)}$ and decreases to
$0^+$? What are the asymptotic properties of a minimizer
$\minloop^\alpha \in \K(\usn)$ when $\alpha \to +\infty$?

\begin{figure}
\centerline{\epsfig{figure=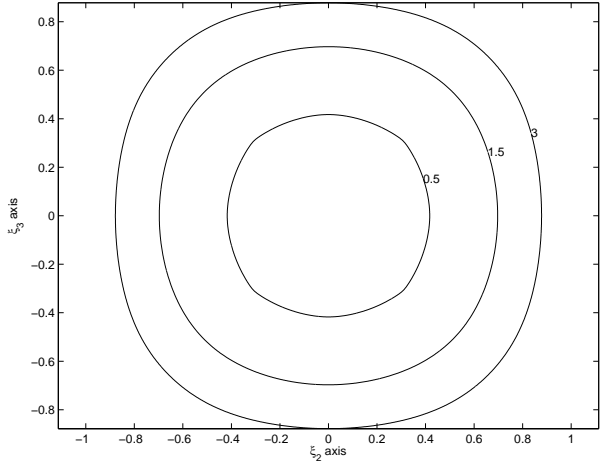,width=8cm}}
\caption{Shrinking of the trajectory of the periodic orbits found for the cone
$\KPdue$, with $P=\Cplato$, as $\alpha$ is decreasing. In the figure we use
$\alpha=3, 1.5, 0.5$}
\label{shrink_traject}
\end{figure}

When $\alpha \rightarrow 0^+$ the radius of action of the attractive forces
between particles converges to zero.  This implies that, at each time $t$, the
generating particle $\partgen$ has to be near some of the other $N-1$
particles. Indeed, otherwise, $\partgen$ could not accelerate to get around
the axes of the rotations in $\mathcal R$ as required by the topological
constraint of membership in $\K(\usn)$. This leads to conjecture that the
diameter of the orbit of the generating particle $\partgen$ converges to zero
as $\alpha \rightarrow 0^+$.  Then it is also to be expected that for
travelling a shorter and shorter trajectory in a fixed time $T$ also the
average speed of $\partgen$ tends to zero:
\begin{equation}
\label{eq:con5}
\left\{
\begin{array}{l}
\lim_{\alpha\rightarrow 0^+} || \minloopgen^\alpha ||_{L^\infty} = 0\\
\lim_{\alpha\rightarrow 0^+} || \dminloopgen^\alpha||_{L^1} = 0
\end{array}
\right.\ .
\end{equation}
In Figure~\ref{shrink_traject} we show the trajectory
$\tau_1^\alpha$ of $\partgen$ for various values of $\alpha$. The
shrinking of $\tau_1^\alpha$ for decreasing $\alpha$ is clearly
visible.

\noindent
If $\alpha>> 1$ the force of attraction between particles is very small when
the interparticle distance $|x_i - x_j|$ is larger than 1 and very large when
it is smaller than 1 (recall that the gravitational constant is
normalized to 1). This observation suggests:

\noindent
1) the limit behavior of minimizers for $\alpha \rightarrow +\infty$ is
   constrained to the sub-region $\mathcal{Y} \subset \mathcal{X}$ of the
   configuration space defined by
\[
\mathcal{Y}= \{x \in \mathcal X: |x_i - x_j| \ge 1, \forall i\ne j\}\ .
\]
Indeed violating this condition generates a large contribution of the
potential term to the action integral;

\noindent
2) in the limit $\alpha \rightarrow +\infty$ the trajectory $\tau_1^\alpha$ of
   $\partgen$, corresponding to a minimizer $\minloop^\alpha\in\K(\usn)$, is
   the shortest possible compatible with the condition $|(R - I)x_1| \ge 1$,
   $\forall R\in \mathcal R\setminus\{I\}$ and with the topological constraint
   that $\minloopgen^\alpha$ is homotopic to $\usngen$.

The rationale behind this is that, in the limit $\alpha\to +\infty$, the
interparticle attraction should act as a `perfect holonomic constraint' and
the limit motion should be a kind of `geodesic motion' with constant kinetic
energy. Then minimizing the action should be equivalent to the minimization of
the length of the trajectory of the generating particle. In conclusion we
advance the following conjecture. To formulate it we first observe that the
characterization of the set $\mathcal{Y}\subset \mathcal{X}$ is equivalent to
the characterization of the set $\mathcal{Y}_1\subset
\R^3$, $\mathcal{Y}_1 = \{x_1: |Rx_1 - x_1|\geq 1, \forall R\in{\cal
R}\setminus\{ I\}\}$. Let ${\cal P} =
\Gamma\cap\{x_1\in\R^3: |x_1|=1\}$, the set of poles of ${\cal R}$.
For each $p\in{\cal P}$ let $Cyl_p\subset\R^3$ be the open cylinder with axis
the line through $O$ and $P$ and radius $r_p= \frac{1}{2\sin(\pi/|C_p|)}$,
where $C_p\subset{\cal R}$ is the maximal cyclic group corresponding to
$p$. Then we have
\[
{\cal Y}_1 = \R^3\setminus \bigcup_{p\in{\cal P}}Cyl_p\ .
\]

\begin{conjecture} 
Let $\minloop^\alpha\in\K(\usn)$ be a minimizer of
$\A^\alpha|_{\K(\usn)}$. Then, in the limit $\alpha\to+\infty$,
$\minloopgen^\alpha$ converges (possibly up to subsequences) in the $C^1$
topology to a minimizer of the problem
\[
\min_{\genp\in \mathcal{U}} \frac{1}{2} \int_0^T |\dgenp (t)|^2\;dt\,,
\]
where $\mathcal{U}$ is the subset of $H^1$ $T$--periodic maps
$\genp:\R \to \mathcal{Y}_1$ such that $\genp$ is homotopic
to $\usngen$.
\label{conjecture}
\end{conjecture}

In other words, in the limit $\alpha\to +\infty$, the speed of the
generating particle is constant and the trajectory coincides with the
path of a wire stretched on the boundary of $\mathcal{Y}_1$ and
homotopic to $\usngen$. Since $\mathcal{Y}_1$ is a unilateral
constraint the limit trajectory is expected to be the union of
`geodesic arcs' lying on the $\partial \mathcal{Y}_1$ and segments in
$\R^3\setminus \mathcal{Y}_1$ touching $\partial \mathcal{Y}_1$ at the
extrema.  If we apply Conjecture~\ref{conjecture} to the case of four
bodies considered in Section~\ref{sec:cones} we get that the limit
trajectory of the generating particle is the union of four circular
arcs of angle $\pi$ and radius $1/2$ and that the constant limit speed
$|\dminloopgen^{\infty}|$ and corresponding action $\A^\infty$ are
given, for $T=1$, by
\begin{equation}
\left\{
\begin{array}{l}
|\dminloopgen^{\infty}| = 2\pi\cr
\A^\infty = 8\pi\cr
\end{array}\right.\ .
\label{action_infty}
\end{equation}

\begin{figure}[ht]
\psfragscanon
\psfrag{x}{$\xi_1$}
\psfrag{y}{$\xi_2$}
\psfrag{z}{$\xi_3$}
\centerline{
\epsfig{figure=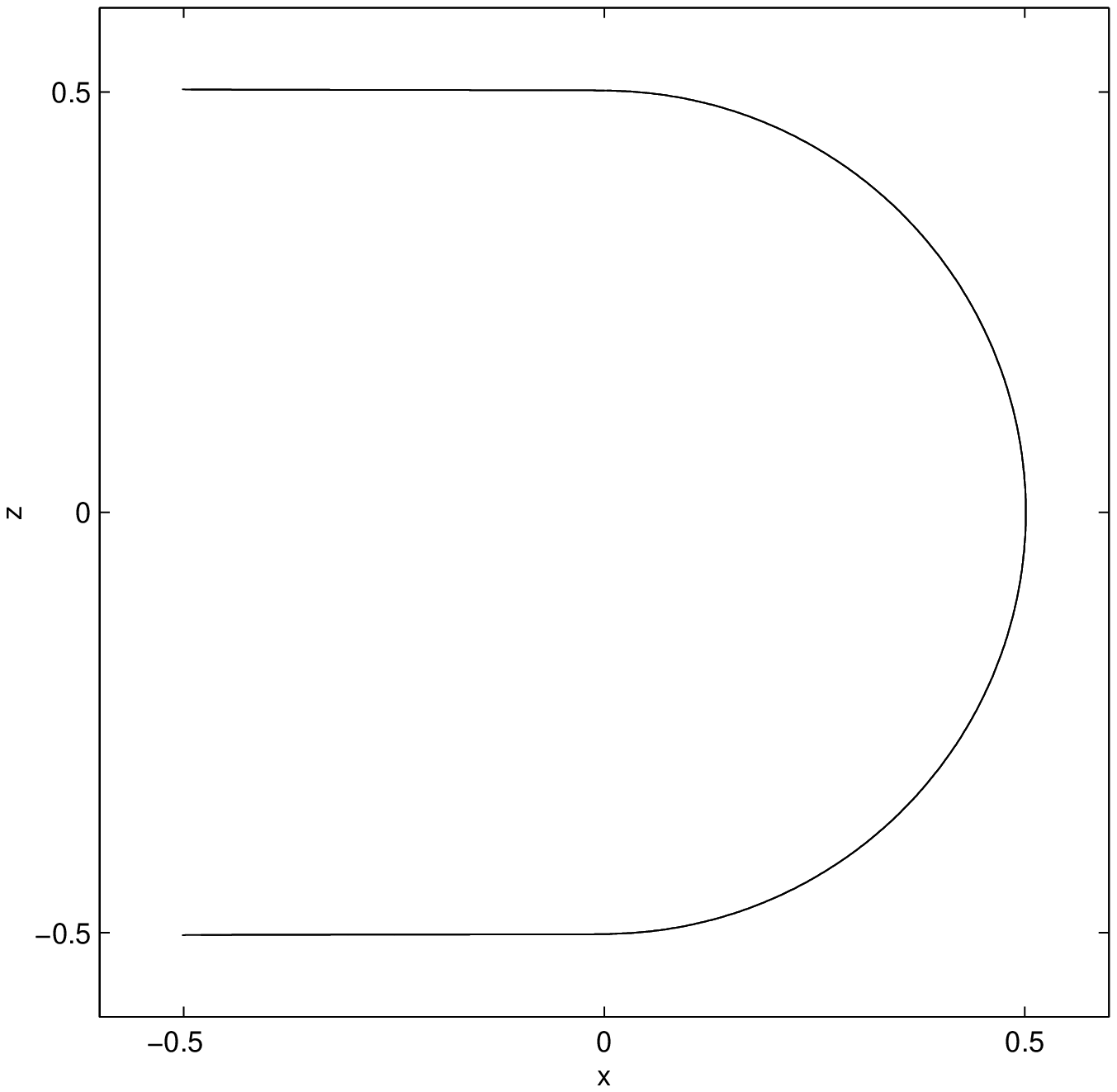,width=4.2cm}
\hskip 0.3cm
\epsfig{figure=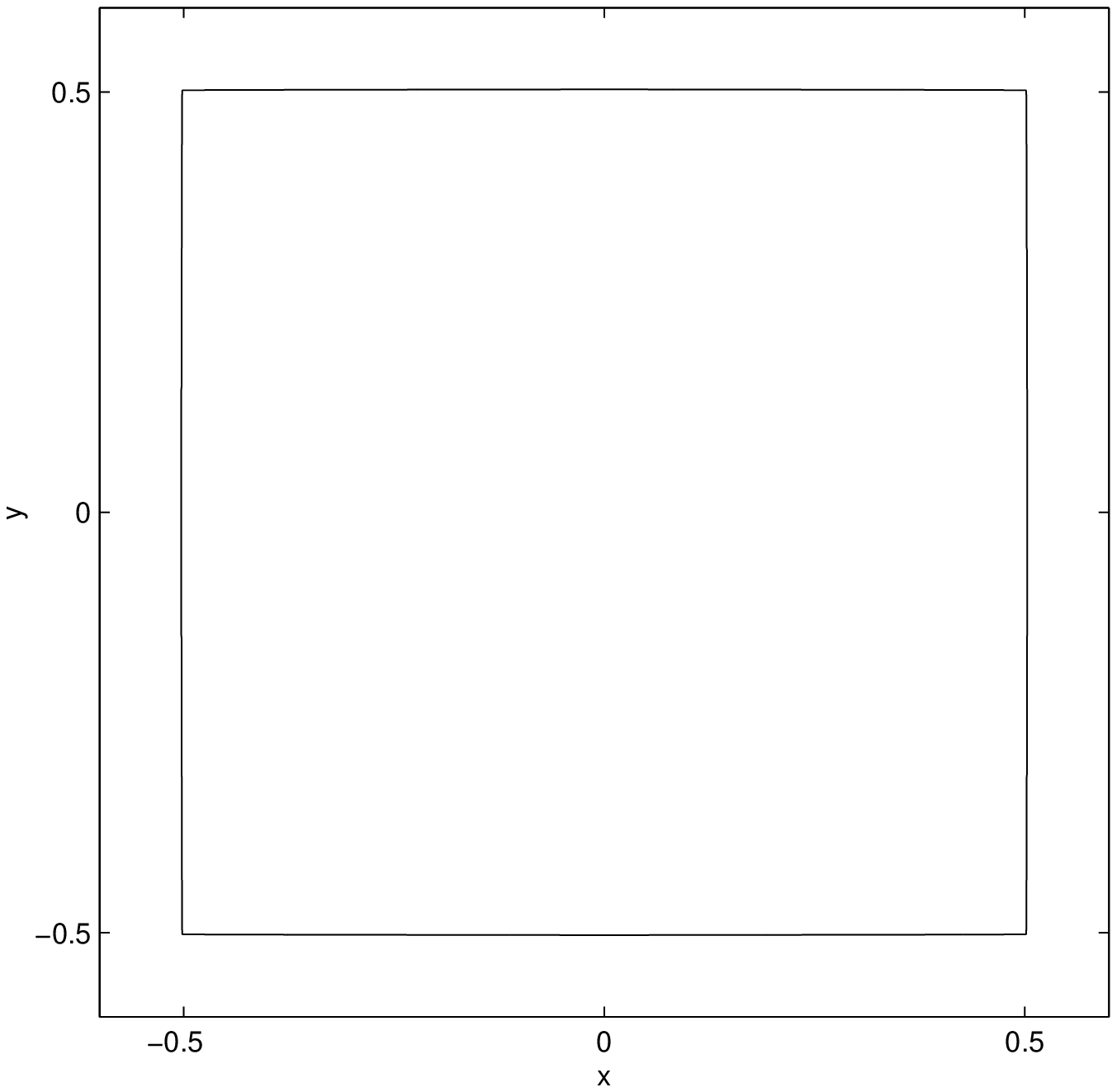,width=4.2cm}}
\centerline{\epsfig{figure=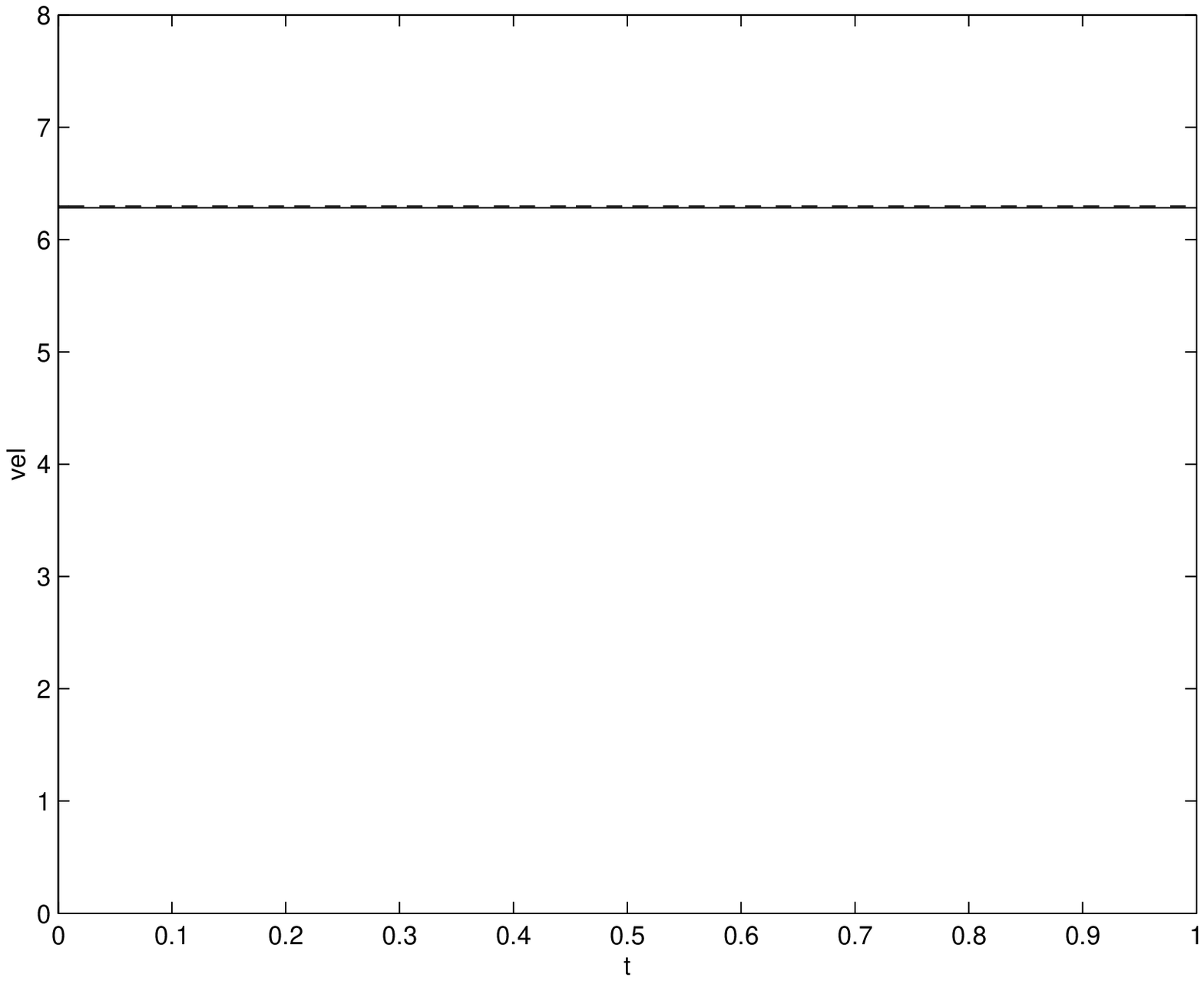,width=5cm}
}
\psfragscanoff
\caption{Projections of $\tau_1^{1000}$ and graph of 
$t\to |\dminloopgen^{1000}(t)|$. The dashed line corresponds to the
conjectured value for $|\dminloopgen^{\infty}(t)|$.}
\label{fig:check_vel}
\end{figure}
In Figure~\ref{fig:check_vel} we show orthogonal projections of
$\tau_1^\alpha$ and the graph of the map $t\to
|\dminloopgen^\alpha(t)|$ computed numerically for $\alpha = 1000$.
The value of $|\dminloopgen^{1000}(t)|$ is practically coincident with
the limit value $2\pi$ given by (\ref{action_infty}). Also the
numerical value of $\A^{1000}$ is practically equal to $8\pi$.

\begin{figure}[ht]
\psfragscanon
\psfrag{a1,1}{$\asseuno$}\psfrag{a1,2}{$\assedue$}
\psfrag{a2,1}{$\asseuno$}\psfrag{a2,2}{}
\psfrag{us0}{$\minloopgen^\infty(0)$}
\psfrag{usT8}{$\minloopgen^\infty(\frac{T}{8})$}
\psfrag{1su2}{$\frac{1}{2}$}
\psfrag{1susq2}{$\frac{1}{\sqrt{2}}$}
\psfrag{1susq3}{$\frac{1}{\sqrt{3}}$}
\psfrag{sq2}{$\sqrt{2}$}
\psfrag{M}{$M$}\psfrag{V}{$V$}
\centerline{\epsfig{figure=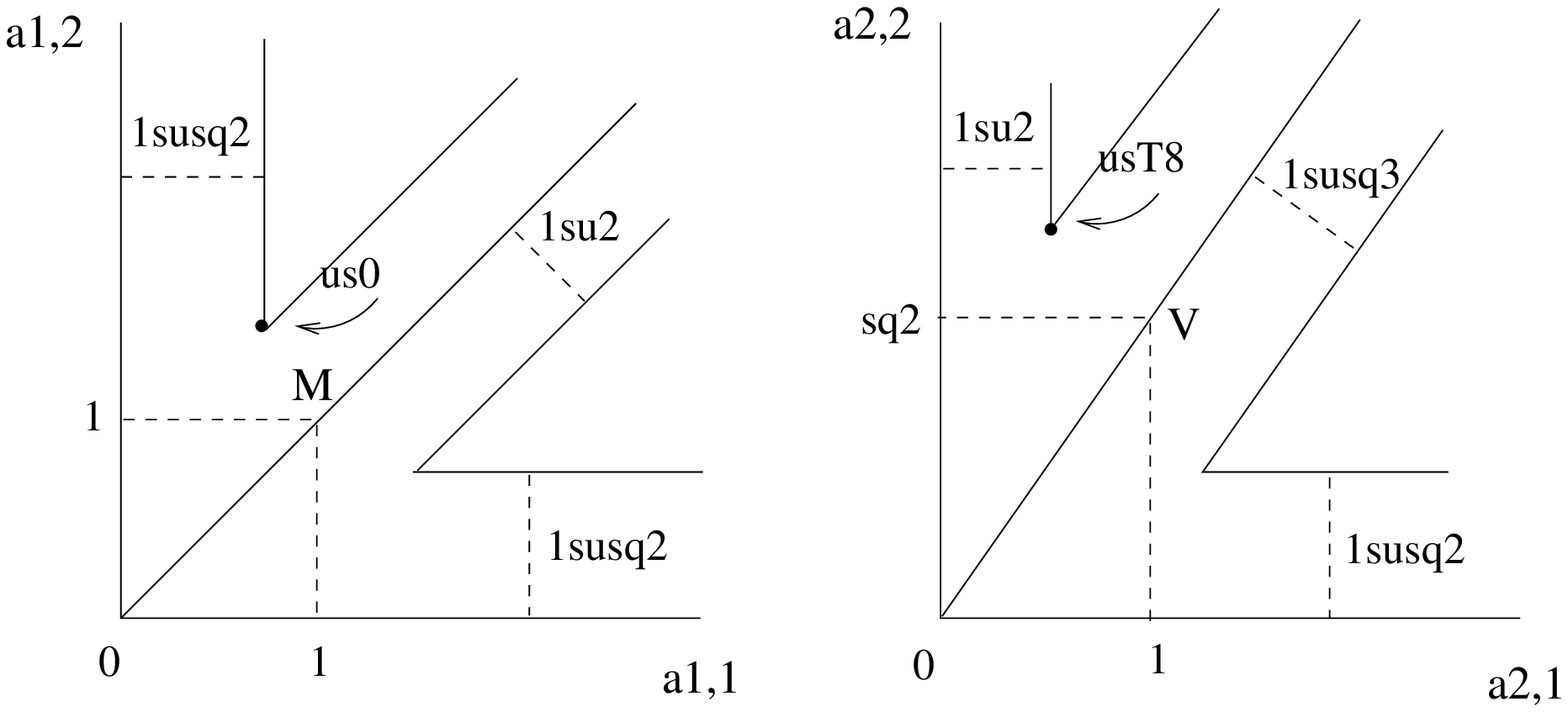,width=13cm}}
\psfragscanoff
\caption{The geometry of ${\cal Y}_1$ and $\minloopgen^\infty(0),
\minloopgen^\infty(T/8)$ for $\KPdue$, $P=\Cplato$. 
We plot the section of ${\cal Y}_1$ with the plane $\assetre=0$ (left) and the
one with the plane determined by $\asseuno$ and $V$ (right).}
\label{missingfig1}
\end{figure}

We found good agreement with this conjecture in several other situations. For
instance, for the cone $\KPdue$ defined by (\ref{cone2})
with $P=\Cplato$, a simple computation (cfr. Figure~\ref{missingfig1}) based on
Conjecture~\ref{conjecture} shows that the limit motion of the generating
particle should satisfy
\begin{equation}
\left\{\begin{array}{l}
{\minloopgen}^\infty(0) = (\frac{1}{\sqrt{2}},\sqrt{2},0)\cr
{\minloopgen}^\infty(\frac{T}{8}) = (
\frac{1}{2},
\frac{1}{\sqrt{2}}(1+\frac{1}{\sqrt{2}}),
\frac{1}{\sqrt{2}}(1+\frac{1}{\sqrt{2}}))\cr
\end{array}
\right.
\label{minloopinfty}
\end{equation}

\begin{figure}[ht]
\psfragscanon
\psfrag{alpha}{$\alpha$}
\psfrag{x}{\framebox[0.6cm]{$\xi_1$}}
\psfrag{y}{\framebox[0.6cm]{$\xi_2$}}
\epsfig{figure=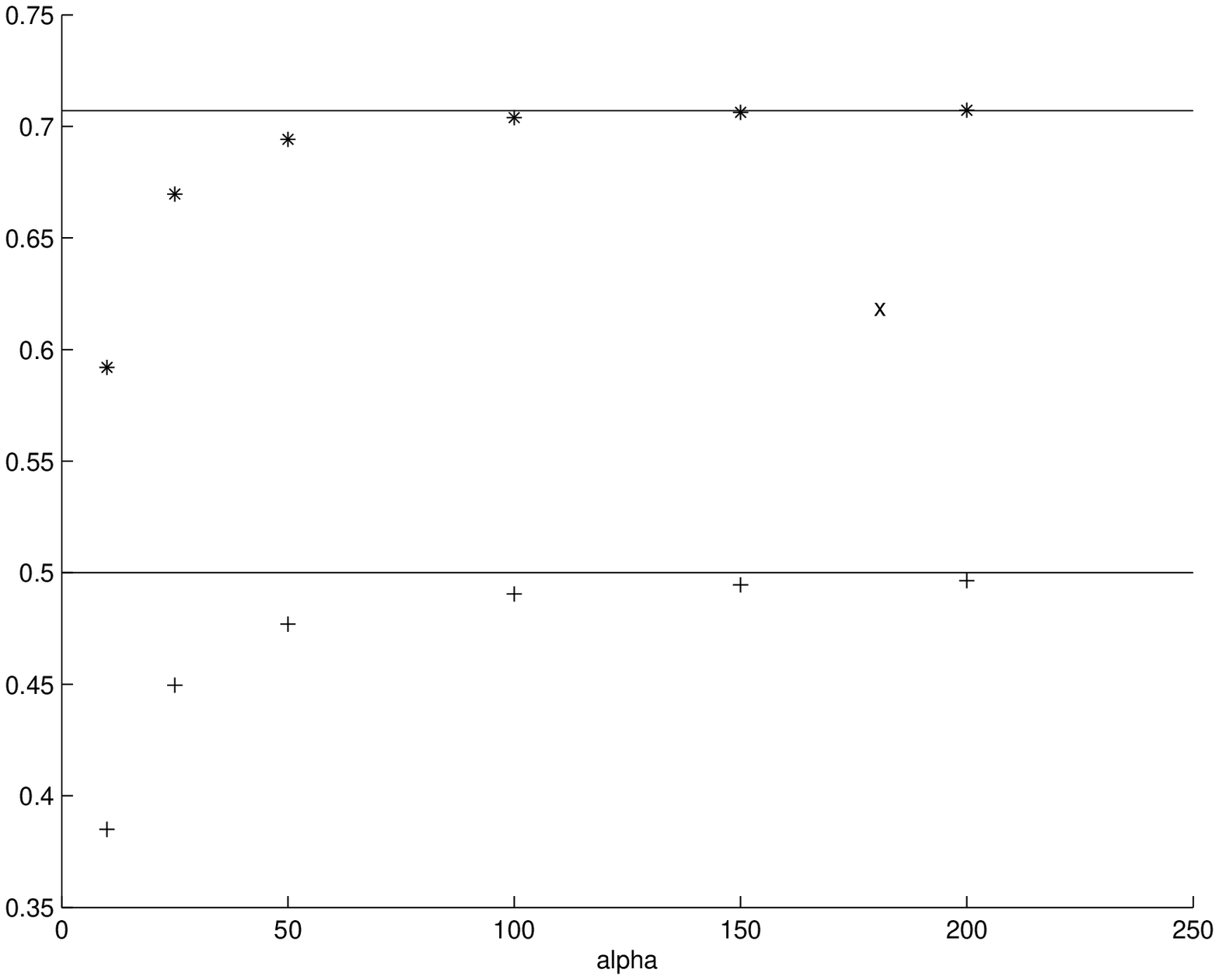,width=5.7cm}
\epsfig{figure=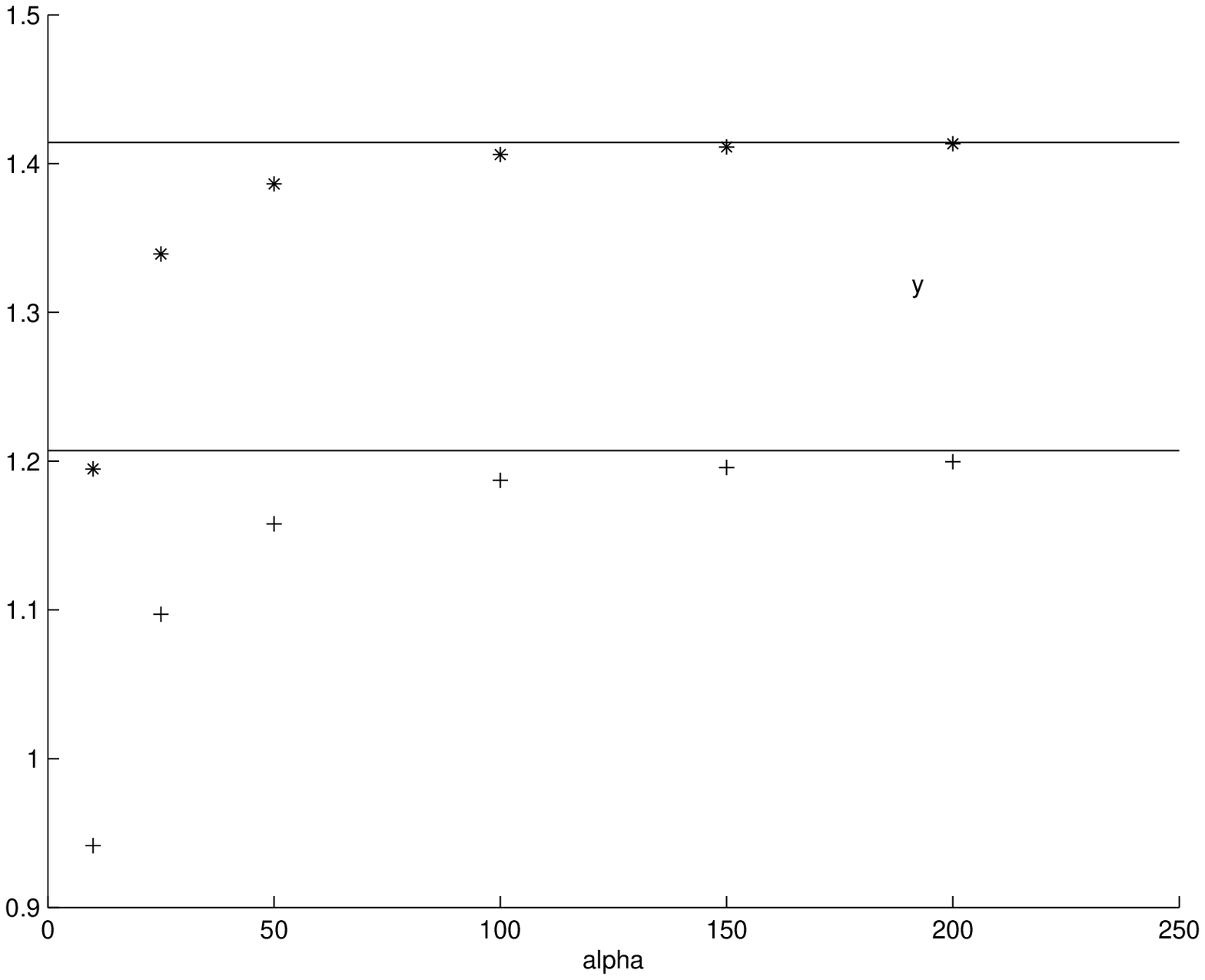,width=5.7cm}
\psfragscanoff
\caption{
Numerical simulations of $\minloopgen^\alpha(0)\cdot\ej$ and
$\minloopgen^\alpha(\frac{T}{8})\cdot\ej$ for $j=1,2$ and $\alpha=
10,25,50,100,150,200$.  The asterisks correspond to $\minloopgen^\alpha(0)$,
the crosses to $\minloopgen^\alpha(\frac{T}{8})$.
The horizontal lines indicate the values in (\ref{minloopinfty}).}
\label{missingfig2}
\end{figure}

In Figure~\ref{missingfig2} we plot as functions of $\alpha$ the
components of $\minloopgen^\alpha(0)$ and
$\minloopgen^\alpha(T/8)$ and the theoretical values given by
(\ref{minloopinfty}). We also plot the speed
$|\dminloopgen^\alpha(t)|$ for several values of $\alpha>>1$. For
$\alpha=200$ the relative error is within $1/100$.

The class of cones $\K$ such that (\ref{bordercone}) holds can be
largely generalized by requiring the sequence $\sigma=\{D_k
\}_{k\in\Z}$ introduced in Subsection~\ref{subsec:chartop} to satisfy only the
conditions (II), (III).  As before we can associate to $\sigma$ a
motion ${\rm u}^\sigma:\R\to\R^{3N}$ 
as in Section~\ref{subsec:chartop}
and consider the cone $\K^\sigma\subset H^1(\R,\R^{3N})$ of the maps
homotopic to $u^\sigma$.  The cones $\K^\sigma$ satisfy
(\ref{bordercone}) and, even if we restrict to sequences $\sigma$ that are
simple in the sense of Definition~\ref{def:simple}, form an
uncountable family.  Now we can not talk about a minimizer
$\minloop\in\K^\sigma$ since the action $\A(u)$ of each
$u\in\K^\sigma$ is infinite. But we can still regard as a minimizer a
function $\minloop\in\K^\sigma$ with the property that the restriction
$\minloop|_{[t_1,t_2]}$ to each compact interval $[t_1,t_2], t_1<t_2$
minimizes the action on the set of maps $\typloop:[t_1,t_2]\to\R^{3N}$
that satisfy $\typloop(t_i) = \minloop(t_i), i=1,2$ and are such that
$\chi_{\R\setminus[t_1,t_2]} \minloop + \chi_{[t_1,t_2]} \typloop$ is
homotopic to $u^\sigma$. For potentials corresponding to $\alpha\geq
2$, such a minimizer $\minloop$ can not have collisions and, depending
on the particular sequence $\sigma$ considered, it may exhibit
several interesting and complex behaviors, including heteroclinic
connections between periodic orbits and chaotic motion. For the
Newtonian potential ($\alpha=1$), at least under the assumption that
the sequence $\sigma$ satisfies condition (\ref{starcond}), the
discussion in Section~\ref{sec:coll} excludes the possibility of
partial collisions. We conjecture that in many cases, actually
infinitely many, the minimizers corresponding to such sequences $\sigma$ do
not have total collisions too and therefore are genuine solutions of
the classical $N$--body problem. 

Returning to the setting of periodic motions, we note that a natural
generalization of the situation discussed in Section~\ref{sec:cones} where we
have just a single generating particle that determines the motion of all the
other particles, is obtained by considering $M\geq 1$ generating particles.

\smallbreak Given $M$ positive constants $m_h, h=1,..,M$ and $M$ cones
$\K_h\subset \lamoa, h=1,...,M$ we can consider maps $U=(u^1,...,u^M)
\in (\lamoa)^M$ and look for minimizers of the action
\begin{equation}
\A(U) = \frac{1}{2} \int_0^T \Biggl( N\sum_{h=1}^M |\dgenp^h(t)|^2 +
\sum_{\stackrel{R,R'\in{\cal R}}{(R',h)\neq(R',k)}} \frac{m_h
m_k}{|R_i\genp^h-\genp^k|}
\Biggr)\,dt
\label{M_action}
\end{equation}
on the cone $\K = \K(u^1)\times\K(u^2)\times\ldots
\K(u^M)\subset(\lamoa)^M$. In (\ref{M_action}) $\genp^h:\R\to\R^3$ is
the motion of the $h$--th generating particle. The analogous of
Proposition~\ref{prop:AKu_coerc} holds: $\A(U)|_\K$ is coercive.

\medbreak
We conclude this Section by a few remarks on the method used for
numerical simulation of the orbits discussed in the present paper. The
method is based on a rather simple and natural idea already used in
(\cite {moore93}): we consider the $L^2$ gradient dynamics defined by
the action functional
\begin{equation}
u_\theta = - grad_{L^2} \A(u)
\label{eq:con8}
\end{equation}
with periodic boundary conditions. Here $\theta$ is a fictitious time
and the true time $t$ plays the role of a space variable. Suppose that
$\K\in\lamoasim$ is such that if $\minloop \in \K$ is a
minimizer then
\begin{equation}
\label{eq:con9}
\A(\minloop) < a_c:= \inf\{ \A(u), \, u \in \Kclo,\,
u \in \mathfrak{S}\}
\end{equation}
It follows that $\minloop$ is collision free and that, if $\overline u
\in \K$ satisfies the condition $\A(\overline u) < a_c$ then
the dynamics defined by (\ref{eq:con8}), with initial datum $\overline u$,
automatically preserves the cone $\K$:
\begin{equation}
\label{eq:con10}
u(\cdot ,\theta, \overline u) \in \K,  \, \forall \theta \ge 0
\end{equation}
$\theta \rightarrow u(\cdot ,\theta, \overline u)$ being the solution
of (\ref{eq:con8}) through $\overline u$, and remains away from
$\mathfrak{S}$. This follows from $\A(u(\cdot ,\theta, \overline u))
\le \A(\overline u) < a_c$ and from (\ref{eq:con9}) which imply
$u(\cdot ,\theta, \overline u)$ can not cross the boundary of $\K$.

Once (\ref{eq:con10}) is established, the general theory of infinite
dimensional dissipative dynamical systems (see Lemma 3.8.2 in \cite{hale})
implies that the $\omega$--limit set of $\overline u$ is contained in the set
$E$ of equilibria of (\ref{eq:con8}), $E:= \{ u\in
\K: {\rm grad}_{L^2} \A(u) = 0\}$, that is, in the set of
periodic solutions of the classical Lagrange's equations which are the objects
of our interest.

Based on (\ref{eq:con8}) one can set up an automatic procedure for the
numerical computation of a minimizer $\minloop\in\K$ for each given
$\K\in\lamoasim$. To do this it suffices to develop a systematic way
for constructing a suitable initial condition
$\bar{\typloop}\in\K$. Our choice is
\[
\bar{\typloop} = \vnun
\]
where $(\nu,n)$ is the pair corresponding to $\K$ in the sense of
Proposition~\ref{prop_nu} and $\vnun$ is defined in (\ref{vnun}). We have
developed a routine that, given a sequence $\nu$ of
vertexes of an Archimedean polyhedron $\mathcal{Q}_{\cal R}$ and a number
$n\in\N$ computes a minimizer $\minloop$ of the action on the cone
$\K$ corresponding to $(\nu,n)$. A sample of motions computed in this way can
be found at {\tt http://adams.dm.unipi.it/\~{}gronchi/nbody/}.



\begin{acknowledgement}
It is a pleasure to express our gratitude to our collegues N. Gavioli
and N. Guglielmi who generously helped us with the numerical and
algebraic aspects of the present work.
\end{acknowledgement}



\end{document}